\documentclass[twoside,11pt]{article}

\usepackage[abbrvbib]{jmlr2e}

\usepackage{array}
\usepackage[normalem]{ulem}

\usepackage{url}            
\usepackage{amsfonts}       
\usepackage{microtype}      
\usepackage{multirow}
\usepackage{rotating}
\usepackage{lastpage}
\usepackage{tabularx}
\usepackage{amsmath}
\usepackage{bm}
\usepackage{enumitem}
\usepackage{xcolor}
\usepackage{algorithm}
\usepackage{algorithmicx,algpseudocode}

\makeatletter
\algnewcommand{\LineComment}[1]{\Statex \hskip\ALG@thistlm \(\triangleright\) #1}
\makeatother

\algdef{SE}[DOWHILE]{Do}{doWhile}{\algorithmicdo}[1]{\algorithmicwhile\ #1}%
\newcommand{\X}{\mathcal{X}}
\newcommand{\bxi}{\boldsymbol{\xi}}
\newcommand{\bzt}{\boldsymbol{\zeta}}

\newcommand{\Pc}{ P }
\newcommand{\Xc}{\mathcal{X}}
\newcommand{\ba}{{\mathbf a}}
\newcommand{\bA}{{\mathbf A}}
\newcommand{\bb}{{\mathbf b}}
\newcommand{\bx}{{\mathbf x}}
\newcommand{\by}{{\mathbf y}}
\newcommand{\bz}{{\mathbf z}}
\newcommand{\bv}{{\mathbf v}}
\newcommand{\bu}{{\mathbf u}}
\newcommand{\vertiii}[1]{{\left\vert\kern-0.25ex\left\vert\kern-0.25ex\left\vert #1
		\right\vert\kern-0.25ex\right\vert\kern-0.25ex\right\vert}}

\DeclareMathOperator*{\argmin}{arg\,min}
\newtheorem{assumption}{Assumption}

\allowdisplaybreaks

\jmlrheading{26}{2025}{1-\pageref{LastPage}}{2/24; Revised 4/25}{9/25}{24-0201}{Qihang Lin, Negar Soheili, Runchao Ma, and Selvaprabu Nadarajah}

\ShortHeadings{A Restarting Level Set Method}{Lin, Soheili, Ma, and Nadarajah}
\firstpageno{1}

\begin{document}

\title{An Adaptive Parameter-free and Projection-free Restarting Level Set Method for Constrained Convex Optimization Under the Error Bound Condition}

\author{\name Qihang Lin \email qihang-lin@uiowa.edu\\
       \addr The University of Iowa, Tippie College of Business,  Iowa City, IA 52242
       \AND
       \name Negar Soheili\footnote{Corresponding author} \email nazad@uic.edu \\
       \addr University of Illinois at Chicago, College of Business Administration, Chicago, IL 60607
       \AND
       \name Runchao Ma \email runchao-ma@uiowa.edu\\
       \addr The University of Iowa, Tippie College of Business,  Iowa City, IA 52242
       \AND
       \name Selvaprabu Nadarajah \email selvan@uic.edu \\
       \addr University of Illinois at Chicago, College of Business Administration, Chicago, IL 60607}

\editor{Martin Jaggi}

\maketitle

\begin{abstract}
    Recent efforts to accelerate first-order methods have focused on convex optimization problems that satisfy a geometric property known as error-bound condition, which covers a broad class of problems, including piece-wise linear programs and strongly convex programs. Parameter-free first-order methods that employ projection-free updates have the potential to broaden the benefit of acceleration. Such a method has been developed for unconstrained convex optimization but is lacking for general constrained convex optimization. We propose a parameter-free level-set method for the latter constrained case based on projection-free subgradient method that exhibits accelerated convergence for problems that satisfy an error-bound condition. Our method maintains a separate copy of the level-set sub-problem for each level parameter value and restarts the computation of these copies based on objective function progress. Applying such a restarting scheme in a level-set context is novel and results in an algorithm that dynamically adapts the precision of each copy. This property is key to extending prior restarting methods based on static precision that have been proposed for unconstrained convex optimization to handle constraints. We report promising numerical performance relative to benchmark methods.
\end{abstract}

\begin{keywords}
 level set method, parameter free, projection free, accelerated methods, constrained convex optimization, error bound condition
\end{keywords}

\section{Introduction}\label{sec:intro}
In this paper, we consider a convex optimization problem with inequality constraints:
	\begin{eqnarray}
	f^* := 
	\min_{\bx \in \Xc} \{f(\bx):= f_0(\bx)\quad\text{s.t.}\quad g(\bx) := \max_{i=1,\dots,m}f_i(\bx)\leq 0 \},
	\label{eq:gco}
	\end{eqnarray}
	where $f_i$ for $i=0,1,\dots,m$ are convex real-valued functions and $\Xc\subset\mathbb{R}^n$ is a closed convex set with an efficiently computable projection mapping. Given $\epsilon > 0$, we call $\bar\bx$ an \emph{$\epsilon$-optimal} solution to \eqref{eq:gco} if $f(\bar \bx)-f^*\leq \epsilon$ and \emph{$\epsilon$-feasible}  if $\bar\bx \in \Xc$ and $g(\bar\bx)\leq \epsilon$. There has been a large volume of literature on constrained convex optimization problems in deterministic settings~\citep{yu2017online,bayandina2018mirror,grimmer2018radial,lin2018level,lin2018feasiblelevel,wei2018solving,aravkin2019level,fercoq2019almost, xu2018primal, xu2017first, xu2017global,ene2021,  ito2023, sujanani2024,deng2024} and significant, but comparatively less, activity in stochastic~\citep{lan2016algorithms,yu2017online,lin2019data} settings.

	Recent developments have focused on accelerating first order methods with linear convergence rates (without assuming strong convexity) based on a more general geometric property known as the ``error bound condition (EBC)'', which is 
\begin{eqnarray}
	\label{eq:errorbound}
\dfrac{\text{dist}(\bx, \mathcal{X}^*)^d}{G}  \leq \max \lbrace f(\bx) -f^*,g(\bx)\rbrace,\quad\forall \bx\in \mathcal{X},
	\end{eqnarray}
	for some parameters $G > 0$ and $d\geq 1$, where $\mathcal{X}^*$ is the set of optimal solutions to \eqref{eq:gco} and $\text{dist}(\bx,\mathcal{X}^*)$ denotes the Euclidean distance of $\bx$ to $\mathcal{X}^*$, which we refer to as \emph{solution distance}. EBC provides a lower bound on the maximum of the optimality $(f(x)-f^*)$ and feasibility $(g(x))$ gaps as a function of the solution distance, $d$, and $G$.  A larger such lower bound is desirable, as it implies that additional progress in reducing solution distance yields meaningful improvements in optimality and feasibility. In contrast, if this lower bound is very small, iterations that move $\bx$ closer to $\mathcal{X}^*$ may not result in significant optimality or feasibility progress. Smaller values of $G$ lead to larger lower bounds all else being the same and are thus preferred. For a fixed $G$, we prefer a smaller $d$ because it corresponds to more favorable geometry closer to the optimal solution set, which is typically when the progress of first-order methods slows down. Specifically, when $\bx$ is close to the set of optimal solutions (i.e., is $\text{dist}(\bx,\mathcal{X}^*) < 1$), smaller $d$ leads to a larger lower bound.
    Many constrained convex optimization problems satisfy EBC with exponents $d = 1$ or $d = 2$. For example, piecewise linear programs typically yield $d = 1$, while strongly convex problems yield $d = 2$. Although less commonly discussed in the literature, the EBC framework can also accommodate exponents larger than two. A concrete example and bound are provided in Appendix \ref{Appendix:Example}.
    
We refer to an algorithm as \emph{adaptive} if its convergence rate improves (i.e., it accelerates) as $d$ and $G$ become smaller\footnote{We note that algorithms that exhibit acceleration for only specific values of $d$ (e.g., under strong convexity, which corresponds to $d=2$) would not be considered adaptive under our definition.}. \looseness=-1
Adaptive methods for unconstrained or simply constrained\footnote{Here, simply constrained means $m = 0$ in~\eqref{eq:gco} and $\X$ is a simple set, e.g., $\mathbb{R}^n$, a box or a ball.} convex optimization problems have been explored. These methods use the knowledge of parameters $d$, $G$, and $f^*$ in their step length \citep{polyak1969,iouditski2014primal,bolte2017error, davis2018, wei2018solving,xu2018frank,zhang2019proximal,grimmer2019convergence,necoara2019linear,johnstone2020faster,zhang2020new,renegar2021FeasProb,grimmer2021general,wang2022,grimmer2024optimal} or in the frequency with which they restart \citep{xu2016homotopy,liu2017adaptive,xu2017adaptive, xu2017stochastic,xu2017admm,yang2018rsg,davis2019stochastic, fercoq2019adaptive,yan2019stochastic,fercoq2020restarting,charisopoulos2022,fercoq2022}. Such parameters are unknown in general, and their accurate estimation is challenging. Adaptive algorithms that are also parameter-free (i.e., do not rely on the knowledge of unknown parameters) would therefore broaden or ease the applicability of accelerated methods.  

While a range of parameter-free adaptive methods exist for unconstrained problems, analogous results for constrained settings remain limited. To the best of our knowledge,  \citet{roulet2017sharpness}, \citet{diaz2021}, \citet{ito2021}, and \citet{renegar2018simple} are the only parameter-free adaptive methods that handle unconstrained problems. The restarted method in \citet{roulet2017sharpness} utilizes a logarithmic grid search over $d$ and $G$ to determine the optimal frequency of restart without knowing $d$ and $G$ parameters in advance. While effective, their method is only applicable to unconstrained problems and does not extend to constrained settings. The method in \citet{renegar2018simple} employs a subgradient scheme to solve $K$ copies of the problem with termination precisions $2^{k}\epsilon$ for $k = -1,0,\ldots,K-1$. Copy $k$ communicates its current solution to copy $k-1$ when a restart is triggered. Since each copy operates at a pre-defined precision, this restarting subgradient method (RSG) can be interpreted as a static precision-based restarting scheme. A key feature of RSG is that it tracks the reduction in the objective function (i.e., progress made) to determine when to trigger a restart, as opposed to relying on the distance from $f^*$ (i.e., remaining progress). This structure enables the method to operate without requiring prior knowledge of problem parameters. However, RSG is not designed to handle convex optimization problems with (potentially) complicated constraints. This method has been adapted to bundle methods in \citet{diaz2021}, where similar restart logic is used in a subgradient bundle framework. \citet{ito2021} apply a similar idea for smooth unconstrained convex optimization problems where the restarts are triggered based on the norm of proximal gradients instead of function values.

 Several works, including \citet{nesterov_nonsmoothminimization}, \citet{lin2015adaptive}, \citet{lan2023optimal}, and \citet{sujanani2024}, study parameter-free methods for smooth, strongly convex, unconstrained problems that satisfy the EBC with exponent \(d = 2\). These methods employ accelerated gradient schemes combined with backtracking line search to estimate the gradient Lipschitz constant, using a guess-and-check procedure to achieve optimal complexity without prior knowledge of $G$. In contrast, \citet{davis2022} propose a subgradient method for unconstrained, non-smooth convex problems under the same EBC assumption $(d = 2)$, which achieves \emph{local} nearly linear convergence. Our work extends these results by providing \emph{global} convergence guarantees for general constrained problems.
 
 There are limited adaptive algorithms that solve the constrained convex optimization problem~\eqref{eq:gco}; see \citet{renegar2016efficient}, \citet{yang2017richer}, and \citet{adcock2025restarts}. The algorithm proposed in \citet{renegar2016efficient} addresses the case of \(d = 1\), but it requires prior knowledge of the optimal value \(f^*\). The method in \citet{yang2017richer} applies to general $d$, but assumes knowledge of both $G$ and $d$, among other problem parameters. Moreover, it requires projection onto the set $\mathcal{X} \cap \{ \bx \mid g(\bx) \leq 0 \}$, which is often computationally intractable.
\citet{adcock2025restarts} propose a restarted method that uses a grid search approach - similar to that of \citet{roulet2017sharpness} - to estimate both $d$ and $G$. Unlike \citet{roulet2017sharpness}, however, they introduce a scheduling criterion function to determine the order in which the grid estimates are applied. Their method can be extended to constrained problems by incorporating an exact penalty term with a sufficiently large penalty parameter. However, it remains unclear how such a parameter can be selected while maintaining the method’s parameter-free nature. Furthermore, even with this extension, the complexity of the method in \citet{adcock2025restarts} is higher than ours. For instance, in the non-smooth case with $d = 1$, our method achieves a complexity of $\mathcal{O}\left( \log^2(1/\epsilon) \right)$, while theirs is $\mathcal{O}\left( c^{-1} \log^{2 + c}(1/\epsilon) \right)$ for any $c > 0$.

\emph{The goal of our paper is to develop a first method that is simultaneously adaptive, parameter-free, and projection-free for general convex constraints.}
\footnote{As is common in the literature when using the term projection-free, we do allow projections onto $\mathcal{X}$, which are inexpensive to compute and often available in closed form.}. Extending RSG to handle general convex constraints would be the natural first avenue to consider to achieve our goal. Indeed, the $k$-th copy in such an extension will solve the problem of finding an $2^{k}\epsilon$-optimal and $2^{k}\epsilon$-feasible solution, which entails balancing optimality and feasibility by defining a single objective. Defining such an objective is possible given optimal Lagrange multipliers or $f^*$, which are both unknown. Thus, extending RSG in this manner appears to be challenging. 

\begin{table}[!ht]
\begin{center}
\renewcommand{\arraystretch}{1.2}
\begin{tabularx}{\linewidth}{|>{\centering\arraybackslash}p{25mm}|c|c|>{\centering\arraybackslash}X|}
\hline
& Smooth & Algorithm & Complexity (Function/Gradient Evaluations)\\
EBC & Functions & Parameters & \\
\hline
\hline
\parbox[t]{25mm}{\centering General \\ {\small [This paper]}} & No & None & $\mathcal{O}\left((mM^2 G^{(2/d)})/\epsilon^{(2-2/d)}\log^2(1/\epsilon)\right)$\\
 & Yes & None & $\mathcal{O}\left( m\sqrt{L}\max \left\{ G^{1/d}/\epsilon^{(1/2 - 1/d)}, 1\right\} \log^2(1/\epsilon) \right)$\\
\hline
\parbox[t]{25mm}{\centering General \\ \citep{yang2017richer}, where m = 1$^{\dagger \dagger}$} & No & $d$, $G$, $M$ & $\mathcal{O}\left((M^{4}G^{(2/d)})/\epsilon^{(2-2/d)}\log\left(1/\epsilon\right)\right)$\\
 & Yes & $d$, $G$, $M$, $L$ & $\mathcal{O}\left((M^2G^{(1/d)})/\epsilon^{(1-1/d)}\log\left(1/\epsilon\right)\right.$\\
  & & & + $\left.(\sqrt{ML} G^{(1/d)})/\epsilon^{(1-1/d)}\right)$ \\
 \hline
  \parbox[t]{25mm}{\centering $d=2$\\ \citep{lin2018feasiblelevel}} & No & $M$, $\mu$  & $\mathcal{O}\left(mM^2/\mu\epsilon\right)$\\
 & Yes & $M$, $L$, $\mu$  & $\mathcal{O}\left(m\left(\sqrt{\log(m)}M/\sqrt{\mu\epsilon} + \sqrt{L/\mu}\log\left(1/\epsilon\right)\right)\right)$\\
\hline
\parbox[t]{25mm}{\centering  $d=2$ \\ \citep{xu2017global}}  & Yes & $M$, $L$, $\mu$  & $\mathcal{O}\left( m\left((M+\sqrt{ML})/\epsilon+\sqrt{L/\epsilon}\right)\log(1/\epsilon) \right)$\\
\hline
\parbox[t]{25mm}{\centering  $d=1$  \\ \citep{renegar2016efficient}} & No & $f^*$  & $\mathcal{O}\left(m G^2\log(1/\epsilon)\right)$\\
\hline
\end{tabularx}
\caption{\label{table:complexity} 
Complexity of subgradient-based methods that accelerate under the error bound condition~\eqref{eq:errorbound} for smooth and non-smooth constrained convex optimization problems. The bounds reflect the total number of subgradient, gradient, or function evaluations (including those from constraint functions) required to obtain an $\epsilon$-optimal and $\epsilon$-feasible solution. Here, $M$ denotes an upper bound on the subgradient norm, and $L$ and $\mu$ are the smoothness (Lipschitz) and strong convexity parameters associated with the functions $f_i$, $i = 0,1,\ldots,m$. $^{\dagger \dagger}$ indicates that the method requires projection onto $\mathcal{X}\cap\{\bx \mid g(\bx)\leq 0\}$.
}
\end{center}
\end{table}

We show that our goal can be achieved by developing a restarting level set method. Level-set methods in~\citet{aravkin2019level} and \citet{lin2018level,lin2018feasiblelevel,lin2019data} convert the solution of the constrained optimization problem \eqref{eq:gco} into the solution of a sequence of unconstrained non-smooth optimization problems that depend on a scalar, referred to as the level parameter. The level parameter is chosen by approximately solving a one-dimensional root-finding problem\footnote{An alternative approach to handle constraints is using the radial duality framework of \citet{grimmer2018radial,renegar2016efficient}, where a constrained convex optimization problem is converted to an unconstrained convex optimization problem, but solving this converted problem relies on being able to do exact line searches to evaluate the gauge of the feasible set.}. \citet{aravkin2019level} provides a detailed discussion and complexity analysis of the level set approach. Subgradient algorithms to solve the aforementioned level-set subproblems are desirable as they are easy to implement. However, level set methods that leverage subgradient algorithms depend on unknown parameters, which may include an upper bound on the subgradient norm and/or a smoothness parameter \citep{lin2018level,lin2018feasiblelevel,lin2019data}. Our contributions are the following: 
\begin{itemize} \setlength{\itemsep}{1pt}
\item We introduce a restarting level set method (RLS) that is both parameter-free and projection-free. It maintains $K$ copies, each distinguished by a level parameter and its associated sub-problem. Each restart of copy $k$ results in an update of the level parameters associated with copies $k, k+1, \ldots, K$. The restart of a copy is triggered by the progress made on the sub-problem objective. To the best of our knowledge, this is the first level set method based on projection-free subgradient oracles that is parameter-free, which may also be of broader interest within the literature on level set methods. 
\item RLS can be interpreted as a restarting scheme with dynamic precision, which is conceptually different from RSG. To elaborate, a fixed level parameter implicitly sets a precision for a subproblem. Since RLS updates the level parameter of multiple copies in every restart, it dynamically and implicitly updates the precision associated with copies. Our key algorithmic insights are that (i) triggering restarts based on the progress of the subproblem objective and (ii) dynamically changing the precision of subproblems at each restart by updating the level parameters side step the issues discussed earlier of extending RSG to handle general constrained convex optimization. 
\item We establish that RLS is adaptive in addition to being parameter-free, that is, RLS exhibits acceleration under EBC. Its iteration complexity for non-smooth problems is $\mathcal{O}((M^2G^{2/d})/\epsilon^{(2-2/d)}\log^2(1/\epsilon))$. The adaptive method in~\citet{yang2017richer} has a complexity of $\mathcal{O}((M^4G^{2/d})/\epsilon^{(2-2/d)}\log(1/\epsilon))$ but it depends on unknown parameters and requires a non-trivial projection onto $\mathcal{X}\cap\{\bx|g(x)\leq0\}$. The additional $\log\left(\frac{1}{\epsilon}\right)$ factor in the RLS complexity can be interpreted as the cost of both relaxing parameter dependence and 
employing simple subgradient oracles. Despite this worsening, the dependence of RLS on $M$ is better. Table \ref{table:complexity} shows that an analogous property holds in the smooth setting when comparing RLS with the approach in~\citet{yang2017richer}. Although perhaps unfair to RLS, Table \ref{table:complexity} also compares its computational complexity to algorithms designed for specific values of $d$, all of which depend on unknown parameters. RLS worsens by the same logarithmic factor of $\log\left(\frac{1}{\epsilon}\right)$ relative to these more specialized algorithms and sometimes improves on the dependence with respect to $G$.
\item  We evaluate the performance of RLS on two types of constrained classification problems: (i) fairness-constrained classification, where we consider both non-smooth and smooth formulations, and (ii) non-smooth Neyman--Pearson classification with constraints on the Type II error control. In the fairness-constrained setting, we compare RLS with three parameter-free but non-adaptive benchmarks in the non-smooth case - the feasible level set method from \citet{lin2018feasiblelevel}, the subgradient method from \citet{yu2017online}, and the switching subgradient method from \citet{bayandina2018mirror} - and with the feasible level set method and the constraint extrapolation (ConEx) method from \citet{boob2023stochastic} in the smooth case. For the Neyman--Pearson classification problem, we compare RLS against the same three non-smooth baselines. These experiments include large-scale real-world datasets with up to 72{,}000 data points and 47{,}000 features. Through the experiments, we find that the adaptive nature of RLS enables faster convergence to an $\epsilon$-optimal and $\epsilon$-feasible solution than existing benchmarks in various regimes.

\end{itemize}

This paper is organized as follows: In Section~\ref{sec:LSM}, we discuss a general level-set method for solving~\eqref{eq:gco} and discuss its convergence rate. In Section~\ref{section:main_algorithm}, we present RLS and analyze its oracle complexity. In Section~\ref{sec:foms}, we provide first order oracles for use with RLS in the non-smooth and smooth settings, also analyzing the total complexity of RLS in each case. In Section~\ref{sec:exp}, we numerically compare RLS to other benchmarks. We conclude in Section~\ref{sec:conc}.

\section{Level Set Method}\label{sec:LSM}
Throughout the paper, we assume that there exists a computable feasible solution to \eqref{eq:gco} that strictly satisfies the constraint $g(\bx)\leq 0$. This assumption is formalized below. 
	\begin{assumption}
		There exists a computable solution $\tilde{\bx}$ such that $\tilde{\bx} \in \mathcal{X}$ and $g(\tilde{\bx}) < 0$. 
		\label{assume:strictfeasible}
	\end{assumption}
	
	We define the \emph{level-set function} corresponding to \eqref{eq:gco} as:
	\begin{align}
	H(r) := \min_{\bx \in \Xc} P(\bx;r),
	\label{preliminaries:levelset_function}
	\end{align}
	where $r \in \mathbb{R}$ is called a  \emph{level parameter} and 
	\begin{align}
	\label{preliminaries:Pfunction}
	P(\bx;r):= \max \lbrace f_0(\bx) -r, f_1(\bx),\dots,f_m(\bx) \rbrace=\max \lbrace f(\bx) -r,g(\bx)\rbrace.
	\end{align}
	
Below, we summarize some of the well-known properties of $H(r)$ that we use in our algorithm and analyses. 
	
	\begin{lemma}[\citet{lin2018feasiblelevel, nesterov2018lectures}]
		The function $H(r)$ defined in \eqref{preliminaries:levelset_function} has the following properties:
		\begin{enumerate} \setlength{\itemsep}{1pt}
			\item $H(r)$ is non-increasing and convex in $r$;
			\item $H(f^*) = 0$;
			\item $H(r) > 0$ when $r < f^*$ and $H(r) < 0$ when $r > f^*$. 
			\item $H(r) - \delta \leq H(r+\delta) \leq H(r)$ for any $\delta \geq 0$.
		\end{enumerate}
		\label{lemma:levelset_properties}
	\end{lemma}
Lemma \ref{lemma:levelset_properties} indicates that $f^*$ is the unique root to the root-finding problem $H(r) = 0$ under Assumption~\ref{assume:strictfeasible}. Notice that Assumption~\ref{assume:strictfeasible} guarantees the existence of strong duality for problem~\eqref{eq:gco}. Without the assurance of strong duality, the root of $H(r)$ might not lead to $f^*$ (refer to~\citet{Friedlander2020} for further insights). Moreover, as shown in Lemma \ref{lemma:solutionquality}, it is possible to find an $\epsilon$-feasible solution and use the level parameter $r$ to bound the optimality gap.
\begin{lemma}
    Given $\epsilon > 0$, if $\bx\in\mathcal{X}$ satisfies $\Pc(\bx;r)\leq \epsilon$ for some $r$, then $\bx$ is $(r-f^*+\epsilon)$-optimal and $\epsilon$-feasible. Consequently, when $r<f^*$, this solution $\bx$ is $\epsilon$-optimal and $\epsilon$-feasible.
\label{lemma:solutionquality}
\end{lemma}
Hence to tackle~\eqref{eq:gco}, a level-set method applies a root-finding scheme to solve $H(r) = 0$ with the goal of finding $r = f^*$. 
Once such $r$ is found, then $\Pc(\bx;r)\leq \epsilon$  implies that $\bx$ is a nearly optimal and nearly feasible solution to \eqref{eq:gco}. Notice that for any $\bx^*\in \mathcal{X}^*:=\argmin_{\bx \in \Xc} \Pc(\bx;f^*)$ we have $\Pc(\bx^*;f^*) = 0$. 
	
Algorithm~\ref{algo:FLS} presents a level set method that generates a sequence of level parameters $\{r_k\}_{k\geq 0}$ converging to $f^*$. 
Given $\alpha\in (0,1)$, each update of $r_{k+1}=r_k+\alpha\Pc(\bx_k;r_k)$ requires a vector $\bx_k\in\mathcal{X}$ satisfying $\alpha\Pc(\bx_k;r_k)<  f^*-r_k$, that is, condition \eqref{eq:xkquality}. The convergence of Algorithm \ref{algo:FLS} relies on $r_k < f^*$. Under this condition, it follows from Lemma~\ref{lemma:levelset_properties}(4) and $\alpha < 1$ that $\bx_k$ can be computed by solving $\min_{\bx\in\mathcal{X}}\Pc(\bx;r_k)$ because $\min_{\bx\in\mathcal{X}}\Pc(\bx;r_k)=H(r_k)\leq f^*-r_k < (f^*-r_k)/\alpha$. Moreover, the level parameter $r_k$ converges to $f^*$ from below causing the quantity $\Pc(\bx_k;r_k)\geq 0$ to converge to zero. Thus, Algorithm~\ref{algo:FLS} terminates and the final solution is both $\epsilon$-optimal and $\epsilon$-feasible by Lemma~\ref{lemma:solutionquality}.

	\begin{algorithm}[t]
		\caption{Level Set Method}
		\label{algo:FLS}
		\begin{algorithmic}[1]
			\State{{\bf Input:} $\alpha\in(0,1)$, $\epsilon >0$, and $r_0<f^*$.}
			\State{\bf Initialize:} $k = -1$.
			\Do
				\State $k = k +1$.
				\State Solve $\min_{\bx\in\mathcal{X}}\Pc(\bx;r_k)$ to obtain $\bx_k\in\mathcal{X}$ such that
				\begin{equation}
				\label{eq:xkquality}
				\alpha\Pc(\bx_k;r_k)< f^*-r_k.
				\end{equation} 
				\State Set $r_{k+1}=r_{k}+\alpha\Pc(\bx_k;r_k)$.
			\doWhile{$\Pc(\bx_k;r_k)> \epsilon$.}
			\State{{\bf Output:} $\bx_k$.}
		\end{algorithmic}
	\end{algorithm}
	
Next, we present analysis that shows the correctness and complexity of Algorithm \ref{algo:FLS}. This analysis depends on a condition measure that is defined as
\begin{eqnarray}\label{eq:conditionnumber}
    \theta:=-\lim_{r\rightarrow f^{*-}}\frac{H(r)-H(f^*)}{r-f^*}=\lim_{r\rightarrow f^{*-}}\frac{H(r)}{f^*-r}.
\end{eqnarray}
The following lemma sheds light on some geometric properties of the condition measure $\theta$. We skip proving this lemma as it directly follows from the definition of $\theta$ in~\eqref{eq:conditionnumber} and the first and fourth properties in Lemma~\ref{lemma:levelset_properties}.
\begin{lemma}
	It holds that $0<\theta \leq \dfrac{H(r)}{f^*-r}\leq 1$ for any $r < f^*$ and $0\leq\dfrac{H(r)}{f^*-r}\leq\theta$ for  any $r > f^*$.\looseness = -1 \label{theorem:levelset_method_theta}
\end{lemma}
Lemma \ref{theorem:levelset_method_theta} establishes bounds on $H(r)$. In particular, when $r_k$ is less than $f^*$, there is a natural lower bound of zero. 
	
If $r_0 < f^*$, Lemma \ref{lem:r.update.guarantee} shows that $r_k$ remains less than $f^*$ at each iteration of Algorithm \ref{algo:FLS} and that it converges to $f^*$ at a geometric rate.
\begin{lemma}\label{lem:r.update.guarantee} Suppose $\bx_k\in\mathcal{X}$ and $r_k$ satisfy $r_k< f^*$ and \eqref{eq:xkquality}. The update $r_{k+1}=r_k+\alpha\Pc(\bx_k;r_k)$ in Algorithm~\ref{algo:FLS} guarantees (i) $r_{k+1}< f^*$; and (ii)  $f^* - r_{k+1} \leq (1-\alpha\theta)(f^*-r_k)$.
\end{lemma}
	
\proof{} Fix $k$. Given $r_k$ and $\bx_k$ satisfying $\alpha\Pc(\bx_k;r_k)<  f^*-r_k$, we have
\begin{align*}
	r_{k+1}  = r_k + \alpha\Pc(\bx_k;r_k)<r_k + f^* - r_k = f^*.
\end{align*}
In addition, it follows that
\begin{align}
	f^*-r_{k+1} = f^*-r_k-\alpha\Pc(\bx_k;r_k) \leq f^*-r_k-\alpha H(r_k) \leq  \left(1-\alpha\theta\right)(f^*-r_k), \label{theroem:levelset_proof}
\end{align}
where the first inequality holds because $0\leq H(r_k)\leq\Pc(\bx_k;r_k)$ and the second one because $H(r_k) \geq \theta(f^* - r_k)$ for $r_k< f^*$ by Lemma~\ref{theorem:levelset_method_theta}. 
\hfill\BlackBox

In Theorem~\ref{theorem:levelset_method}, we characterize the number of outer iterations needed by Algorithm~\ref{algo:FLS} to find an $\epsilon$-optimal and $\epsilon$-feasible solution to~\eqref{eq:gco}. 
\begin{theorem}
	Given $\epsilon >0$, $\alpha\in(0,1)$, and $r_0 < f^*$, Algorithm~\ref{algo:FLS} terminates in at most 
	\begin{align}\label{eq:trueK}
		\hat K:=\left\lceil \frac{1}{\alpha\theta} \ln \left( {\frac{f^*-r_0}{\alpha\epsilon}} \right)  \right\rceil,
	\end{align}
	iterations and returns an $\epsilon$-optimal and $\epsilon$-feasible solution to \eqref{eq:gco}.
	\label{theorem:levelset_method}
\end{theorem}
\proof{}
		Since Algorithm~\ref{algo:FLS} starts with $r_0< f^*$, Lemma~\ref{lem:r.update.guarantee} indicates that $r_k <f^*$ and~\eqref{theroem:levelset_proof} holds at each iteration $k\geq 0$. Since $\alpha\theta\in (0,1)$, recursively applying~\eqref{theroem:levelset_proof} for $k=0,1,\dots$ shows that $0\leq f^*-r_k \leq \left(1-\alpha\theta\right)^k(f^*-r_0)$ for any $k\geq 0$. Hence, it follows from inequality~\eqref{eq:xkquality} that $\Pc(\bx_k;r_k)<\frac{1}{\alpha}(f^*-r_k)\leq \frac{1}{\alpha}\left(1-\alpha\theta\right)^k(f^*-r_0)$. When $k = \left\lceil \frac{1}{\alpha\theta} \ln \left( {\frac{f^*-r_0}{\alpha\epsilon}} \right)  \right\rceil$, we get $\Pc(\bx_k;r_k)\leq\epsilon$  and the algorithm stops.  Since $r_k < f^*$,  Lemma~\ref{lemma:solutionquality} then guarantees that the solution $\bx_k$ returned at termination is an $\epsilon$-optimal and $\epsilon$-feasible solution.
\hfill\BlackBox

\begin{remark}
The $1/\theta$ factor in the definition of $\hat{K}$ arises because we leverage the natural lower bound of zero on $H(r)$, as discussed after Lemma \ref{theorem:levelset_method_theta}. This choice results in the specific form of condition \eqref{eq:xkquality}, which is a simplification of the analogous condition in \citet{aravkin2019level}, where the authors use non-trivial lower and upper bounds on $H(r)$. The lower bound used in \citet{aravkin2019level} is easily computable only when the problem has certain structures. While the number of iterations in the algorithm from \citet{aravkin2019level} becomes independent of $\theta$ when such bounds are efficiently calculable, obtaining such bounds for a broader class of problems poses a formidable challenge. We opt for our simpler approach, which proves implementable across a more general class of problems and contributes to a clearer exposition of the fundamental concepts underlying RLS in the subsequent section.
\end{remark}
 
\section{Restarting Level Set Method}	\label{section:main_algorithm}
We define a sequence of pairs $\{(\bx_k,r_k)\}_{k=0}^K$ with $r_0<f^*$ and $\bx_k\in\mathcal{X}$  as a \emph{level set sequence} of length $K+1$ if $r_{k+1}=r_{k}+\alpha\Pc(\bx_k;r_k)$ holds for some $\alpha\in(0,1)$ and each $k=0,1,\dots,K-1$. It follows from Theorem~\ref{theorem:levelset_method} that the level set sequence generated in Algorithm~\ref{algo:FLS} converges to an $\epsilon$-optimal and $\epsilon$-feasible solution when (i) $\bx_k$ satisfies \eqref{eq:xkquality} and (ii) $K\geq \hat K$, where $\hat K$ is defined in~\eqref{eq:trueK}. 
One way to find $\bx_k\in\mathcal{X}$ satisfying \eqref{eq:xkquality} is to solve the \emph{level-set subproblem}~\eqref{preliminaries:levelset_function} with $r=r_k$, which we denote by $\mathtt{LSP}(r_k)$:\looseness = -1
	\begin{align}
	\min_{\bx\in\mathcal{X}}\Pc(\bx;r_k).
	\label{preliminaries:levelset_function_k}
	\end{align}
Let $\mathtt{fom}$ represent the first-order method applied to \eqref{preliminaries:levelset_function_k}. Although there are many choices for $\mathtt{fom}$, it is difficult to numerically verify \eqref{eq:xkquality} because $f^*$ is unknown. As a result, we are not able to terminate $\mathtt{fom}$ at the right time. If $\mathtt{fom}$ is terminated too soon, the returned solution $\bx_k$ will not satisfy \eqref{eq:xkquality} and Algorithm~\ref{algo:FLS} may not converge. If $\mathtt{fom}$ is terminated too late, Algorithm~\ref{algo:FLS} will converge but consume longer run time than it actually needs. 

To address this issue, we parallelize the sequential approach in Algorithm \ref{algo:FLS} by simultaneously solving multiple level-set subproblems with potentially different level parameters and restarting the solution of each such subproblem if a predetermined amount of progress in reducing its objective function has been made. We refer to this procedure as the \emph{restarting level set} (RLS) method and describe its core idea next. It requires maintaining $K+1$ (with $K \geq \hat{K}$) first order method instances denoted by $\mathtt{fom}_0, \mathtt{fom}_1, \ldots, \mathtt{fom}_K$. Given  $\bx_{\text{ini}}\in\mathcal{X}$ and $r_{\text{ini}}< f^*$, we initialize the level set sequence $\{(\bx_k^{(0)},r_k)\}_{k=0}^K$ with $\bx_{k}^{(0)}=\bx_{\text{ini}}$ for $k=0,1,\dots,K$, $r_{0}=r_{\text{ini}}$, and $r_{k+1}=r_k+\alpha\Pc(\bx_k^{(0)};r_k)$ for $k=1,\dots,K-1$. We then apply $\mathtt{fom}_k$ to solve $\mathtt{LSP}(r_k)$ starting from $\bx_{k}^{(0)}$ for $k=0,1,\dots,K$ and denote by $\bx_{k}^{(t_k)}$ the solution computed by $\mathtt{fom}_k$  after $t_k$ iterations.  
The $\mathtt{fom}$ instances communicate with each other via restarts. We initiate a \emph{restart} at instance $k'$ when $\Pc(\bx_{k'}^{(0)};r_{k'})\geq 0$, and the objective function of $\Pc(\bx, r_{k'})$ experiences a significant reduction. Specifically, the restart occurs when the solution $\bx_{k'}^{(t_{k'})}$ satisfies the condition
	\begin{equation}
	\label{eq:reduceB}
	\Pc(\bx_{k'}^{(t_{k'})};r_{k'})\leq B \Pc(\bx_{k'}^{(0)};r_{k'}),
	\end{equation}
where $B\in(0,1)$ is a constant factor. The sequence of steps associated with a restart is shown in Definition \ref{def:restart}. 
\begin{definition}[Restart at $k'$] Set $\bx_{k'}^{(0)} = \bx_{k'}^{(t_{k'})}$ and $r_{k+1}=r_k+\alpha\Pc(\bx_k^{(0)};r_k)$ for all $k = k', k'+1,\ldots,K-1$. Restart $\mathtt{fom}_{k}$ from $\bx_{k}^{(0)}$ for  $k = k', k'+1,\ldots,K$.
\label{def:restart}
\end{definition}

In other words, initiating a restart at $k' < K$, we reset $\mathtt{fom}_{k'}$ to begin at an updated initial solution equal to $\bx_{k'}^{(t_{k'})}$ and solve the same level set subproblem $\mathtt{LSP}(r_{k'})$ as before the restart because $r_{k'}$ is not updated. In contrast, for indices $k = k'+1,\ldots,K$, $\mathtt{LSP}(r_{k})$ changes because the restart changes $r_k$. This new level set subproblem is solved by $\mathtt{fom}_{k}$ starting from the same initial solution $\bx_{k}^{(0)}$ that was used before the restart. If a restart is initiated at $K$, we only perform the update $\bx_{K}^{(0)} = \bx_{K}^{(t_K)}$. Restarts initiated at index $k'$ have no effect on $\mathtt{fom}_{k}$ for $k=0,1,\dots,k'-1$, so those $\mathtt{fom}$s will continue their iterations. 
Lemma \ref{lem:lssequence} follows from our restarting updates. 
\begin{lemma}\label{lem:lssequence}
	The level set sequence $\{(\bx_k^{(0)},r_k)\}_{k=0}^K$ remains a level set sequence after each restart. \looseness=-1
\end{lemma}	
\begin{remark}We highlight that restarts dynamically modify the precision with which subproblems need to be solved by $\mathtt{fom}$ instances to satisfy \eqref{eq:xkquality}. Specifically, the precision of a solution at instance $k$ can be viewed as the deviation of $\Pc(\bx_{k}^{(t_{k})};r_{k})$ from the optimal subproblem objective of $H(r_k)$. This difference satisfies 
\[\Pc(\bx_{k}^{(t_{k})};r_{k}) - H(r_k) \leq \dfrac{1}{\alpha} (f^*-r_k) - \theta (f^*-r_k) \leq (1 - \alpha\theta) \Pc(\bx_{k}^{(0)};r_{k}),\] 
where the first inequality follows from Lemma~\ref{theorem:levelset_method_theta} and the fact that the solution $x_k^{(t_k)}$ satisfies \eqref{eq:xkquality} and the second inequality holds since the initial solution $\bx_{k}^{(0)}$ does not satisfy \eqref{eq:xkquality}. i.e. $\alpha \Pc(\bx_{k}^{(0)};r_{k}) > f^* - r_k.$
It follows from these inequalities that the precision $\Pc(\bx_{k}^{(t_{k})};r_{k}) - H(r_k)$ is $\Theta(\Pc(\bx_{k}^{(0)};r_{k}))$. A restart at $k'$ changes $\bx_{k}^{(0)}$ or $r_{k}$ of instances $k \geq k'$ and hence the value of $\Pc(\bx_{k}^{(0)};r_{k})$ at these instances. Therefore, RLS can be interpreted as an approach that dynamically modifies the precision required to satisfy \eqref{eq:xkquality} at a subset of instances, each time a restart is initiated as the result of there being sufficient progress in the subproblem objective function at some instance, that is, condition \eqref{eq:reduceB} is satisfied. 
\end{remark}

Well-known choices of $\mathtt{fom}$ (e.g., subgradient method) can be shown to meet the restart condition \eqref{eq:reduceB} after a finite number of iterations starting from $\bx^{(0)}$ when $\Pc(\bx^{(0)};r)\geq 0$ is sufficiently large. In particular, when $\Pc(\bx^{(0)};r)> \epsilon$, the number of iterations taken by $\mathtt{fom}$ to find a desirable solution can be bounded by an integer $n^{\mathtt{fom}}$ that is independent of $r$ and $\bx^{(0)}$. We summarize this property in Assumption \ref{assum:FOM} and take it to be true in this section. We will present $\mathtt{fom}$ choices that satisfy this assumption in \S\ref{sec:foms}.
	
\begin{assumption}\label{assum:FOM}
Consider $\alpha$,  $B$, $\epsilon$, $r$, and $\bx^{(0)}\in\mathcal{X}$ such that $0<\alpha<B<1$, $\epsilon>0$, $r < f^*$, and $\Pc(\bx^{(0)};r)> \epsilon$. There exists an $\mathtt{fom}$ initiated at $\bx^{(0)}$ which guarantees that the inequality $\Pc(\bx^{(t)};r) \leq B  \Pc(\bx^{(0)};r)$ holds after $t\leq n^{\mathtt{fom}}$ iterations if \eqref{eq:xkquality} is not satisfied by then, where $n^{\mathtt{fom}}$ is an integer depending on $\alpha$, $B$ and $\epsilon$ and independent of  $\bx^{(0)}$ and $r$.
\end{assumption} 

To interpret RLS, it is useful to think of instance $\mathtt{fom}_k$ as being a proxy for the first-order method applied to solve $\mathtt{LSP}(r_k)$ in iteration $k$ of Algorithm \ref{algo:FLS}. By Lemma~\ref{lem:lssequence}, the relationship between the $r_k$ values across iterations $k$ in Algorithm \ref{algo:FLS} are also maintained for the $r_k$ values across instances $\mathtt{fom}_k$ in RLS but some $r_k$ may not be computed based on a solution that satisfies \eqref{eq:xkquality}. As a result, even though the initial level set parameter is less than $f^*$, subsequent updates can push $r_k \geq f^*$. Theorem \ref{thm:keythm} states key properties of RLS that allow it to nevertheless converge to a near-optimal and near-feasible solution, as we will discuss shortly.
	
\begin{theorem}
	\label{thm:keythm}
	Given $\epsilon >0$ and $\alpha\in (0,1)$, consider a level set sequence $\{(\bx_k^{(0)},r_k)\}_{k=0}^K$ with $K\geq \hat K$ where $\hat K$ is defined in \eqref{eq:trueK}. At least one of the following two statements holds:
	    \begin{itemize}
			\item[A.] There exists an index $k^*\in\{0,1,\dots,K\}$ such that $r_{k^*}< f^*$ and $\alpha P(\bx^{(0)}_{k^*}; r_{k^*})\geq f^*- r_{k^*}$;\looseness = -1
			\item[B.] The solution $\bx^{(0)}_{K}$ is an $\epsilon$-optimal and $\epsilon$-feasible solution.
		\end{itemize}
		Moreover, if statement A holds, it follows that the index $k^*$ is unique and we have $r_0<\dots<r_{k^*}< f^*$ and $\alpha P(\bx_{k}^{(0)}; r_{k})< f^*- r_{k}$ for $k=1,2,\ldots,k^*-1$.
\end{theorem}
\proof{}
Suppose statement A does not hold for any index $k\leq K$. We show the statement $B$ holds. By the definition of a level set sequence, we have $r_0<f^*$ and $r_{k+1}= r_k +\alpha\Pc(\bx_{k}^{(0)} ;r_{k})$ for all $k$. Since $\alpha\Pc(\bx_k^{(0)};r_k) < f^* - r_k$ (by our assumption that $A$ does not hold), Lemma~\ref{lem:r.update.guarantee} guarantees that $r_1< f^*$. Using the same argument and by induction, Lemma~\ref{lem:r.update.guarantee} implies  $r_k< f^*$ and $f^*-r_k \leq \left(1-\alpha\theta\right)^k(f^*-r_0)$ for any $k$. Therefore, we have $\Pc(\bx_{K}^{(0)};r_{K})<1/\alpha(f^*-r_{K})\leq \frac{1}{\alpha}\left(1-\alpha\theta\right)^{{K}}(f^*-r_0)$. Since $K\geq \hat K$ with $\hat K$ defined in \eqref{eq:trueK}, we have $ P(\bx_{K}^{(0)}; r_{K})\leq \epsilon$ which indicates $\bx_{K}^{(0)}$ is  an $\epsilon$-optimal and $\epsilon$-feasible solution by Lemma~\ref{lemma:solutionquality}.
		
Next, we prove the rest of the conclusion. Let $k^*$ denote the smallest index satisfying the conditions of statement A. If $k^* = 0$ the conclusion is trivial since $r_0 < f^*$. Suppose $k^* >0$. Since $r_0 < f^*$, following the definition of $k^*$, we must have $\alpha P(\bx_0^{(0)};r_0) < f^* - r_0$. Lemma~\ref{lem:r.update.guarantee} then guarantees that $r_1 < f^*$. Applying the same argument, we can show $r_2, ..., r_{k^*-1}$ are all less than $f^*$ ( i.e. $r_k< f^*$) and hence the inequality $\alpha P(\bx_k^{(0)};r_k) < f^* - r_k$ must hold for any $k< k^*$. In addition, recall that 
$\Pc(\bx;r)\geq H(r)> 0$ for any $r<f^*$. Hence, the relationship $r_{k+1} = r_k +\alpha P(\bx_k^{(0)};r_k)$ shows $r_{k+1} > r_k$ for any $k < k^*$.  The index $k^*$ in $A$ is unique because it is straightforward to see $r_{k^*+1} \geq f^*$ which follows from the definition of $r_{k^*+1}$ and the inequality $\alpha P(x_{k^*}^{(0)}, r_{k^*}) \geq f^* - r_{k^*}$. Moreover, once there exists an index $k$ such that $r_k\geq  f^*$, we must have $r_i \geq f^*$ for all $i\geq k$ since 
		$$
		r_{i+1}-f^*=r_{i}-f^*+\alpha P(\bx_i^{(0)};r_i)\geq r_{i}-f^*+\alpha H(r_i)\geq  (1-\alpha\theta)(r_i-f^*) \geq 0,
		$$
		where the first inequality holds by the definition of $H(\cdot)$ and the second by Lemma~\ref{theorem:levelset_method_theta}. Therefore, given that $r_{k^*+1}\geq f^*$, we must have $r_k\geq f^*$ for any $k\geq k^*+1$. 
\hfill\BlackBox

Given a level set sequence, we refer to the unique index $k^*$ in Theorem \ref{thm:keythm} as the \emph{critical index}. At this index, we have $r_{k^*} < f^*$, but condition \eqref{eq:xkquality}  does not hold. For all indices $k$ that precede the critical index, we have $r_k < f^*$ and condition \eqref{eq:xkquality} holds. Moreover, if  $K\geq \hat K$ and $k^*$ does not exist (intuitively $k^* \geq K+1$), it implies the solution $\bx_{K}^{(0)}$ is an $\epsilon$-optimal and $\epsilon$-feasible. Therefore, we would like to increase the critical index with restarts so that it eventually exceeds $K$. Ideally, if one could initiate restarts at the critical index $k^*$ sufficiently many times, it would reduce $\Pc(\bx_{k^*}^{(0)};r_{k^*})$ by updating $\bx_{k^*}^{(0)}$ while ensuring $r_{k^*} < f^*$ (since $r_{k^*}$ is not updated when a restart is initiated at $k^*$ by Definition \ref{def:restart}). Then one would expect condition \eqref{eq:xkquality} will hold at $k^*$. Moreover, when this happens, $r_{k^*+1} < f^*$ by Lemma \ref{lem:r.update.guarantee}. Hence, repeated restarts at the critical index should intuitively increase its value. 

Unfortunately, this intuition is of little algorithmic use since the critical index is defined in terms of $f^*$, which is unknown. Therefore, RLS instead executes restarts at the smallest index at which condition \eqref{eq:reduceB} is satisfied. This index could be smaller than or greater than the critical index $k^*$. We refer to as {\emph{desirable restarts}} the ones initiated at index $k$ less than or equal to $k^*$ because such restarts may update $r_{k^*}$ or $\bx_{k^*}^{(0)}$. By  Assumption~\ref{assum:FOM} and Theorem~\ref{thm:keythm}, a desirable restart must be initiated at $\mathtt{fom}_{k}$ with some $k\leq k^*$ in no more than $n^{\mathtt{fom}}$ iterations, unless $\Pc(\bx_{k^*}^{(0)};r_{k^*})\leq\epsilon$. In the latter case the solution $\bx_{k^*}^{(0)}$ is $\epsilon$-optimal and $\epsilon$-feasible. The number of desirable restarts can be bounded by $\mathcal{O}(\log^2(1/\epsilon))$ as we will establish later. Such restarts that occur before $k^*$ may already increase the critical index. If this does not happen, desirable restarts will start happening at $k^*$ and increase the critical index.

Algorithm \ref{algo:RestartingLS} formalizes the steps of RLS. The inputs to this algorithm include the number of $\mathtt{fom}$ instances $K$, a total budget $I$ on the number of $\mathtt{fom}$ iterations, a level parameter $r_{\text{ini}} < f^*$, a strictly feasible solution $\bx_{\text{ini}}\in\mathcal{X}$, and parameters $\alpha$, $B$, and $\epsilon$. By Assumption \ref{assume:strictfeasible}, we can set $\bx_{\text{ini}}=\tilde\bx$. For $K$, an ideal value to use is $\hat K$ defined in \eqref{eq:trueK}. However, since the parameters $\theta$ and $f^*$ in this bound are unknown, we compute a bound on $\hat K$ by using a strictly feasible solution. In particular, we use
\begin{equation}\label{eq:tildeK}
\tilde K := \left\lceil \frac{1}{\alpha\tilde\theta(r)} \ln \left( {\frac{\tilde r-r_{0}}{\alpha\epsilon}} \right)  \right\rceil,
\end{equation}
 with $\tilde r := f(\tilde \bx) - g(\tilde\bx)$ and $\tilde\theta(r) : = g(\tilde\bx)/(r - \tilde r)$ for $r<f^*$ and $\tilde \bx$ from Assumption~\ref{assume:strictfeasible}. We show below $\tilde K$ is an upper bound on $\hat K$. We set $r  = r_\text{ini}$ in definition of $\tilde \theta(r)$ when executing Algorithm~\ref{algo:RestartingLS}. \looseness = -1

	\begin{lemma}\label{lemma:approximate_condition_number} Let $\tilde \bx$ be the strictly feasible solution in Assumption \ref{assume:strictfeasible} and $\tilde r:=f(\tilde{\bx})-g(\tilde{\bx})$ and $\tilde \theta(r):= g(\tilde{\bx})/(r-\tilde r)$ for $r<f^*$. Then $\tilde K \geq \hat K$ for $\tilde K$ and $\hat K$ respectively defined in \eqref{eq:tildeK} and \eqref{eq:trueK}.\looseness = -1
\end{lemma}
The proof of this lemma is provided in the Appendix to maintain a coherent flow here.

The initialization step in Algorithm~\ref{algo:RestartingLS} assigns $\bx_{\text{ini}}$ as starting solutions to all $\mathtt{fom}$ instances (i.e. $\bx_{k}^{(0)} = \bx_\text{ini}$ for $k = 0, \ldots, K$) and sets the level parameters used in each $\mathtt{fom}_k$ using $r_{k+1} = r_k + \alpha\Pc(\bx_{k}^{(0)}; r_k)$ starting at $r_0 = r_\text{ini}$.  Algorithm \ref{algo:RestartingLS} runs for a pre-specified number of $\mathtt{fom}$ iterations. At each iteration $i$, it runs $K+1$ $\mathtt{fom}$ iterations simultaneously, one for each instance, until the condition \eqref{eq:reduceB} holds for some index. It then finds a smallest index $k'$ for which $\Pc(\bx_{k'}^{(0)}; r_{k'})\geq 0$ and \eqref{eq:reduceB} hold. A restart is executed at $k'$. If the solution $\bx_{k'}^{(t_{k'})}$ is $\epsilon$-feasible and has a better objective function value than the best $\epsilon$-feasible solution at hand, this solution is updated. Once the total of $I$ iterations have been reached, the $\epsilon$-feasible solution with the least $f_0(\bx)$ is output.
	
 	\begin{algorithm}[h!]
		\caption{Restarting Level Set Method}
		\label{algo:RestartingLS}
		\begin{algorithmic}[1]
			\State{{\bf Input:} $K = \tilde K, I\in \mathcal{Z}_+$, $\alpha\in(0,1)$, $B\in (\alpha,1)$, $\epsilon >0$,  $r_{\text{ini}}< f^*$ and a strictly feasible $\bx_{\text{ini}}\in\mathcal{X}$}.
			\State{\bf Initialization:} Set  $\bx_k^{(0)} = \bx_{\text{ini}}$, $t_k = 0$ for $k=1,\dots,K$, $r_0 = r_{\text{ini}}$, $r_{k+1}=r_k+\alpha\Pc(\bx_k^{(0)};r_k)$ for $k=0,1,\dots,K-1$, and $\bx_{\text{best}} = \bx_{\text{ini}}$.\looseness = -1
			\For{$i=1,2,\dots,\lceil I/(K+1)\rceil$}
				{\color{gray}\LineComment{----- Execute $K$ iterations, one in each $\mathtt{fom}$ instance -----}}
			 	\State Run $\mathtt{fom}_k$, $k = 0,1,\ldots,K$, for one iteration each to obtain $\bx^{(t_1+1)}_1, \bx^{(t_2+1)}_2, \ldots, \bx^{(t_K+1)}_K$.
			 	\State $t_k = t_k + 1$, $k = 1,\ldots, K$.
			 	{\color{gray}\LineComment{----- Check if~\eqref{eq:reduceB} holds and find index $k'$ to initiate restart -----}}
			 \If{$\Pc(\bx_k^{(0)};r_k) \geq 0$ and $\Pc(\bx^{(t_k)}_k;r_k) \leq  B\Pc(\bx_k^{(0)};r_k)$ for some $k = 1,\ldots,K$}
			 \State Find the smallest index $k'\geq 1$ such that $\Pc(\bx_{k'
			 }^{(0)};r_{k'}) \geq 0$ and $\Pc(\bx^{(t_{k'})}_{k'};r_{k'}) \leq B\Pc(\bx_{k'}^{(0)};r_{k'})$. 
		 	\State Set $\bx'_{k'}:=\argmin_{\bx\in\left\{\bx^{(t_{k'})}_{k'}, \bx^{(0)}_{1}, \bx^{(0)}_{2},\dots, \bx^{(0)}_{K} \right \}}\Pc(\bx;r_{k'})$.
			 \EndIf
			 {\color{gray}\LineComment{----- Execute restart at $k'$ -----}}
			 \State $\bx_{k'}^{(0)} \leftarrow \bx'_{k'} $.
			 \State $t _k= 0$ for $k=k',k'+1,\dots,K$.
			\For{$k=k',...,K-1$}
			\State Set $\bx_{k}^{(0)}:=\argmin_{\bx\in\left\{\bx^{(0)}_{1}, \bx^{(0)}_{2},\dots, \bx^{(0)}_{K} \right \}}\Pc(\bx;r_{k})$.
			\State $r_{k+1}\leftarrow r_k +\dfrac{\alpha}{2}\Pc(\bx_{k}^{(0)} ;r_{k})$.
			\EndFor 
			{\color{gray}\LineComment{----- Update best $\epsilon$-feasible solution -----}}
			 \State If $\bx'_{k'}$ is $\epsilon$-feasible and $f_0(\bx'_{k'}) < f_0(\bx_{\text{best}})$ then $\bx_{\text{best}} = \bx'_{k'}$. 
			\EndFor
			\State{{\bf Output:} $\bx_{\text{best}}$}.
		\end{algorithmic}
	\end{algorithm}

By Lemma~\ref{lem:lssequence}, the sequence $\{(\bx_k^{(0)},r_k)\}_{k=0}^K$ generated  at each iteration of  Algorithm~\ref{algo:RestartingLS} is a level set sequence. In addition, since $K=\tilde K\geq\hat K$,  Theorem~\ref{thm:keythm} guarantees that the critical index $k^*$ must exist unless an $\epsilon$-optimal and $\epsilon$-feasible solution is obtained. The correctness and complexity of Algorithm~\ref{algo:RestartingLS} rely on the following lemma, which we will prove in the Appendix.
\begin{lemma}\label{lem:CInondecreasing}
The critical index $k^*$  is non-decreasing throughout Algorithm~\ref{algo:RestartingLS}.\end{lemma}

 We count the total number of $\mathtt{fom}$ iterations taken by the RLS algorithm as the total number of iterations taken by $\mathtt{fom}$ instances between two consecutive desirable restarts times the total number of desirable restarts. Lemma~\ref{lem:CInondecreasing} shows the critical index never decreases after desirable restarts. Since the critical index $k^*$ is at most $\tilde K$ (by its definition), this index varies between one to $\tilde K$ in an increasing order until it reaches $\tilde K$. If the latter case happens, Assumption~\ref{assum:FOM} and Theorem~\ref{thm:keythm} ensure that in finite number of iterations condition~\eqref{eq:xkquality} holds at all $\mathtt{fom}$ instances and hence an $\epsilon$-optimal and $\epsilon$-feasible solution must be found. Let $D$ denote the total number of desirable restarts before an $\epsilon$-optimal and $\epsilon$-feasible solution is found. Using the argument above,  we can compute  $D$ as follows. Suppose $D_{k,k^*}$ shows the number of desirable restarts starting at $\mathtt{fom}_k$ before an $\epsilon$-optimal and $\epsilon$-feasible solution is found with $k\leq k^*$, assuming $k^*$ is the current critical index. Then
\begin{eqnarray}\label{eq:TotalDRcount}
D &:= \sum_{k^* = 0}^{\tilde K} D_{k^*,k^*} + \sum_{k^* = 0}^{\tilde K}\sum_{k = 0}^{k^*-1} D_{k,k^*}.
\end{eqnarray}

The first sum in the above equation accounts for the total number of desirable restarts that occur at each $\mathtt{fom}_{k^*}$ for all possible critical index $k^* = 0, 1, \ldots, \tilde K$. The second term however is the total number of desirable restarts that start at each $\mathtt{fom}_k$ instance with $k < k^*$, where $k^*$ is the critical index and varies from zero to $\tilde K$. Notice that some of the values between zero and $\tilde K$ in the above sums cannot be critical index. This could happen if the critical index is strictly positive at the beginning of Algorithm~\ref{algo:RestartingLS}. In this case, all the smaller indices will never be a critical index since by Lemma~\ref{lem:CInondecreasing}, critical index does not decrease throughout the algorithm. We assume $D_{kk^*} = D_{k^*k^*} = 0$ for such cases. The equation~\eqref{eq:TotalDRcount} can be re-written as 
\begin{eqnarray}\label{eq:TotalDRcount_rearranged}
D &:= \sum_{k^* = 0}^{\tilde K} D_{k^*,k^*} + \sum_{k = 0}^{\tilde K-1}\sum_{k^* = k+1}^{\tilde K} D_{k,k^*}.
\end{eqnarray}
The equivalence between the second sum in~\eqref{eq:TotalDRcount} and the second one in~\eqref{eq:TotalDRcount_rearranged} can be explained by considering each $D_{k,k^*}$ as an element of an upper-triangle matrix where $k$ is a row index and $k^*$ is a column index. Then the sums of elements over rows and columns can be exchanged. 

We provide upper bounds for $D_{k^*,k^*}$ and $\sum_{k^* = k+1}^{\tilde K} D_{k,k^*}$ respectively in Lemma~\ref{lem:BoundDk'k'} and Lemma~\ref{lem:BoundDkk'} in the Appendix. Using these bounds in~\eqref{eq:TotalDRcount_rearranged}, we obtain the total complexity of Algorithm \ref{algo:RestartingLS}
in Theorem~\ref{thm:correctnessAndComplexity}.

\begin{theorem}\label{thm:correctnessAndComplexity}
Consider $\alpha$, $\beta$, $\epsilon$, and  $r_{\text{ini}}$ such that $0<\alpha<B<1$, $\epsilon > 0$, and $r_{\text{ini}} < f^*$. Suppose Algorithm \ref{algo:RestartingLS} is executed for $K = \tilde K$ where $\tilde K$ is defined in \eqref{eq:tildeK} with $r=r_{\text{ini}}$. This algorithm returns an $\epsilon$-optimal and $\epsilon$-feasible solution in at most 
\begin{equation}\label{eq:total.complexity}
(\tilde K+1)\cdot n^{\mathtt{fom}}\cdot \tilde D = \mathcal{O}\left( n^{\mathtt{fom}}\log^2\left(\frac{1}{\epsilon}\right)\right),
\end{equation}
$\mathtt{fom}$ iterations where $n^\mathtt{fom}$ is as in Assumption~\ref{assum:FOM} and $\tilde D$ is an upper bound on $D$ that can be written as\looseness = -1
{\small
\begin{eqnarray}\label{eq:Dtilde}
\tilde D :=  \left(\tilde K+1\right)\ln(2/\alpha\theta)/\ln(1/B)  
+  \tilde K \cdot \ln\left(4-\alpha\theta/3\alpha\theta\right)/\ln\left(1+\alpha(1-B)\theta/4\right)
= \mathcal{O}\left(\log\left(\frac{1}{\epsilon}\right)\right).
\end{eqnarray}
}
\end{theorem}

\proof{}
Suppose $C$ and $D$ respectively denote the number of $\mathtt{fom}$ iterations between two consecutive desirable restarts and the total number of desirable restarts performed by Algorithm~\ref{algo:RestartingLS}  until an $\epsilon$-optimal and $\epsilon$-feasible solution is found. It then follows that the total number of $\mathtt{fom}$ iterations required by Algorithm~\ref{algo:RestartingLS} is $C\cdot D$, where $D$ is expressed in~\eqref{eq:TotalDRcount_rearranged}. To obtain the iteration complexity, we only require to bound $D$ and $C$. 

\emph{An upper bound on $C$:} Consider a desirable restart and assume the critical index is $k^*$ after this restart. Suppose an $\epsilon$-optimal and $\epsilon$-feasible solution has not been found. We must have $P(x_{k^*}^{(0)}; r_{k^*})>\epsilon$ by Lemma~\ref{lemma:solutionquality}. Assumption~\ref{assum:FOM} indicates that in $t_{k^*} \leq n^{\mathtt{fom}}$ iterations, the inequality $P(x_{k^*}^{(t_{k^*})}; r_{k^*})\leq BP(x_{k^*}^{(0)}; r_{k^*})$ holds. This means that the next desirable restart must occur at one of the $\tilde K+1$ $\mathtt{fom}$ instances in at most $n^{\mathtt{fom}}$ iterations, and $n^{\mathtt{fom}}$ is independent of $\bx^{(0)}$ and $r$ by Assumption \ref{assum:FOM}. Since each iteration of Algorithm~\ref{algo:RestartingLS} simultaneously runs $\tilde K+1$ $\mathtt{fom}$ iterations, we get
\begin{equation}\label{eq:boundonC}
C \leq (\tilde K+1)\cdot n^{\mathtt{fom}}.
\end{equation}

\emph{An upper bound on $D$:} We claim that for any $k=0,1,\dots, \tilde K-1$,
	\begin{align}\label{eq:Dkk'bound}
	\sum_{k^* = k+1}^{\tilde K} D_{k,k^*} \leq \ln\left(4-\alpha\theta/3\alpha\theta\right)/\ln\left(1+\alpha(1-B)\theta/4\right)
	\end{align}
	and for any $k^*=0,1,\dots, \tilde K$,
	\begin{align}\label{eq:Dk'k'bound}
	D_{k^*,k^*}\leq \ln(2/\alpha\theta)/\ln(1/B).
	\end{align} 
	We formally prove these claims in Lemma~\ref{lem:BoundDkk'} and Lemma~\ref{lem:BoundDk'k'} in the Appendix.
From the above inequalities, it follows that the total number of desirable restarts can be bounded by 
	\begin{eqnarray}\label{eq:boundonD}
	 D && = \sum_{k^* = 0}^{\tilde K} D_{k^*,k^*} + \sum_{k = 0}^{\tilde K-1}\sum_{k^* = k+1}^{\tilde K} D_{k,k^*} \nonumber \\[10pt]
	  && \leq \left(\tilde K+1\right)\ln(2/\alpha\theta)/\ln(1/B)  
+  \tilde K \cdot \ln\left(4-\alpha\theta/3\alpha\theta\right)/\ln\left(1+\alpha(1-B)\theta/4\right) = \tilde D.\nonumber\\
	\end{eqnarray}	
	The bound~\eqref{eq:total.complexity} then can be obtained by multiplying \eqref{eq:boundonC} and \eqref{eq:boundonD}.
\hfill\BlackBox

The key remaining question to implement Algorithm \ref{algo:RestartingLS} is the following: What $\mathtt{fom}$ satisfies Assumption \ref{assum:FOM} and what is $n^{\mathtt{fom}}$? In the following sections, we will present different $\mathtt{fom}$s that satisfy Assumption \ref{assum:FOM} for the non-smooth and smooth problems and provide a total complexity bound for Algorithm~\ref{algo:RestartingLS} under each case. 
	
\section{First Order Subroutines}\label{sec:foms} In this section, we provide two different first order methods ($\mathtt{fom}$s) for smooth and non-smooth problems that can be used as subroutines in Algorithm \ref{algo:RestartingLS} to solve $\min_{\bx\in\X}\Pc(\bx; r)$. In particular, when the functions $f_i$, $i = 0,1,\ldots,m$ defining $\Pc(\bx; r)$ are non-smooth, we use the standard subgradient method as an $\mathtt{fom}$, which is an optimal algorithm for solving general non-smooth problems. When the functions $f_i$, $i = 0,1,\ldots,m$ are smooth, we still need to solve the minimization of a non-smooth function since $\Pc(\bx; r)$ is non-smooth because of the $\max\{\cdot\}$ operator in its definition. In this case, 
we discuss the accelerated proximal linear method~\citep{CompositionsofConvexandSmooth}, which is an extension of proximal gradient algorithm by~\citet{beck2009fast} and~\citet{nesterov_nonsmoothminimization}. This algorithm harnesses the inherent smoothness of the functions enclosed within the max operator and yields a better convergence result.

Since all $\mathtt{fom}_k$, $k = 0,1, \ldots,\tilde K$ follow the same steps, we drop the index $k$ in $\bx_k^{(t_k)}$ and $r_k$ used in the first order methods described in this section to simplify our notations. 
	
\subsection{Non-Smooth Case}\label{nonsmoothcase}
First, let's assume the functions $f_i$, $i = 0,1,\ldots,m$ are non-smooth and $\partial f_i(\bx)$ denotes the set of subgradient of $f_i$ at $\bx$. Suppose   Assumption~\ref{assume:strictfeasible} and the error bound condition	\eqref{eq:errorbound} hold. 
We next show that the subgradient method using a specific step length rule satisfies Assumption \ref{assum:FOM} and can be used to solve the level-set subproblem $\min\limits_{\bx\in\mathcal{X}}\Pc(\bx;r)$.  
 
Let $\bx^{(0)}\in\X$ be an initial solution, the subgradient updates can be presented as  
	\begin{equation}
	\label{eq:sgd}
	\bx^{(t+1)} = \text{Proj}_{\mathcal{X}}\left(\bx^{(t)} - \eta^{(t)} \bxi_{r}^{(t)}\right),\quad t=0,1,\dots,
	\end{equation}
where $\text{Proj}_{\mathcal{X}}(\cdot)$ denotes the projection onto $\X$, $\eta^{(t)}>0$ a step size, and $\bxi_{r}^{(t)} \in \partial P(\bx^{(t)};r)$ a subgradient of $P(\bx^{(t)};r)$ with respect to $\bx^{(t)}$. Recall that the projection mapping $\text{Proj}_{\mathcal{X}}(\cdot)$ is easily computable since we assume the set $\mathcal{X}$ is simple (e.g., $\mathbb{R}^n$, a box, or a ball). The output of the subgradient method can be chosen as the historically best iterate, i.e., 
\begin{equation}
	\label{eq:sgdoutput}
	\bar\bx^{(t)}:=\argmin_{s = 0,1,\ldots,t} P(\bx^{(s)};r).
\end{equation}
Proposition~\ref{prop:subgradient_method_convergence} below presents a well-known convergence result of the subgradient method.
	
\begin{proposition}[See Theorem 3.2.2 in \citet{nesterov2018lectures}]\label{prop:subgradient_method_convergence}
Consider $r < f^*$. Let  $\bx^*=\text{Proj}_{\X^*}(\bx^{(0)})$ and $\bar\bx^{(t)}$ be defined as in \eqref{eq:sgdoutput}. The subgradient method in \eqref{eq:sgd} guarantees that for any $t\geq 0$, 
		\begin{align}
		P (\bar\bx^{(t)};r) -  P (\bx^*;r) \leq \frac{\text{dist}(\bx^{(0)}, \X^*)^2 + \sum_{s=0}^{t}(\eta^{(s)})^2 \|\bxi_{r}^{(s)}\|^2 }{2\sum_{s=0}^{t} \eta^{(s)}}. \label{eq:SGDResult}
		\end{align}
	\end{proposition}
Notice that $\bx^*$ in the above proposition could be the projection of any $\bx^{(t)}$, for $t\geq 0$, onto $\mathcal{X}^*$ since by definition of $\mathcal{X}^*$ and $\Pc(\bx;r)$, it follows that $\Pc(\text{Proj}_{\X^*}(\bx^{(t)}); r) = f^* - r$ for all $t\geq 0$.\looseness = -1
	
We will shortly show that the subgradient method discussed above can be used as an $\mathtt{fom}$ in Algorithm~\ref{algo:RestartingLS} if the stepsizes are carefully chosen. However, before proceeding with that, we make the following assumptions to ensure the rate of change of the function $\Pc(\cdot; r)$ does not become arbitrarily steep at the points generated by the subgradient algorithm.
\begin{assumption}\label{assumptions:nonsmooth}
 Consider the initial solution $\bx^{(0)}\in\X$, there exists $M$ such that $\max\limits_{\bxi \in\partial f_i(\bx)} \| \bxi\| \leq M$ for $i=0,1,\dots,m$ and any $\bx\in\X$ that satisfies $\text{dist}(\bx, \mathcal{X}^*)\leq\text{dist}(\bx^{(0)}, \mathcal{X}^*)$.
\end{assumption}
Notice that this assumption is more lenient than the standard assumption of bounded subgradients across the entire domain which is commonly used in the literature. This standard assumption implies that $\text{dist}(\bx, \mathcal{X}^*)$ is bounded for any $\bx\in\X$ when $d > 1$.  In contrast, we only require subgradients to be bounded when we are in proximity to the set of optimal solutions. 

Proposition~\ref{prop:SGD_nonexpansion} below shows that the inequality $\text{dist}(\bx^{(t)}, \mathcal{X}^*)\leq\text{dist}(\bx^{(0)}, \mathcal{X}^*)$ holds at any solution $\bx^{(t)}$ generated by the subgradient algorithm. Therefore, under Assumption~\ref{assumptions:nonsmooth}, it easily follows that $\left\|\bxi_{r}^{(t)}\right\|\leq M$ for any $\bxi_{r}^{(t)}\in \partial P(\bx^{(t)};r)$ and any $r < f^*$.

\begin{proposition}\label{prop:SGD_nonexpansion}
Consider $\alpha$, $B$, $\epsilon$, $r$, and $\bx^{(0)}\in\mathcal{X}$ such that $0<\alpha<B<1$, $\epsilon>0$, $r < f^*$, and $\Pc(\bx^{(0)}; r)> \epsilon$. Let $\eta^{(t)} := \frac{(B-\alpha)  P( \bx^{(0)}; r)  }{\|\bxi_{r}^{(t)}\|^2}$. Suppose $\bx^{(0)}$ does not satisfy~\eqref{eq:xkquality}, i.e. $\alpha\Pc(\bx^{(0)};r) > f^*-r$. At any iteration $t$ of the subgradient method, we have
$$
\text{dist}(\bx^{(t)}, \mathcal{X}^*)\leq \text{dist}(\bx^{(0)}, \mathcal{X}^*),
$$
until a desirable solution is found. i.e. $\Pc(\bar\bx^{(t)};r) \leq BP (\bx^{(0)};r)$.
	\end{proposition}
\proof{}
	Consider $t\geq 1$ and suppose $\Pc(\bar\bx^{(t)};r) > BP (\bx^{(0)};r)$ which indicates that $\Pc(\bx^{(s)};r) > BP (\bx^{(0)};r)$ for $s=0,\dots,t$. Let $\bx^* = \text{Proj}_{\X^*}(\bx^{(t-1)})\in\X^*$.  Since $r< f^*$, from definition of $P(\bx;r)$ it follows that 
\begin{equation}\label{eq:pxstarR}
		P (\bx^*;r) =f^*-r. 
		\end{equation}
Hence, for any $s=0,\dots,t$, we have
	\begin{equation}\label{eq:pxstarRLessthanPxt}
		P (\bx^*;r) =f^*-r<\alpha\Pc(\bx^{(0)};r) <  \frac{\alpha}{B} \Pc(\bx^{(s)};r). 
	\end{equation}
	By the standard analysis of subgradient methods, we have
	\begin{eqnarray}
		\nonumber
		\text{dist}(\bx^{(t)}, \mathcal{X}^*)^2&\leq&\text{dist}(\bx^{(t-1)}, \mathcal{X}^*)^2-2\eta^{(t-1)}\left\langle \bxi_{r}^{(t-1)}, \bx^{(t-1)}-\bx^*\right\rangle+\left(\eta^{(t-1)}\right)^2\left\|\bxi_{r}^{(t-1)}\right\|^2\\\nonumber
		&\leq&\text{dist}(\bx^{(t-1)}, \mathcal{X}^*)^2-2\eta^{(t-1)}\left(\Pc (\bx^{(t-1)};r) - \Pc (\bx^*;r)\right)+\left(\eta^{(t-1)}\right)^2\left\|\bxi_{r}^{(t-1)}\right\|^2 \\\nonumber
		&\leq&\text{dist}(\bx^{(t-1)}, \mathcal{X}^*)^2-2\eta^{(t-1)}(B-\alpha)\Pc (\bx^{(0)};r) +\eta^{(t-1)} (B-\alpha) \Pc( \bx^{(0)}; r)\\\nonumber
		&=&\text{dist}(\bx^{(t-1)}, \mathcal{X}^*)^2-\eta^{(t-1)}\left(B-\alpha\right)\Pc( \bx^{(0)}; r)\\\nonumber
		&\leq&\text{dist}(\bx^{(t-1)}, \mathcal{X}^*)^2-\left(B-\alpha\right)^2\frac{\Pc( \bx^{(0)}; r)^2}{\|\bxi_{r}^{(t-1)}\|^2},
	\end{eqnarray}
	where the second inequality holds by the convexity of $\Pc(\bx;r)$ at $\bx$ and the third inequality by \eqref{eq:pxstarRLessthanPxt} and the definition of $\eta^{(t)}$. Since $\left(B-\alpha\right)^2\frac{\Pc( \bx^{(0)}; r)^2}{\|\bxi_{r}^{(t-1)}\|^2} \geq 0$, the above inequality indicates that
	$
	\text{dist}(\bx^{(t)}, \mathcal{X}^*)^2\leq \text{dist}(\bx^{(t-1)}, \mathcal{X}^*)^2 \leq \ldots \leq \text{dist}(\bx^{(0)}, \mathcal{X}^*)^2.
	$ \hfill\BlackBox
	
Proposition~\ref{prop:SGD.complexity} indicates that the subgradient method with a specific choice of stepsize satisfies Assumption~\ref{assum:FOM} and therefore it can be used as an $\mathtt{fom}$ in Algorithm~\ref{algo:RestartingLS}. Notice that some of the arguments used in the proof of Proposition~\ref{prop:SGD.complexity} are borrowed from the proof of Proposition 3.2 in \citet{freund2018new}.
\begin{proposition}
		\label{prop:SGD.complexity}
		Consider $\alpha$, $B$, $\epsilon$, $r$, and $\bx^{(0)}\in\mathcal{X}$ such that $0<\alpha<B<1$, $\epsilon>0$, $r < f^*$, and $\Pc(\bx^{(0)}; r)> \epsilon$. Let $\eta^{(t)} := \frac{(B-\alpha)  P( \bx^{(0)}; r)  }{\|\bxi_{r}^{(t)}\|^2}$. The subgradient method satisfies Assumption~\ref{assum:FOM} with
        \begin{eqnarray}\label{eq:peritercost}
        	 n^{\mathtt{fom}}= \left\lceil\frac{ M^2 G^{2/d}}{(B-\alpha)^2\epsilon^{2-2/d}}\right\rceil -1.
        \end{eqnarray} 
\end{proposition}
\proof{}
	Suppose $\bx^{(0)}$ does not satisfy~\eqref{eq:xkquality}, i.e. $\alpha\Pc(\bx^{(0)};r) > f^*-r$. We show $\Pc(\bar\bx^{(t)};r) \leq BP (\bx^{(0)};r)$ for $t\leq n^{\mathtt{fom}}$.
		
Let $\bx^* = \text{Proj}_{\X^*}(\bx^{(0)})\in\X^*$. Plugging the definition of $\eta^{(t)}$ into~\eqref{eq:SGDResult} we get
		\begin{eqnarray}
		P (\bar\bx^{(t)};r) &\leq& P (\bx^*;r) + \dfrac{1}{2\sum_{s=0}^{t} \eta^{(s)}}\left[\text{dist}(\bx^{(0)}, \X^*)^2 + \sum_{s=0}^{t}\eta^{(s)} \dfrac{(B-\alpha)  P( \bx^{(0)}; r)  }{\|\bxi_{r}^{(t)}\|^2} \|\bxi_{r}^{(s)}\|^2 \right]\nonumber \\
&\leq&
P (\bx^*;r)  + \dfrac{\text{dist}(\bx^{(0)}, \mathcal{X}^*)^2}{2(B-\alpha)  P( \bx^{(0)};r ) \sum_{s=0}^{t}\|\bxi_{r}^{(s)}\|^{-2}} + \dfrac{B-\alpha}{2}P (\bx^{(0)};r)
		\nonumber \\
		&\leq& f^*-r  + \dfrac{M^2\text{dist}(\bx^{(0)}, \mathcal{X}^*)^2}{2(t+1)(B-\alpha) P (\bx^{(0)};r)} + \dfrac{B-\alpha}{2}P (\bx^{(0)};r)\nonumber \\
		&\leq& f^*-r + \dfrac{M^2 G^{2/d}  P (\bx^{(0)};f^*)^{2/d}}{2(t+1)(B-\alpha)P (\bx^{(0)};r)}+ \dfrac{B-\alpha}{2}P (\bx^{(0)};r) \nonumber\\
		&\leq& f^*-r  + \dfrac{M^2 G^{2/d}  P (\bx^{(0)};r)^{2/d-1}}{2(t+1)(B-\alpha)}+ \dfrac{B-\alpha}{2}P (\bx^{(0)};r) \label{eq:ineq1},
		\end{eqnarray}
where the third inequality follows from \eqref{eq:pxstarR} and the inequality $\|\bxi_{r}^{(t)}\|\leq M$ which holds by Assumption~\ref{assumptions:nonsmooth}. The fourth inequality follows from the error bound condition inequality	\eqref{eq:errorbound} and definition of $P(\bx; r)$, and the last from the inequality $P(\bx^{(0)};f^*)\leq  P(\bx^{(0)};r)$ which holds because $r\leq f^*$. Let 	
\begin{eqnarray}
\label{eq:tkbound}
t =\left\lceil\frac{ M^2 G^{2/d} }{(B-\alpha)^2 P (\bx^{(0)};r)^{2-2/d}}\right\rceil-1.
\end{eqnarray}
Using \eqref{eq:tkbound} in \eqref{eq:ineq1}, we get
$
		P (\bar\bx^{(t)};r) \leq f^*-r  +\frac{B-\alpha}{2}P (\bx^{(0)};r) + \frac{B-\alpha}{2}P (\bx^{(0)};r) 
		\leq B P (\bar\bx^{(0)};r), 
		$
		where the second inequality is by the assumption that $\alpha P (\bx^{(0)};r)> f^*-r$. Since $\Pc(\bx^{(0)};r)> \epsilon$, we have $t\leq n^{\mathtt{fom}}$.
	\hfill\BlackBox

	\vspace{0.1in}
    \begin{remark}\label{rem:stepLengthFOM}
    The step length $\eta = (B-\alpha)\Pc(\bx^{(0)}; r)/\|\bxi_{r}\|^2$ used in subgradient algorithm is likely to reduce at the instance that is restarted since $\Pc(\bx^{(0)}; r)$ becomes smaller as a result of updating the initial solution with a solution that satisfies \eqref{eq:reduceB}. For all the subsequent instances, the step length after the restart would likely be larger since $\Pc(\bx^{(0)}; r)$ is a decreasing function of $r$. 
	\end{remark}   

Corollary~\ref{cor:non.smooth.complexity} below provides the overall complexity of Algorithm~\ref{algo:RestartingLS} by counting the total number of subgradient and function evaluations required in this algorithm.
	\begin{corollary}
		\label{cor:non.smooth.complexity}
		Suppose $\alpha$, $\beta$, $\epsilon$, and $r_{\text{ini}}$ are such that $0<\alpha<B<1$, $\epsilon > 0$, and $r_{\text{ini}} < f^*$. Consider $\tilde K$ defined in \eqref{eq:tildeK} for $r = r_{\text{ini}}$. Suppose we execute Algorithm~\ref{algo:RestartingLS} with subgradient method as $\mathtt{fom}_k$ instances for $k =0, 1,\ldots, \tilde K$. This algorithm finds an $\epsilon$-optimal and $\epsilon$-feasible solution in at most 
		$$
		(m+1)\cdot (\tilde K+1) \cdot \left\lceil\frac{ M^2 G^{2/d} }{(B-\alpha)^2 \epsilon^{2-2/d}}\right\rceil \cdot \tilde D =\mathcal{O}\left(\log^2\left(\frac{1}{\epsilon}\right)\frac{1}{\epsilon^{2-2/d}}\right),
		$$
		subgradient or function evaluations, where $\tilde D$ is defined in \eqref{eq:Dtilde}.
		
	\end{corollary}
\proof{} Theorem~\ref{thm:correctnessAndComplexity} guarantees that an $\epsilon$-feasible and $\epsilon$-optimal solution can be found in at most\looseness=-1
	\begin{equation}\label{eq:1}
	(\tilde K+1)\cdot n^{\mathtt{fom}}\cdot \tilde D,
	\end{equation}
	iterations where $\tilde K$, $\tilde D$, and $n^{\mathtt{fom}}$ are respectively defined in \eqref{eq:tildeK}, \eqref{eq:Dtilde}, and \eqref{eq:peritercost}.
	
It is easy to verify that the subgradient algorithm requires $m+1$ subgradients and function evaluations (one function evaluation for each function in $\Pc(\bx;r)$). Hence, the expression \eqref{eq:1} indicates that the total number of subgradient computations required in Algorithm~\ref{algo:RestartingLS} is
\begin{equation}\label{eq:no.subgradient}
	(m+ 1) \cdot (\tilde K+1) \cdot n^{\mathtt{fom}} \cdot \tilde D.
\end{equation}
The proof is then completed using \eqref{eq:peritercost} in \eqref{eq:no.subgradient}.\hfill\BlackBox

%
\subsection{Smooth Case}\label{sec:smooth}
In this section, we assume Assumption~\ref{assume:strictfeasible}  and the error bound condition~\eqref{eq:errorbound} hold. Moreover, we make an additional assumption about the smoothness of the functions $f_i$, $i = 0,1, \ldots,m$ in~\eqref{eq:gco}.\looseness = -1
	\begin{assumption}
		\label{assumption:smoothness}
		Functions $f_i$, $i = 0,1,\ldots,m$, are $L$-smooth on $\X$ for some $L\geq0$. In other words, the functions $f_i$, $i=0,1,\ldots,m$, are differentiable and  
		$f_i(\bx)\leq f_i(\by)+\left\langle\nabla f_i(y),\bx-\by\right\rangle+\dfrac{L}{2}\|\bx-\by\|^2$ for any $\bx$ and $\by$ in $\X$.
	\end{assumption}
	
It is important to note that while Assumption 2 guarantees the smoothness of the functions $f_i$, the function $\Pc(\bx;r)$ may not possess this property due to the presence of the $\max\lbrace\cdot\rbrace$ operator in its definition. To address this concern, we use the prox-linear method introduced in~\citet{CompositionsofConvexandSmooth}. In particular, we define a composite linear approximation of $\Pc(\bx;r)$ at a specific point $\by\in\mathcal{X}$ that is:
\begin{eqnarray}
	\label{eq:P_proxlinear}
\Pc(\bx;\by,r)&:=&\max_{i=1,\dots,m}\left\{f_0(\by)+\left\langle \bx-\by, \nabla f_0(\by)\right\rangle-r,f_i(\by)+\left\langle \bx-\by, \nabla f_i(\by)\right\rangle\right\}.
\end{eqnarray}

Given the $L$-smoothness of the functions $f_i$, $i = 0,1,\ldots,m$ in Assumption~\ref{assumption:smoothness} and the convexity of these functions in $\bx$, it is straightforward to verify that for any $\bx,\by\in\X$,
\begin{equation}\label{eq:PlessthanPy}
\Pc(\bx; \by; r)\leq \Pc(\bx; r) \leq \Pc(\bx; \by; r) + \dfrac{L}{2} \|\bx - \by\|^2.
\end{equation}

This inequality implies that to solve $\min_{\bx\in\X} \Pc(\bx; r)$, one can apply a majorization-minimization procedure known as the \emph{prox-linear method}, which solves the following subproblem at each iteration:
\begin{equation}\label{prox-linear:subproblem}\bx^{(t+1)}\in\argmin_{\bx\in\Xc}\Pc(\bx;\bx^{(t)},r)+\dfrac{1}{2\eta}\|\bx-\bx^{(t)}\|^2.\end{equation}
Here, $\eta >0 $ is a step size parameter (typically chosen as $\eta = \frac{1}{L}$). We assume that this strongly convex subproblem can be efficiently solved. There are various methods available to tackle this problem; for instance, one can employ the simplex method \citep{wolfe1959simplex} or utilize the algorithm outlined in \citet{chambolle2011first}.

While the basic prox-linear method in~\eqref{prox-linear:subproblem} can be used to minimize \( \Pc(\bx; r) \), it does not achieve the optimal convergence rate. A faster alternative is the \emph{accelerated prox-linear method}, which introduces two auxiliary sequences \( \by^{(t)} \) and \( \bz^{(t)} \), in addition to the primary sequence \( \bx^{(t)} \).

Starting from an initial point \( \bx^{(0)} \in \mathcal{X} \) and setting \( \bz^{(0)} = \bx^{(0)} \), the algorithm proceeds in iterations. At each iteration \( t \), it computes an extrapolated point:
$
\by^{(t)} = a^{(t)} \bz^{(t)} + (1 - a^{(t)}) \bx^{(t)},
$
where $ a^{(t)} \in (0,1]$ is a weighting parameter. The next iterate $\bx^{(t+1)}$ is then obtained by solving the prox-linear subproblem:
\begin{equation}
\label{acc-prox-linear:subproblem}
\bx^{(t+1)} \in \argmin_{\bx \in \mathcal{X}} \Pc(\bx; \by^{(t)}, r) + \dfrac{1}{2\eta} \|\bx - \by^{(t)}\|^2,
\end{equation}
where $\eta = 1/L$. Finally, the auxiliary sequence is updated as
$
\bz^{(t+1)} = \bx^{(t)} + \frac{1}{a^{(t)}} (\bx^{(t+1)} - \bx^{(t)}),
$
enabling the method to exploit acceleration.

In practice, the step size \( \eta = \frac{1}{L} \) may not be known a priori or may be conservatively estimated. A standard remedy is to incorporate a backtracking line search~\citep{beck2009fast,nesterov_nonsmoothminimization}, as shown in Algorithm~\ref{alg:backtracking}, which dynamically adjusts \( \eta \). The procedure starts with a guess $ \eta^{(t)}$ for $1/L$, solves~\eqref{acc-prox-linear:subproblem} for $\eta = \eta^{(t)}$, and checks whether the inequality
\[
\Pc(\bx^{(t+1)}; r) \leq \Pc(\bx^{(t+1)}; \by^{(t)}, r) + \dfrac{1}{2\eta^{(t)}} \|\bx^{(t+1)} - \by^{(t)}\|^2,
\]
holds. By inequality~\eqref{eq:PlessthanPy}, this condition is satisfied whenever \( \eta^{(t)} \leq \frac{1}{L} \). If the condition fails, then \( \eta^{(t)} > \frac{1}{L} \), and the estimate is decreased by a factor \( \beta_{\text{dec}} < 1 \). The subproblem is resolved with the new \( \eta^{(t)} \), and the process continues until the condition is met.

Once the line search is successful, $ \bx^{(t+1)}$ and the corresponding $\eta^{(t)}$ are returned. Additionally, to allow for potentially larger step sizes in future iterations, the step size is increased by a factor $\beta_{\text{inc}} > 1 $ and used as the initial estimate in the next iteration. This is implemented in Line 3 of Algorithm~\ref{algorithm:accelerated_prox_linear}.

\begin{algorithm}[ht]
\caption{BacktrackingLineSearch}
\label{alg:backtracking}
\begin{algorithmic}[1]
\Require Outer iterate $\bx^{(t)}$, extrapolated point $\bz^{(t)}$, previous weight $a^{(t-1)}$, level parameter $r < f^*$, previous step size $\eta^{(t-1)}$, initial trial step size $\bar \eta$, backtracking parameter $\beta_{\text{dec}} \in (0,1)$.
\State \textbf{Initialize:} $\eta \leftarrow \beta_{\text{dec}}^{-1}\bar\eta$ 
\Repeat
    \State $\eta \leftarrow \beta_{\text{dec}} \cdot \eta$
    \State Compute $a$ by solving:
    \[
    \frac{(a^{(t-1)})^2}{\eta^{(t-1)}} \left(1 - a\right) = \frac{a^2}{\eta}
    \]
    \State Compute extrapolated point:
    \[
    \by = a \cdot \bz^{(t)} + (1 - a) \cdot \bx^{(t)}
    \]
    \State Compute candidate update:
    \[
    \bx^{+} = \arg\min_{\bx \in \Xc} \left\{ \Pc(\bx; \by, r) + \frac{1}{2 \eta} \left\|\bx - \by\right\|^2 \right\}
    \]
\Until{$P(\bx^{+}; r) \leq P(\bx^{+}; \by, r) + \dfrac{1}{2 \eta} \left\|\bx^{+} - \by\right\|^2$}
\State       \Return $(\bx^{+}, \by, \eta, a)$
\end{algorithmic}
\end{algorithm}

\begin{algorithm}[ht]
\caption{Accelerated Prox-Linear Method with Backtracking Line Search}
\label{algorithm:accelerated_prox_linear}
\begin{algorithmic}[1]
\Require Level parameter $r\le f^*$, initial point $\bx^{(0)}\in X$, initial step size $\eta^{(\text{ini})}>0$, and backtracking parameters $\beta_{\text{dec}}\in (0,1)$ and $\beta_{\text{inc}}>1$.
\State \textbf{Initialize}: set $\bz^{(0)}\leftarrow  \bx^{(0)}$, $a^{(-1)} \leftarrow 1$ and $\eta^{(-1)} \leftarrow \eta^{(\text{ini})}$.
\For{$t=0,1,2,\dots$}
   \State Set $\bar\eta \leftarrow\min\{\beta_{\text{inc}}\,\eta^{(t-1)},\eta^{(\text{ini})}\}$.
    \State $(\bx^{(t+1)},\by^{(t)},\,\,\eta^{(t)},\,a^{(t)}) \leftarrow$ \textbf{BacktrackingLineSearch}$(\bx^{(t)},\bz^{(t)},a^{(t-1)}, r, \eta^{(t-1)},\bar\eta, \beta_{\text{dec}})$.
    \State Update the extrapolated solution:
    \[
    \bz^{(t+1)} = \bx^{(t)} + \frac{1}{a^{(t)}}\Bigl(\bx^{(t+1)} - \bx^{(t)}\Bigr).
    \]
\EndFor
\State \Return $\bx^{(t+1)}$.
\end{algorithmic}
\end{algorithm}

\begin{theorem}
Given $\eta^{(\text{ini})}>0$, $\beta_\text{inc} > 1$, $\beta_{\text{dec}}\in (0,1)$, and $r < f^*$, Algorithm~\ref{algorithm:accelerated_prox_linear} guarantees that for any $\bx \in \X$ and $t \geq 0$,
\begin{align}\label{ineq:approx_linear}
\Pc(\bx^{(t)};r) - \Pc(\bx;r) &\leq \frac{4\max\left\{\frac{1}{\eta^{(\text{ini})}},\frac{L}{\beta_{\text{dec}}}\right\}}{(t+1)^2} \left\{ \eta^{(\text{ini})}\left[\Pc(\bx^{(0)};r)-\Pc(\bx;r)\right] + \frac{1}{2}\left\|\bx^{(0)} - \bx\right\|^2 \right\}.
\end{align}
Moreover, the total number of prox-linear subproblems solved in Algorithm \ref{alg:backtracking} is at most
$$
\left[1 + \frac{\ln (\beta_{\text{inc}})}{\ln (\beta_{\text{dec}}^{-1})} \right] t + \frac{1}{\ln (\beta_{\text{dec}}^{-1})} \ln \left( \max \left\{ 1, \frac{\eta^{(\text{ini})} L}{\beta_{\text{dec}}} \right\} \right) = \mathcal{O}(t).
$$
\label{theorem:APLM}
\end{theorem}

\proof{}
Let $t \geq 0$. In this proof, we denote by $\eta^{(t)}$, $a^{(t)}$, $\bx^{(t+1)}$, and $\by^{(t)}$ the final step size, weight, solution, and extrapolated point returned by the backtracking line search at the $t$-th outer iteration of Algorithm~\ref{algorithm:accelerated_prox_linear}, i.e., after the internal subroutine, Algorithm~\ref{alg:backtracking}, has terminated.

Since $\by^{(t)} = a^{(t)}\bz^{(t)} + (1 - a^{(t)}) \bx^{(t)}$, we have
$\bz^{(t)} - \by^{(t)} = -\frac{1 - a^{(t)}}{a^{(t)}} (\bx^{(t)} - \by^{(t)}).$ Let $\hat\bx = a^{(t)} \bx + (1 - a^{(t)}) \bx^{(t)}$ for $\bx \in \X$. Then:
$
\hat\bx - \by^{(t)} = a^{(t)}(\bx - \by^{(t)}) + (1 - a^{(t)})(\bx^{(t)} - \by^{(t)}),
$
which implies
\begin{equation}\label{eq:xhatminusy}
    \hat\bx - \by^{(t)} = a^{(t)} (\bx - \bz^{(t)}).
\end{equation} 
Also, using $\bz^{(t+1)} = \bx^{(t)} + \frac{1}{a^{(t)}}(\bx^{(t+1)} - \bx^{(t)})$ and $\hat\bx = a^{(t)}\bx + (1 - a^{(t)})\bx^{(t)}$, it follows that:
\begin{equation}\label{eq:xhatminusxt+1}
    \|\hat\bx - \bx^{(t+1)}\|^2 = (a^{(t)})^2 \|\bx - \bz^{(t+1)}\|^2.
\end{equation}
Next, by the stopping condition in Algorithm~\ref{alg:backtracking}, the definition of $\bx^{(t+1)}$, and strong convexity of the subproblem, we obtain:
\begin{align*}
\Pc(\bx^{(t+1)}; r)
&\leq \Pc(\bx^{(t+1)}; \by^{(t)}, r) + \frac{1}{2\eta^{(t)}} \|\bx^{(t+1)} - \by^{(t)}\|^2 \\
&\leq \Pc(\hat\bx; \by^{(t)}, r) + \frac{1}{2\eta^{(t)}} \|\hat\bx - \by^{(t)}\|^2 - \frac{1}{2\eta^{(t)}} \|\hat\bx - \bx^{(t+1)}\|^2.
\end{align*}
Using convexity of $\Pc(\cdot; \by^{(t)}, r)$, \eqref{eq:xhatminusy}, and \eqref{eq:xhatminusxt+1}, we get
\begin{align*}
\Pc(\bx^{(t+1)}; r)
&\leq (1 - a^{(t)}) \Pc(\bx^{(t)}; \by^{(t)}, r) + a^{(t)} \Pc(\bx; \by^{(t)}, r) \\
&\quad + \frac{(a^{(t)})^2}{2\eta^{(t)}} \|\bx - \bz^{(t)}\|^2 - \frac{(a^{(t)})^2}{2\eta^{(t)}} \|\bx - \bz^{(t+1)}\|^2 \\
&\leq (1 - a^{(t)}) \Pc(\bx^{(t)}; r) + a^{(t)} \Pc(\bx; r) + \frac{(a^{(t)})^2}{2\eta^{(t)}} \|\bx - \bz^{(t)}\|^2 - \frac{(a^{(t)})^2}{2\eta^{(t)}} \|\bx - \bz^{(t+1)}\|^2,
\end{align*}
where the last inequality uses~\eqref{eq:PlessthanPy}.
By adding $\Pc(\bx; r)$ to both sides and rearranging the terms, we obtain
\small
\begin{align*}
\frac{\eta^{(t)}}{(a^{(t)})^2} \left[ \Pc(\bx^{(t+1)}; r) - \Pc(\bx; r) \right]
&\leq \frac{\eta^{(t)}(1 - a^{(t)})}{(a^{(t)})^2} \left[ \Pc(\bx^{(t)}; r) - \Pc(\bx; r) \right]+ \frac{1}{2} \|\bx - \bz^{(t)}\|^2 - \frac{1}{2} \|\bx - \bz^{(t+1)}\|^2.
\end{align*}
\normalsize
In addition, by construction (Algorithm~\ref{alg:backtracking}, Line 4), we have:
\begin{align*}
\frac{\eta^{(t)}}{(a^{(t)})^2} \,\left[\, \Pc(\bx^{(t+1)}; r) \,-\, \Pc(\bx; r)\, \right]
&\leq \dfrac{\eta^{(t-1)}}{(a^{(t-1)})^2} \left[\, \Pc(\bx^{(t)}; r)\, - \,\Pc(\bx; r)\, \right] \\
&\quad + \frac{1}{2}\, \|\,\bx\, - \,\bz^{(t)}\,\|^2 \,-\, \frac{1}{2}\, \|\,\bx\, - \,\bz^{(t+1)}\,\|^2.
\end{align*}
Summing the inequality above over $t = 0, \dots, T - 1$ (with $a^{(-1)} = 1$, $\eta^{(-1)} = \eta^{(\text{ini})}$, and $\bz^{(0)} = \bx^{(0)}$), we get:\looseness = -1
\begin{equation}\label{eq:boundonPdifference}
\Pc(\bx^{(t)}; r) - \Pc(\bx; r)
\,\leq\, \frac{(a^{(t-1)})^2}{\eta^{(t-1)}} \left( \eta^{(\text{ini})}\, \left[\,\Pc(\bx^{(0)}; r)\, -\, \Pc(\bx; r)\,\right] + \frac{1}{2} \,\|\bx - \bx^{(0)}\|^2 \right).
\end{equation}
Since Algorithm~\ref{alg:backtracking} is terminated if $\eta^{(t)} \leq \frac{1}{L}$, it must hold that $\eta^{(t)} \geq \beta_{\text{dec}}/L$ if Algorithm~\ref{alg:backtracking} runs for more than one iteration. If Algorithm~\ref{alg:backtracking} is terminated after one iteration, we have $\eta^{(t)}=\min\{\beta_{\text{inc}}\,\eta^{(t-1)},\eta^{(\text{ini})}\}$ (Algorithm~\ref{algorithm:accelerated_prox_linear}, Line 3). Combining these two cases, we have:
$$
\eta^{(t)} \geq \min \left\{ \eta^{(\text{ini})}, \frac{\beta_{\text{dec}}}{L} \right\}\quad \text{ for } t = -1, 0, 1, \ldots
$$
Also, the recurrence relation
$
\frac{\sqrt{\eta^{(t)}}}{a^{(t)}} = \frac{\sqrt{\eta^{(t)}}}{2} + \frac{1}{2}\sqrt{\eta^{(t)}+\frac{4\eta^{(t-1)}}{(a^{(t-1)})^2}} \geq\frac{\sqrt{\eta^{(t)}}}{2} + \frac{\sqrt{\eta^{(t-1)}}}{a^{(t-1)}}
$ leads to:
$$
\frac{(a^{(t-1)})^2}{\eta^{(t-1)}} \leq \frac{4}{(t+1)^2} \max \left\{ \frac{1}{\eta^{(\text{ini})}}, \frac{L}{\beta_{\text{dec}}} \right\}.
$$
Substituting this into \eqref{eq:boundonPdifference} completes the proof of~\eqref{ineq:approx_linear}. Finally, the total number of prox-linear subproblems solved follows by the same argument as in Lemma 4 of \citet{nesterov_nonsmoothminimization}.
\hfill\BlackBox

We next show Algorithm~\ref{algorithm:accelerated_prox_linear} can be used as an $\mathtt{fom}$ in Algorithm~\ref{algo:RestartingLS}.
\begin{proposition}
	\label{prop:APLM.complexity}
	Consider $\alpha$,  $B$, $\epsilon$, $r$, $\eta^{(\text{ini})}$, $\beta_\text{dec}$ and $\bx^{(0)}\in\mathcal{X}$ such that $0<\alpha<B<1$, $\epsilon>0$, $r < f^*$, $\eta^{(\text{ini})}>0$, $\beta_{\text{dec}}\in (0,1)$, and $\Pc(\bx^{(0)}; r) > \epsilon$.The accelerated prox-linear method described in Algorithm~\ref{algorithm:accelerated_prox_linear}  satisfies Assumption~\ref{assum:FOM}  with
	\begin{eqnarray}
		\label{eq:peritercostsmooth}
		 n^{\mathtt{fom}}:=\max\left\{\sqrt{\frac{2 \max\left\{\frac{1}{\eta^{(\text{ini})}},\frac{L}{\beta_{\text{dec}}}\right\}G^{2/d}}{\alpha(B-\alpha)\epsilon^{1-2/d}}},
\sqrt{\frac{4\max\left\{1,\frac{\eta^{(\text{ini})}L}{\beta_{\text{dec}}}\right\}}{B-\alpha}}
\right\}-1.
	\end{eqnarray}
\end{proposition}	
\proof{}
Suppose $\bx^{(0)}$ does not satisfy~\eqref{eq:xkquality}, i.e. $\alpha\Pc(\bx^{(0)};r)> f^*-r$. We only need to show 	$\Pc(\bar\bx^{(t)};r) \leq BP (\bx^{(0)};r)$ for $t\leq n^{\mathtt{fom}}$.
Let $\bx^* = \text{Proj}_{\X^*}(\bx^{(0)})$, and consider the approximation bound from Theorem~\ref{theorem:APLM} for any $t \geq 1$ and $\bx = \bx^*$: 
\begin{align}
\Pc(\bx^{(t)};r)
&\leq \Pc(\bx^*;r) + \frac{4C_1}{(t+1)^2} \left[ \Pc(\bx^{(0)}; r) - \Pc(\bx^*; r) \right] + \frac{2 C_2 G^{2/d} \cdot \Pc(\bx^{(0)}; f^*)^{2/d}}{(t+1)^2}, \nonumber \\
&\leq f^* - r + \frac{4C_1}{(t+1)^2} \left[ \Pc(\bx^{(0)}; r) - (f^* - r) \right] + \frac{2 C_2 G^{2/d} \cdot \Pc(\bx^{(0)}; r)^{2/d}}{(t+1)^2}, \label{eq:ineq2_apl}
\end{align}
where for simplifying the equations, we defined
$$
C_1 := \max\left\{ 1, \frac{ \eta^{(\text{ini})} L }{ \beta_{\text{dec}} } \right\}, \qquad
C_2 := \max\left\{ \frac{1}{\eta^{(\text{ini})}}, \frac{L}{\beta_{\text{dec}}} \right\}.
$$ The first inequality follows from the error bound condition~\eqref{eq:errorbound} and the last from~\eqref{eq:pxstarR} and the fact that $\Pc(\bx^{(0)};f^*)\leq  \Pc(\bx^{(0)};r)$ which holds because $r<f^*$.	Note that $\Pc(\bx^{(0)};r)-(f^*-r)\geq (1/\alpha-1)(f^*-r)\geq 0$.
Now, define
\begin{eqnarray}
\label{eq:tkbound_apl}
t :=  \max\left\{\,\sqrt{\frac{2 \,C_2\, G^{2/d}}{\alpha\,(B\,-\,\alpha)\, \Pc(\bx^{(0)};r)^{1-2/d}}},\;
\sqrt{\frac{4\,C_1}{B\,-\,\alpha}}
\,\right\}- 1\geq0.
\end{eqnarray}
Plugging this $t$ into~\eqref{eq:ineq2_apl} gives:
\begin{align*}
\Pc(\bx^{(t)};r)\leq& f^*-r  + (B-\alpha)  [\Pc(\bx^{(0)};r)-(f^*-r)]+\alpha(B-\alpha) \Pc(\bx^{(0)};r)\\
=&(1-B+\alpha)(f^*-r)  + (B-\alpha)  \Pc(\bx^{(0)};r)+\alpha(B-\alpha) \Pc(\bx^{(0)};r)\\
\leq&(1-B+\alpha)\alpha \Pc(\bx^{(0)};r)  + (B-\alpha)  \Pc(\bx^{(0)};r)+\alpha(B-\alpha) \Pc(\bx^{(0)};r)\\
=&  B\Pc(\bx^{(0)};r),\end{align*} where the second inequality holds by the assumption that $\alpha\Pc(\bx^{(0)};r)> f^*-r$. 	Our proof is complete because $t\leq n^{\mathtt{fom}}$ which follows from our assumption that $\Pc(\bx^{(0)};r)> \epsilon$. \hfill\BlackBox

\begin{remark}
When $d < 2$ and $\epsilon$ is sufficiently small, $n^{\mathtt{fom}}$ drops to $\mathcal{O}(1)$. That is because the condition $\Pc(\bx^{(t)}; r) \le B \Pc(\bx^{(0)}; r)$ is satisfied after a constant number of prox-linear steps - as long as $\Pc(\bx^{(0)}; r) > \epsilon$. This phenomenon is due to the error bound condition: when $d < 2$, the landscape of $\Pc(\bx; r)$ becomes sharp enough near the optimal solutions set, and the first term in the maximum in \eqref{eq:tkbound_apl} decays rapidly as a function of $\Pc(\bx^{(0)}; r)$. As a result, the overall complexity bound in~\eqref{eq:total.complexity} remains valid and logarithmic in $1/\epsilon$. 
\end{remark}

\begin{corollary}\label{cor:smooth_alp.complexity}
	Suppose $\alpha$, $B$, $\gamma$, $\epsilon$, and $r_{\text{ini}}$ are such that $\alpha<B<1$, $\gamma >1$, $\epsilon > 0$, and $r_{\text{ini}} < f^*$. Consider $\tilde K$ defined in \eqref{eq:tildeK} for $r= r_{\text{ini}}$ and let $\eta^{(\text{ini})} > 0$, $L > 0$, and $\beta_{\text{dec}} \in (0,1)$. Suppose Algorithm~\ref{algo:RestartingLS} is executed with the accelerated prox-linear method described in Algorithm~\ref{algorithm:accelerated_prox_linear} as $\mathtt{fom}_k$ instances for $k = 0,1,\ldots, \tilde K$. Then, the algorithm finds an $\epsilon$-optimal and $\epsilon$-feasible solution in at most 

    \begin{align*}
    & (m+1)(\tilde{K}+1)\cdot\mathcal{O}\left(
    \max\left\{
    \sqrt{ \frac{ \max\left\{ \frac{1}{\eta^{(\text{ini})}}, \frac{L}{\beta_{\text{dec}}} \right\} G^{2/d} }{ \alpha(B - \alpha) \epsilon^{1 - 2/d} } },
    \sqrt{ \frac{ \max\left\{ 1, \frac{ \eta^{(\text{ini})} L }{ \beta_{\text{dec}} } \right\} }{ B - \alpha } }
    \right\}
    \right) \cdot \tilde{D},\\
    &= \mathcal{O}\left(\log^2\left(\dfrac{1}{\epsilon}\right)\max\left\{\frac{1}{\epsilon^{1/2-1/d}},1\right\}\right),
    \end{align*} 
    gradient and function evaluations, where $\tilde D$ is defined in \eqref{eq:Dtilde}.
\end{corollary}
\proof{}
The proof of this corollary is similar to the proof of Corollary~\ref{cor:non.smooth.complexity} except that \eqref{eq:peritercostsmooth} is used as $n^{\mathtt{fom}}$ and follows from the number of prox-linear subproblems to solve, which equals $\mathcal{O}(n^{\mathtt{fom}})$ (see the second conclusion of Theorem~\ref{theorem:APLM}).
\hfill\BlackBox

\section{Numerical Experiments}\label{sec:exp}
We assess the numerical effectiveness of Algorithm~\ref{algo:RestartingLS} on two types of binary classification problems: those with fairness constraints and those with Type~II error constraints. This evaluation serves as a benchmark, comparing our proposed method against existing approaches in the literature. For the fairness-constrained setting, we consider both (1) a non-smooth formulation, using the standard subgradient algorithm as the $\mathtt{fom}$ in Algorithm~\ref{algo:RestartingLS}, and (2) a smooth formulation, using Algorithm~\ref{algorithm:accelerated_prox_linear} as the $\mathtt{fom}$.  For the Type~II error-constrained problems, we focus on the non-smooth formulation and use the subgradient method as the $\mathtt{fom}$, with the goal of evaluating performance on larger-scale instances.

\subsection{Non-smooth Classification Problem With Fairness Constraints}
\label{sec:num_non-smooth}
Consider a set of $n$ data points $\mathcal{D}=\{(\ba_i,b_i)\}_{i=1}^n$ where $\ba_i\in\mathbb{R}^p$ denotes a feature vector, and $b_i\in\{1,-1\}$ represents the class label for $i=1,2,\dots,n$. Let $\mathcal{D}_M\subset \mathbb{R}^p$ and $\mathcal{D}_F\subset \mathbb{R}^p$ be two different sensitive groups of instances. We want to find a linear classifier $\bx\in\mathbb{R}^p$ that not only predicts the labels well (i.e., minimizes a loss function) but also treats each sensitive instance from $\mathcal{D}_M$ and $\mathcal{D}_F$ fairly. Such classification problems are suitable for applications such as loan approval or hiring decisions. For example, in the context of loan approval, we can ensure that loans are equally provided for different applicants independent of their group memberships (e.g., home-ownership or their gender). A correct classifier $\bx$ satisfies $b_i \ba_i^\top \bx > 0$ for all $i$. One can train such a classifier by solving the following optimization problems that minimizes a non-increasing convex loss function $\phi(\cdot)$ subject to fairness constraints:
	\begin{eqnarray}
		\label{eq:DDCproblem}
		\min_{\bx} &&\frac{1}{n}\sum_{i=1}^n\phi(-b_i\ba_i^\top\bx)\\\label{eq:firstconstraint}
		\text{s.t.}&&\frac{1}{n_F}\sum_{\ba\in\mathcal{D}_F}\sigma(\ba^\top\bx) \geq \frac{\kappa}{n_M}\sum_{\ba\in\mathcal{D}_M}\sigma(\ba^\top\bx),\\\label{eq:secondconstraint}
		&&\frac{1}{n_M}\sum_{\ba\in\mathcal{D}_M}\sigma(\ba^\top\bx) \geq \frac{\kappa}{n_F}\sum_{\ba\in\mathcal{D}_F}\sigma(\ba^\top\bx),
	\end{eqnarray}
	where  $\kappa\in(0,1]$ is a constant, $n_M$ and $n_F$ are respectively the number of instances in $D_M$ and $D_F$, and finally	$\sigma(\ba^\top\bx) :=\max\{0,\min\{1,\{0.5+\ba^\top\bx\}\}\in[0,1]$ is the probability of predicting $\ba$ in class $+1$. The constraint~\eqref{eq:firstconstraint} guarantees that the percentage of the instances in $\mathcal{D}_F$ that are predicted in class $+1$ is at least a $\kappa > 0$ fraction of that in $\mathcal{D}_M$. The constraint~\eqref{eq:secondconstraint} has a similar interpretation for instances in $\mathcal{D}_M$. Choosing an appropriate $\kappa$ ensures that the obtained classifier is fair to both $\mathcal{D}_M$ and $\mathcal{D}_F$ groups. An analogous model was considered in \citet{goh2016satisfying}.
	
	Constraints \eqref{eq:firstconstraint} and \eqref{eq:secondconstraint} are non-convex, so we approximate this problem by a convex optimization problem following the approach described in~\citet{lin2019data}. In particular,  we reformulate the first constraint as 
	$
	\frac{\kappa}{n_M}\sum_{\ba\in\mathcal{D}_M}\sigma(\ba^\top\bx)+\frac{1}{ n_F}\sum_{\ba\in\mathcal{D}_F}\sigma(-\ba^\top\bx)\leq 1
	$
	by using $\sigma(\ba^\top\bx)=1-\sigma(-\ba^\top\bx)$. 	In addition, the function $\sigma(\ba^\top\bx)$ can be approximated by $(0.5+\ba^\top\bx)_+ = \max\{0,0.5+\ba^\top\bx\}$. Hence, we can rewrite the constraint~\eqref{eq:firstconstraint} as a convex constraint:
	$
	\frac{\kappa }{n_M}\sum_{\ba\in\mathcal{D}_M}(\ba^\top\bx+0.5)_++\frac{1}{n_F}\sum_{\ba\in\mathcal{D}_F}(-\ba^\top\bx+0.5)_+\leq 1. 
	$
	Applying a similar convex approximation to~\eqref{eq:secondconstraint}, we obtain the following convex reformulation of~\eqref{eq:DDCproblem}:
	\begin{eqnarray}
		\label{eq:DDCproblem_convex}
		\min_{\bx} &&\frac{1}{n}\sum_{i=1}^n\phi(-b_i\ba_i^\top\bx)\\\nonumber
		\text{s.t.}&&\frac{\kappa }{n_M}\sum_{\ba\in\mathcal{D}_M}(\ba^\top\bx+0.5)_++\frac{1}{n_F}\sum_{\ba\in\mathcal{D}_F}(-\ba^\top\bx+0.5)_+\leq 1,\\\nonumber
		&&\frac{\kappa }{n_F}\sum_{\ba\in\mathcal{D}_F}(\ba^\top\bx+0.5)_++\frac{1}{n_M}\sum_{\ba\in\mathcal{D}_M}(-\ba^\top\bx+0.5)_+\leq 1.
	\end{eqnarray}
The above reformulation contains piece-wise linear functions (non-smooth) in its objective and constraints. Hence, this problem satisfies the error bound condition~\eqref{eq:errorbound} with $d=1$ and an unknown value of $G$. We implemented RLS (Algorithm~\ref{algo:RestartingLS}) to solve~\eqref{eq:DDCproblem_convex} using subgradient method in each $\mathtt{fom}_k$ instance to solve the subproblem $\min_{\bx\in\mathcal{X}}\Pc(\bx_k, r_k)$. We compare the performance of our method with three parameter-free but non-adaptive benchmarks: (1) the feasible level-set (FLS) method described in~\citet{lin2018feasiblelevel}; (2) the stochastic subgradient method by Yu, Neely, and Wei (YNW) in~\citet{yu2017online}; and (3) the switching gradient (SWG) method  in~\citet{bayandina2018mirror} which is based on a mirror descent algorithm and a switching sub-gradient scheme \footnote{We do not benchmark against adaptive methods for constrained convex optimization as they require knowledge of unknown parameters. Recall that an adaptive algorithm can harness the EBC when it holds, leading to improved convergence results.}. We choose these benchmarks as they are the most representative first-order methods that can solve ``non-smooth'' convex ``constrained'' optimization problems. The FLS is a level-set method that guarantees the feasibility of the solutions generated at each iteration, but does not exploit the error bound condition as our method does. Hence, it has a higher complexity bound for the class of problems that satisfy the error bound condition~\eqref{eq:errorbound}. YNW is a primal-dual method with dual variables updated implicitly, while SWG is a pure primal method. We used the number of equivalent data passes performed by each algorithm as a measure of complexity. \medskip

\noindent{\bf Instances:} To check the performance of the above mentioned methods, we considered three classification problems using (i) ``LoanStats'' data set from lending club\footnote{\url{https://www.lendingclub.com/info/statistics.action}} which is a platform that allows individuals lend to other individuals; (ii) ``COMPAS'' data set from ProPublica\footnote{\url{https://github.com/propublica/compas-analysis}}; and (iii) ``German'' data set from UCI Machine Learning Repository\footnote{\url{https://archive.ics.uci.edu/ml/datasets/statlog+(german+credit+data)}}.

The LoanStats data set contains information of $128,375$ loans issued on the loan club platform in the fourth quarter of 2018. The goal in this application is to predict whether a loan request will be approved or rejected. After creating dummy variables, each loan in this data set is represented by a vector of 250 features. We randomly partitioned this dataset into two sets: one of $63,890$ examples to construct the objective function and another of $64,485$ examples to construct the constraints. The set $D_M$ in this application denotes the set of instances with ``home-ownership = Mortgage''. All the other instances are considered in the set $D_F$. The fairness constraints in this application guarantee that customers with mortgages will not be disproportionately affected in receiving new loans. 

The COMPAS dataset includes criminal history, jail and prison time, demographics, and other factors for 6,172 instances. The goal of this application is to predict whether a person will be rearrested within two years after their initial arrest. This information, for example, helps judges make bail
decisions by predicting the criminal recidivism of defendants. We used $4,115$ examples from this dataset to formulate the objective function, and the remaining examples to formulate the constraints. We chose sets $D_M$ and $D_F$ to respectively be the sets of male and female instances to ensure fair treatment among different genders. 

The German credit dataset describes financial details of customers and is used to determine whether the customer should be granted credit (i.e., is a good customer) or not (i.e., is a bad customer). This dataset contains 1,000 instances and 50 variables. We used $667$ of the examples to formulate the objective function and the remaining to formulate the constraints. Similar to COMPAS data, the sets $D_M$ and $D_F$ respectively denote male and female examples. 
	
In all above applications, we chose $\kappa=0.9$, and $\phi(z)=(1-z)_+$. We terminated RLS and the other three benchmarks when the number of data passes reached 20,000. However, we ran the two level-set methods, i.e., RLS and FLS, longer for 100,000 data passes to find the optimal value $f^*$. In particular, we selected the smallest objective value among all feasible solutions as a close approximation of $f^*$. We used this value to compute and plot $P(x;f^*)$.\medskip

\noindent {\bf Algorithmic configurations:} We next describe how we selected each of the parameters used in RLS, FLS, YNW, and SWG.
\begin{itemize}
\item RLS: We choose a value for $\epsilon$ from $\{1,0.1,0.01,0.001\}$, for $\alpha$ from $\{0.4,0.5,0.7,0.9\}$ and for $B$ from $\{0.5,0.9,0.95,0.99\}$ to minimize the value of $\Pc(\bx; f^*)$ at termination (we only consider the combination with $\alpha\leq B$). Since the optimal value of~\eqref{eq:DDCproblem_convex} is positive, we set $r_{\text{ini}} = 0$ to satisfy the condition $r_{\text{ini}} < f^*$. Although the vector of all zeros could be used as an initial solution $\bx_{\text{ini}}$, we used a different value in our algorithm. In particular, the complexity bound of RLS shows that the algorithm requires fewer iterations if $g(\bx_{\text{ini}}) < 0$ and is as small as possible. Hence, we followed a heuristic to obtain a better quality solution. We applied the subgradient algorithm to minimize the function $g(\bx)$ and used the returned solution after 40 iterations. Finding such a solution accounted for less than 1\% of our total-run time. Finally, we used $K = \tilde K$ where $\tilde K$ is computed as in~\eqref{eq:tildeK} for $r=r_{\text{ini}}$. 
\item FLS: The FLS algorithm described in~\citet{lin2018feasiblelevel} requires two input parameters: $r_{\text{ini}} > f^*$ and $\alpha > 1$. To obtain $r_{\text{ini}} > f^*$, we first followed the same step as in RLS to obtain a feasible solution $\bx_{\text{ini}}$ such that $g(\bx_{\text{ini}}) < 0$ and then set $r_{\text{ini}} = f(\bx_{\text{ini}})$. Since $\bx_{\text{ini}}$ is feasible, using Lemma~\ref{lemma:levelset_properties} we can guarantee that $r_{\text{ini}} > f^*$. We set $D=5\|\bx_{\text{ini}}\|$ and chose $\alpha$ from $\{100,200,500\}$ that produces the smallest $\Pc(\bx; f^*)$ at termination. 
\item YNW: We followed the guidance in \citet{yu2017online} to setup YNW. Specifically, we chose the control parameters $V$ and $\alpha$ as $V =\sqrt{T}$ and 
	$\alpha=T$, respectively, as a function of the total number of iterations $T$, where $V$ is the weight of the gradient of the objective function and $\alpha$ is the weight of the proximal term in the updating equation of $\bx$ in YNW. We chose the number of iterations $T$ such that the total number of data passes is $20,000$. We used the all-zero vector as an initial solution to YNW. 
\item SWG: The only input of SWG algorithm is the precision $\epsilon$. We selected $\epsilon$ from $\{0.01, 0.001, 0.0001\}$ to minimize the value of $\Pc(x; f^*)$ achieved by SWG at termination. The initial solution is chosen as a feasible solution $\bx_{\text{ini}}$ found by the same step as in RLS and FLS. Of course, for a fair comparison, the number of data passes that RLS, FLS and SWG spend in searching for the initial solution is included when comparing the performance of these algorithms. 
\medskip
\end{itemize}
	
\begin{figure}[t!]
    \centering
    \resizebox{\textwidth}{!}{%
        \begin{tabular}{@{}c|ccc@{}}
            & LoanStats & COMPAS & German \\
            \hline
            \raisebox{8ex}{\small\rotatebox[origin=c]{90}{$\Pc(\bx; f^*)$}}
            & \includegraphics[width=0.2\textwidth]{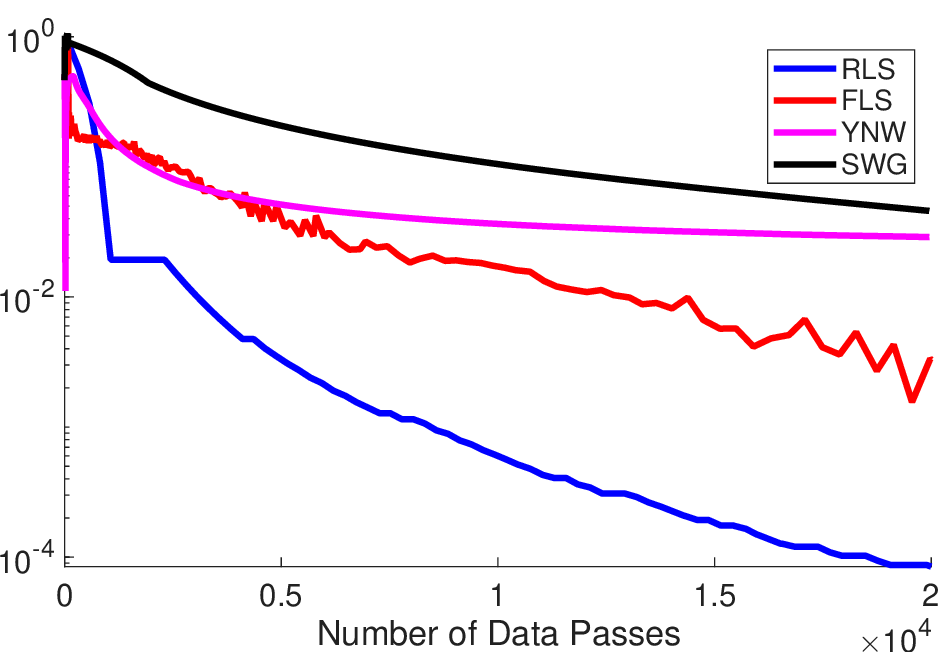}
            & \includegraphics[width=0.2\textwidth]{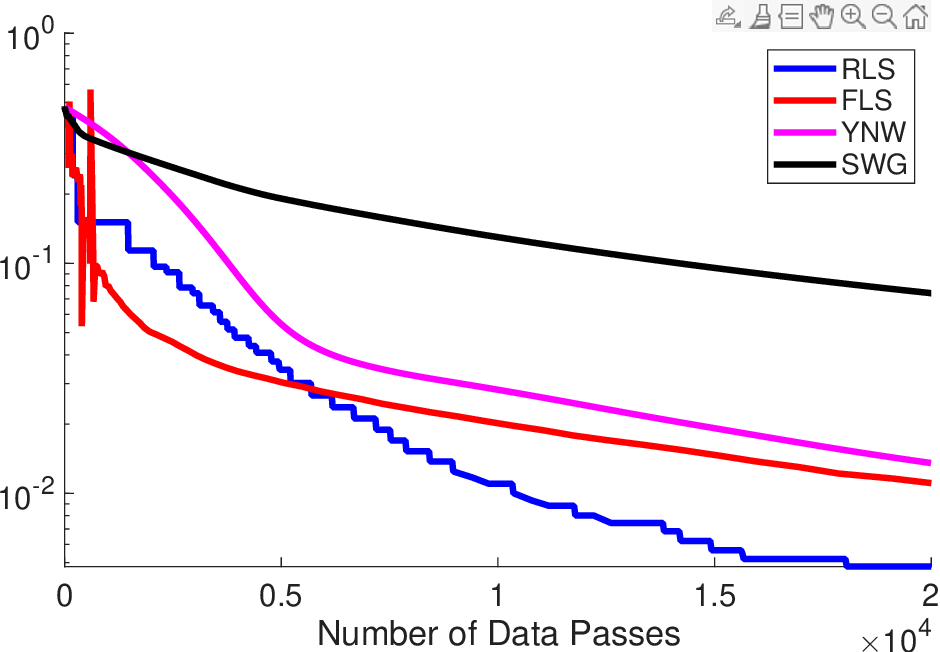}
            & \includegraphics[width=0.2\textwidth]{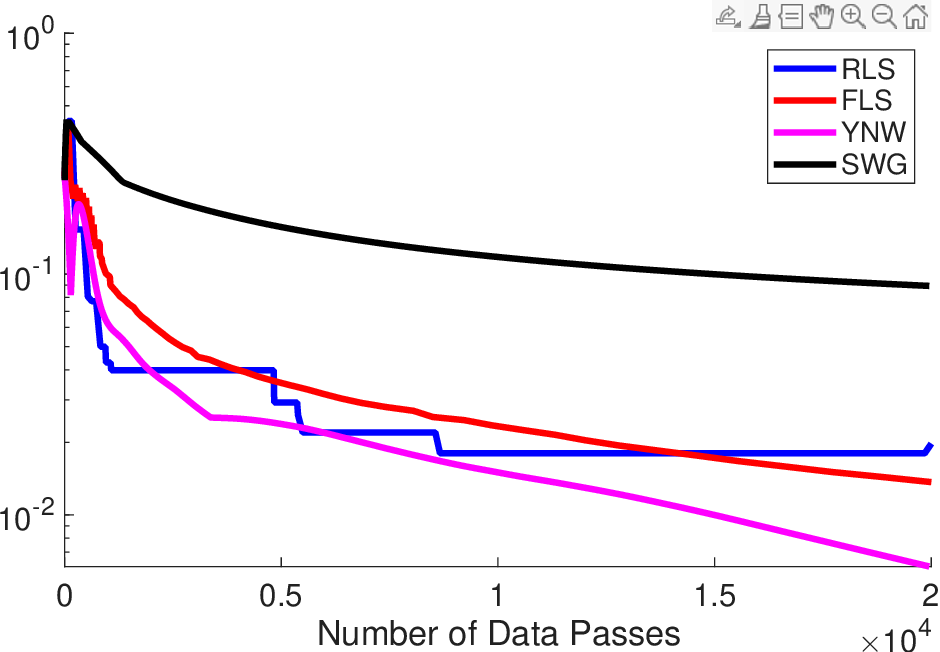} \\
            \raisebox{8ex}{\small\rotatebox[origin=c]{90}{$f(x)-f^*$}}
            & \includegraphics[width=0.2\textwidth]{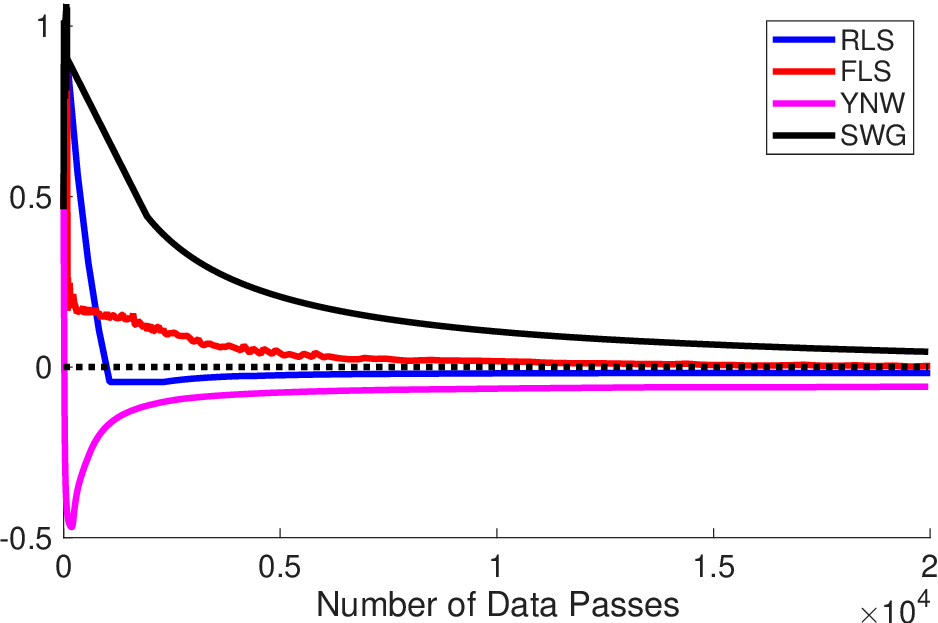}
            & \includegraphics[width=0.2\textwidth]{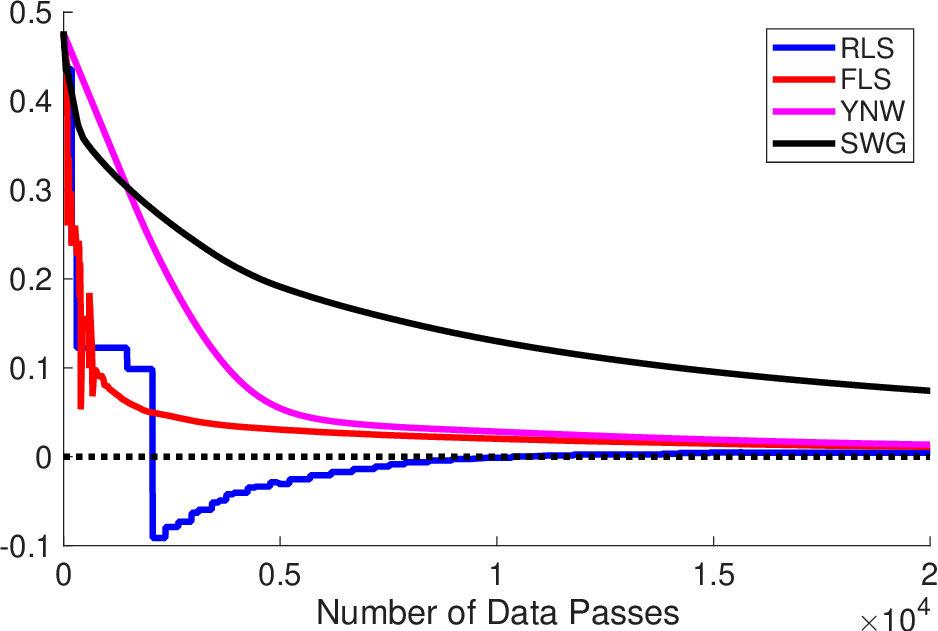}
            & \includegraphics[width=0.2\textwidth]{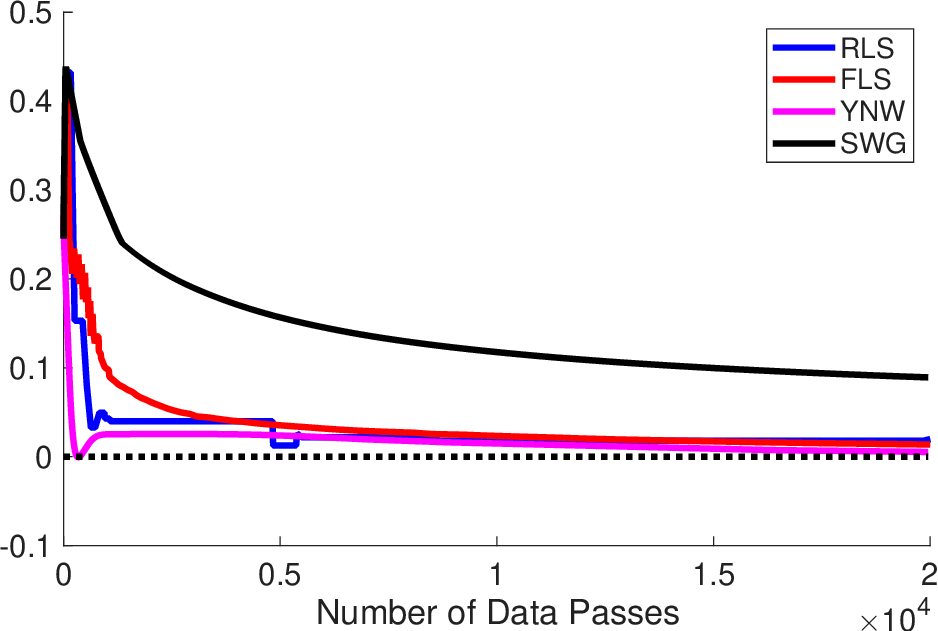} \\
            \raisebox{8ex}{\small\rotatebox[origin=c]{90}{$ \max\limits_{i=1,2,\dots,m}\{f_i(\bx)\}$}}
            & \includegraphics[width=0.2\textwidth]{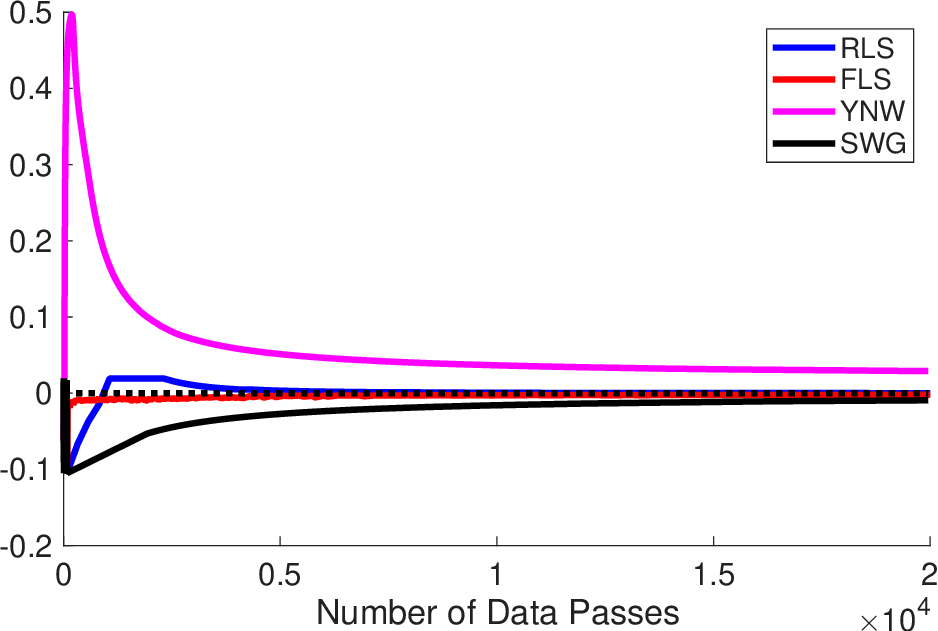}
            & \includegraphics[width=0.2\textwidth]{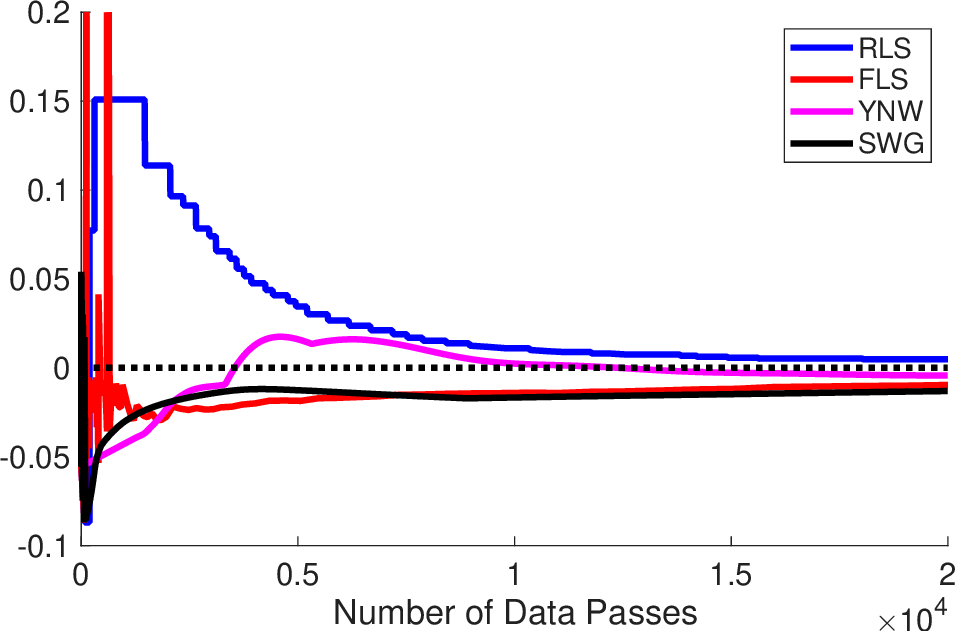}
            & \includegraphics[width=0.2\textwidth]{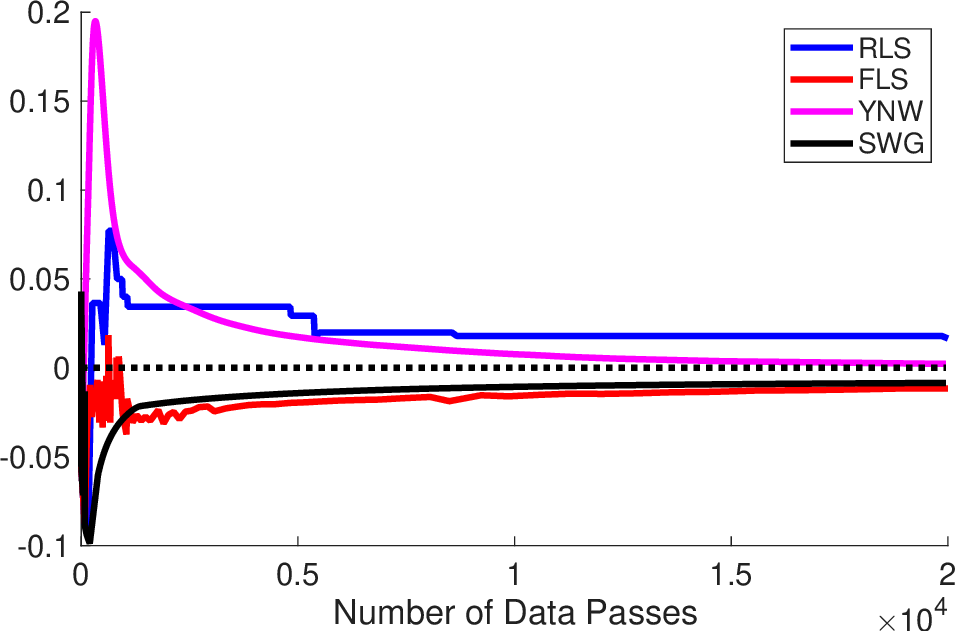} \\
        \end{tabular}%
    }
    \vspace{2ex}
    \caption{Performance of RLS, FLS, YNW, and SWG for solving non-smooth binary classification problems with fairness constraints.}
    \label{fig:fairness}
\end{figure}

\noindent {\bf Results:} Figure \ref{fig:fairness} displays the performance of RLS, FLS, YNW and SWG as a function of data passes. In particular, the $x$-axis in this figure shows the number of data passes each algorithm performed while the $y$-axis represents the feasibility and optimality of solutions returned by each method. More specifically, the $y$-axis in the first row of this plot represents $\Pc(\bx;f^*)$. As explained before, this quantity is zero when the solution $\bx$ is both feasible and optimal and positive otherwise. The speed at which $\Pc(\bx;f^*)$ converges to zero indicates how fast the solution $\bx$ becomes feasible and optimal at each iteration. The $y$-axis in the second row represents $f(\bx) - f^*$ that measures the optimality gap at $\bx$ and in the third row, $g(\bx) = \max_{i = 1, \ldots, m} \{f_i(\bx)\}$ checks whether the solution $\bx$ is feasible or not (if $g(\bx) < 0$ then $\bx$ is feasible).  

We observe that RLS achieves the fastest reduction in the function $\Pc(\bx; f^*)$ across two out of the three datasets. Recall that $\Pc(\bx; f^*)$ equals the maximum of infeasibility and suboptimality of $\bx$. Our computational results reveal that our algorithm compares favorably with the benchmarks on these two datasets. However, on the German credit dataset, RLS doesn't exhibit the best performance, possibly due to the problem satisfying the error bound condition ($d=1$) with a large $G$, causing our method to have higher complexity than FLS and YNW. It can also be seen that RLS and YNW may generate infeasible solutions before convergence, while FLS and SWG maintain feasibility at each iteration. 

\subsection{Smooth Classification Problem With Fairness Constraints}
\label{sec:num_smooth}
In this section, we revisit the binary classification problem with fairness constraints outlined in \eqref{eq:DDCproblem}-\eqref{eq:secondconstraint}. However, this time, we employ a different approximation for the function $\sigma(\ba^\top\bx)$ by using $\log_4(1+\exp(\ba^\top\bx))$, and adopt $\phi(z)=\ln(1+\exp(-z))$. The choice of the base-four logarithm ensures that $\log_4(1+\exp(\ba^\top\bx))\geq 0.5$ when $\ba^\top\bx \geq 0$, aligning with our prior selection of $\sigma(\ba^\top\bx) = (0.5+\ba^\top\bx)_+$.

Therefore, we can formulate our classification problem with fairness constraints as
\begin{eqnarray}
	\label{eq:DDCproblem_convex_smooth}
	\min_{\bx} &&\frac{1}{n}\sum_{i=1}^n\ln(1+\exp(-b_i\ba_i^\top\bx))\\\nonumber
	\text{s.t.}&&\frac{\kappa }{n_M}\sum_{\ba\in\mathcal{D}_M}\log_4(1+\exp(\ba^\top\bx))+\frac{1}{n_F}\sum_{\ba\in\mathcal{D}_F}\log_4(1+\exp(-\ba^\top\bx))\leq 1,\\\nonumber
	&&\frac{\kappa }{n_F}\sum_{\ba\in\mathcal{D}_F}\log_4(1+\exp(\ba^\top\bx))+\frac{1}{n_M}\sum_{\ba\in\mathcal{D}_M}\log_4(1+\exp(-\ba^\top\bx))\leq 1.
\end{eqnarray}

The convex optimization problem presented above features differentiable objectives and constraints with Lipschitz continuous gradients. While the satisfaction of the error bound condition~\eqref{eq:errorbound} is generally uncertain, our RLS method remains applicable and automatically leverages the error bound condition for any values of $d$ and $G$, when \eqref{eq:errorbound} is valid, thanks to its adaptive nature.

Notably, even in the absence of a strongly convex regularization term like $\frac{1}{2}\|\bx\|_2^2$ in the objective and constraints, when the dataset size significantly exceeds the dimension of $\ba$, it is plausible that both the objective and constraints exhibit $\mu$-strong convexity. This implies satisfaction of \eqref{eq:errorbound} with $d=2$ and $G=\mu^{-1}$. In such instances, our method boasts a theoretical complexity of $\tilde{\mathcal{O}}(1/\sqrt{\epsilon})$ without necessitating knowledge of $d$ and $\mu$.

In contrast, existing non-adaptive methods require prior knowledge of $d$ and $\mu$ to achieve a similar complexity. Without this information, they are limited to a complexity of $\tilde{\mathcal{O}}(1/\epsilon)$.

We employed the Restarting Level Set (RLS) algorithm, detailed in Algorithm~\ref{algo:RestartingLS}, to solve the convex and smooth problem specified in \eqref{eq:DDCproblem_convex_smooth}. In each instance $\mathtt{fom}_k$, we utilized the accelerated proximal linear method outlined in Algorithm~\ref{algorithm:accelerated_prox_linear} to solve the subproblem $\min_{\bx\in\mathcal{X}}\Pc(\bx, r_k)$, where the subproblem in Step 6 in Algorithm~\ref{alg:backtracking} is solved by applying $\mathtt{quadprog}$ function in MATLAB to solve its dual problem, a convex quadratic program in $\mathbb{R}^3$, and recovering the primal solution.

Our method's performance is benchmarked against two non-adaptive approaches: (1) Feasible Level Set (FLS) from~\citet{lin2018feasiblelevel}, introduced previously; (2) Constraint Extrapolation (ConEx) method from~\citet{boob2023stochastic}, a single-loop primal-dual algorithm considered a variant of well-known primal-dual accelerated gradient methods for min-max optimization~\citep{chambolle2011first,hamedani2021primal}. We selected these benchmarks as they represent prominent first-order methods capable of solving ``smooth" convex ``constrained" optimization problems, achieving optimal complexity for both convex and strongly convex problems.\medskip

\noindent{\bf Instances:} To evaluate these methods, we tackled three classification problems using the same datasets as in the non-smooth case. Across all applications, we set $\kappa=0.9$. Termination criteria for RLS and the two benchmarks were based on reaching 80,000 data passes. However, we extended the runs for the two level-set methods, RLS and FLS, to 100,000 data passes to pinpoint the optimal value $f^*$. Specifically, we selected the smallest objective value among all feasible solutions as a close approximation of $f^*$, using this value to compute and plot $\Pc(x;f^*)$.\medskip

\noindent {\bf Algorithmic configurations:} 
The parameters for RLS and FLS methods were chosen by the same procedure as in the non-smooth case, except for RLS, where we employed the accelerated proximal linear method with line search outlined in Algorithm~\ref{algorithm:accelerated_prox_linear} to solve the subproblem $\min_{\bx\in\mathcal{X}}\Pc(\bx, r_k)$. We set $\eta^{(\text{ini})}=1$, $\beta_{\text{dec}}=0.8$ and $\beta_{\text{inc}}=1.2$ in Algorithm~\ref{algorithm:accelerated_prox_linear}.

For the ConEx method, without assuming strong convexity, we adhered to the parameter settings recommended in Theorem 2 of~\citet{boob2023stochastic}. Specifically, we chose $\gamma_t=1$, $\theta_t=1$, $\tau_t=\tau$, and $\eta_t=\eta$ for each iteration $t$ in ConEx. The values of $\tau$ and $\eta$ were selected from the set $\{1, 10, 100\}$ to minimize the value of $\Pc(\bx; f^*)$ at termination.

\medskip

\begin{figure}[t!]
    \centering
    \resizebox{\textwidth}{!}{%
        \begin{tabular}{@{}c|ccc@{}}
            & LoanStats & COMPAS & German \\
            \hline
            \raisebox{8ex}{\small\rotatebox[origin=c]{90}{$\Pc(\bx; f^*)$}}
            & \includegraphics[width=0.28\textwidth]{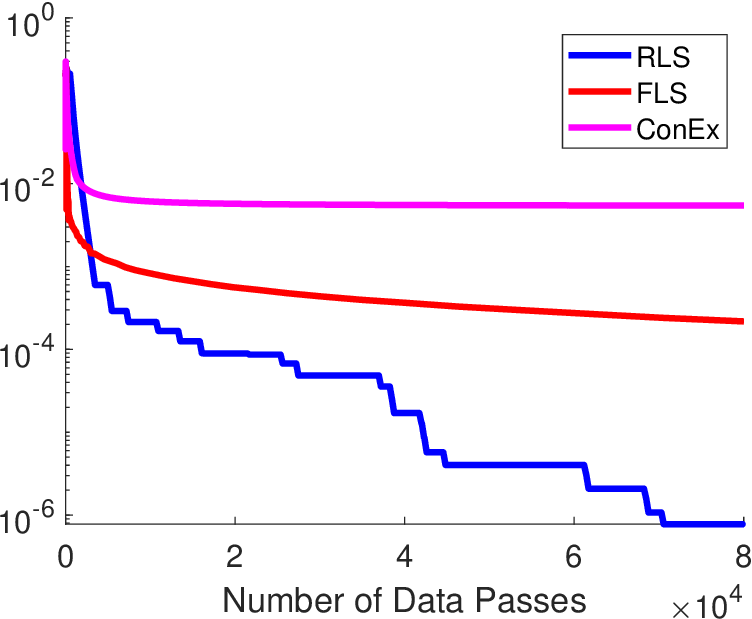}
            & \includegraphics[width=0.28\textwidth]{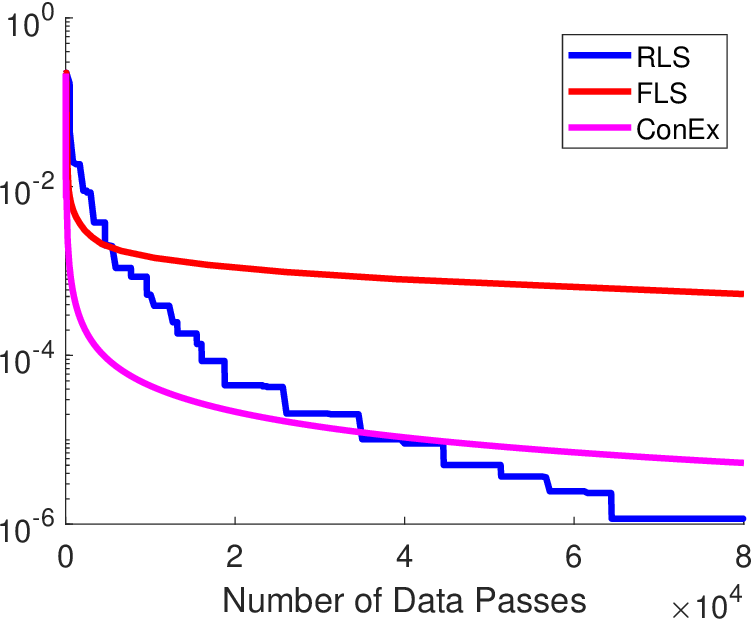}
            & \includegraphics[width=0.28\textwidth]{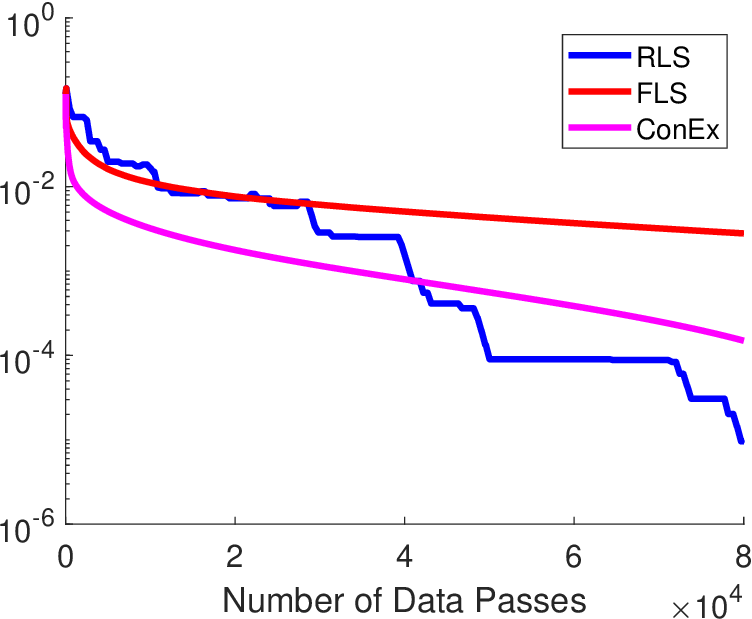} \\
            \raisebox{8ex}{\small\rotatebox[origin=c]{90}{$f(x)-f^*$}}
            & \includegraphics[width=0.28\textwidth]{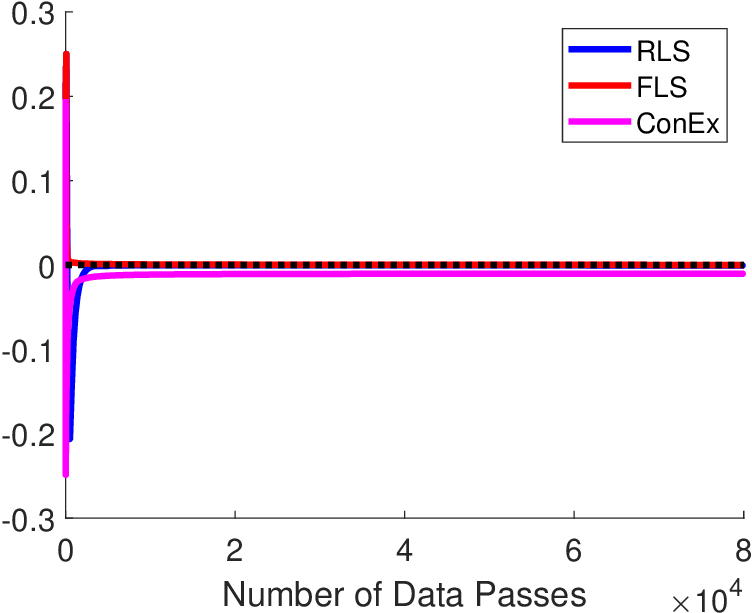}
            & \includegraphics[width=0.28\textwidth]{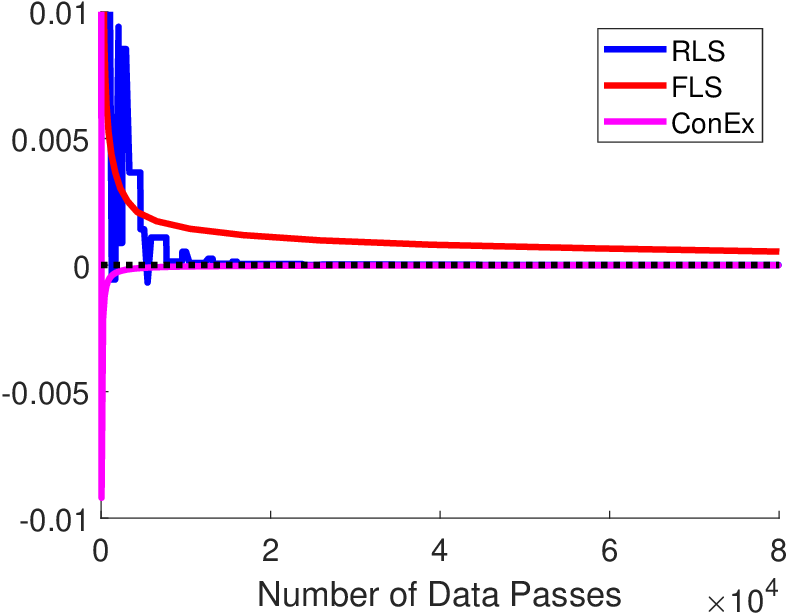}
            & \includegraphics[width=0.28\textwidth]{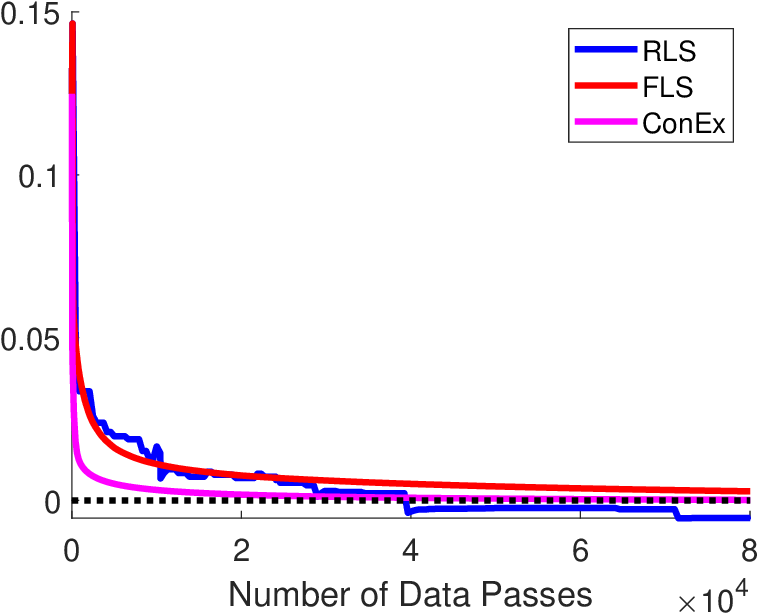} \\
            \raisebox{8ex}{\small\rotatebox[origin=c]{90}{$ \max\limits_{i=1,2,\dots,m}\{f_i(\bx)\}$}}
            & \includegraphics[width=0.28\textwidth]{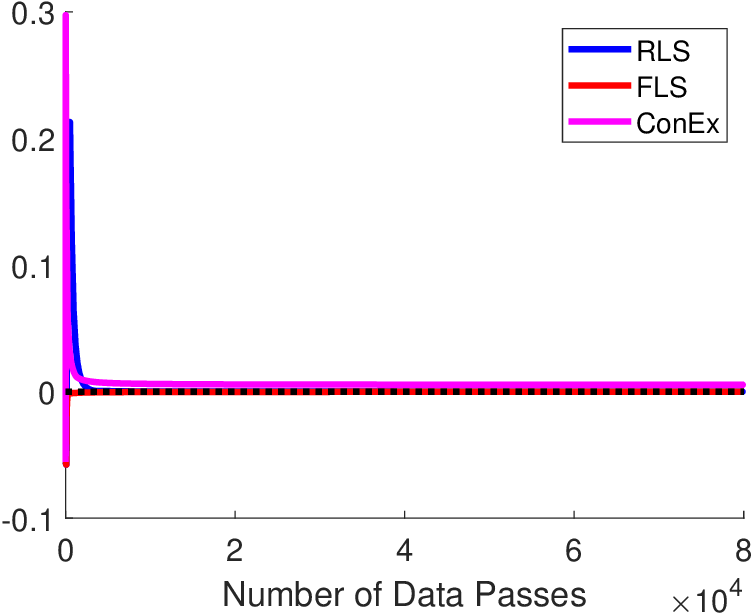}
            & \includegraphics[width=0.28\textwidth]{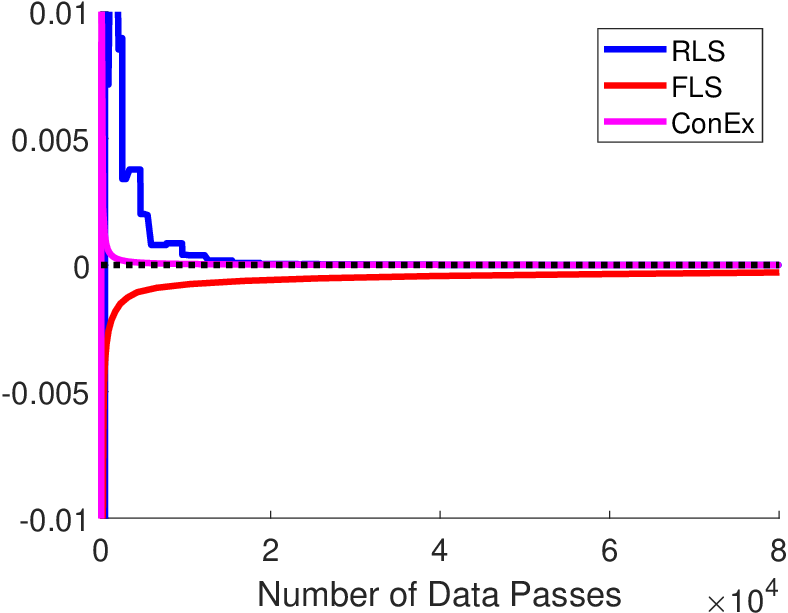}
            & \includegraphics[width=0.28\textwidth]{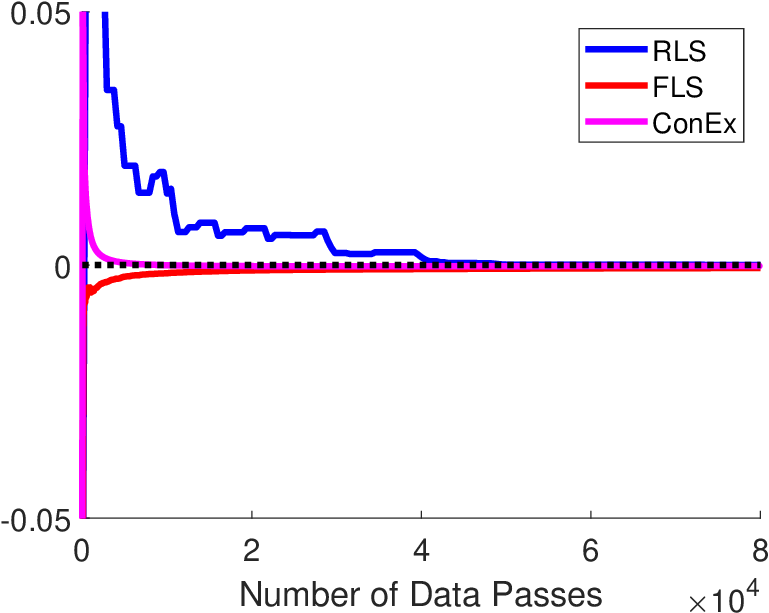} \\
        \end{tabular}%
    }
    \vspace{2ex}
    \caption{Performance of RLS, FLS, and ConEx for solving smooth binary classification problems with fairness constraints.}
    \label{fig:fairness_smooth}
\end{figure}

\noindent {\bf Results:} Figure \ref{fig:fairness_smooth} illustrates the performance of RLS, FLS, and ConEx over varying data passes. The axes' meanings remain consistent with Figure \ref{fig:fairness}.

Notably, RLS demonstrates the swiftest reduction in the function $\Pc(\bx; f^*)$ to a small scale (e.g. $10^{-5}$) on the three datasets. FLS stands out as the sole algorithm maintaining feasibility at each iteration, while ConEx diminishes the infeasibility faster than RLS on the COMPAS dataset and the German credit dataset. RLS, despite not being the fastest in reducing infeasibility, promptly achieves a high-quality solution in the sense of having a small $\Pc(\bx; f^*)$, showcasing its favorable comparison with existing approaches in the literature.

\subsection{Non-smooth Neyman–Pearson Classification with Type II Error Control}
\label{sec:num_nonsmooth_NP_new}

In this section, we evaluate the performance of our RLS method on a Neyman-Pearson classification problem, where the goal is to minimize the Type I error (i.e., the misclassification of negative samples as positive) subject to a Type II error constraint (i.e., ensuring that the misclassification rate for positive samples remains below a given threshold). In contrast to our previous experiments on fairness constraints, this instance emphasizes the control of Type II error and involves larger problem instances.

Consider a training set of $n$ data points $\mathcal{D}=\{(\ba_i,b_i)\}_{i=1}^n$, where $\ba_i\in\mathbb{R}^p$ denotes a feature vector and $b_i\in\{1,-1\}$ is the corresponding class label. Let 
$\mathcal{D}_+ = \{ (\ba,b)\in\mathcal{D} \mid b=1 \}$ and $ \mathcal{D}_- = \{ (\ba,b)\in\mathcal{D} \mid b=-1 \}.$
The Neyman–Pearson linear classification problem aims to find a vector $\bx\in \mathbb{R}^p$ that minimizes the empirical Type I error (computed over $\mathcal{D}_-$) while ensuring that the empirical Type II error (computed over $\mathcal{D}_+$) does not exceed a specified threshold $\kappa$. By approximating both error types with a non-increasing convex surrogate loss function $\phi(\cdot)$, the problem is formulated as
\begin{equation}
\label{eq:NPclassification_new}
\min_{\|\bx\|_2\leq R}\left\{\frac{1}{n_-}\sum_{(\ba,b)\in\mathcal{D}_-}\phi(-\ba^\top\bx)\right\} \quad \text{s.t.} \quad \frac{1}{n_+}\sum_{(\ba,b)\in\mathcal{D}_+}\phi(\ba^\top\bx)-\kappa\leq0,
\end{equation}
where $R>0$ is a regularization parameter, $n_+=|\mathcal{D}_+|$, $n_-=|\mathcal{D}_-|$, and the constraint enforces a bound on the Type II error. Problem \eqref{eq:NPclassification_new} is an instance of the general constrained convex optimization model discussed in Section~\ref{nonsmoothcase} (with $m=1$ and $\mathcal{X}=\{\bx: \|\bx\|_2\leq R\}$).

We consider three large-scale datasets: (i) ``gisette'' \citep{guyon2004result}, (ii) ``rcv1.binary'' \citep{lewis2004rcv1}, and (iii) ``real-sim''\footnote{\url{https://people.cs.umass.edu/~mccallum/code-data.html}}. These datasets were obtained from the LIBSVM library\footnote{\url{https://www.csie.ntu.edu.tw/~cjlin/libsvmtools/datasets/}}. The gisette dataset involves handwritten digit recognition and comprises $6,000$ data points with $5,000$ features; the rcv1.binary dataset, used for text categorization, has $20,242$ training examples with $47,236$ features; and the real-sim dataset consists of $72,309$ data points with $20,958$ features. Each feature vector is normalized to have unit Euclidean norm. For all instances, we set $\kappa=0.5$, choose $\phi(z)=(1-z)_+$, and set $R=1$ for gisette and $R=5$ for rcv1.binary and real-sim.

Since all these instances are non-smooth, we implemented Algorithm~\ref{algo:RestartingLS} using a subgradient method in each $\mathtt{fom}_k$ instance. Similar to Section~\ref{sec:num_non-smooth}, we compare RLS with the FLS, the YNW method, and the switching subgradient method (SWG) on these larger instances. In our experiments, 
the parameters of all methods are chosen by the same procedure as in Section~\ref{sec:num_non-smooth}. All methods are terminated when the total number of data passes reaches 40,000, while RLS is run for 100,000 data passes to obtain an accurate approximation of the optimal value $f^*$. Specifically, the smallest objective value among all feasible solutions produced by RLS is selected as an approximation of $f^*$ and is used to compute and plot the residual function $P(\bx;f^*)$.

\begin{figure}[t!]
    \centering
    \resizebox{\textwidth}{!}{%
        \begin{tabular}{@{}c|ccc@{}}
            & gisette & rcv1.binary & real-sim \\
            \hline
            \raisebox{8ex}{\small\rotatebox[origin=c]{90}{$P(\bx;f^*)$}}
            & \includegraphics[width=0.28\textwidth]{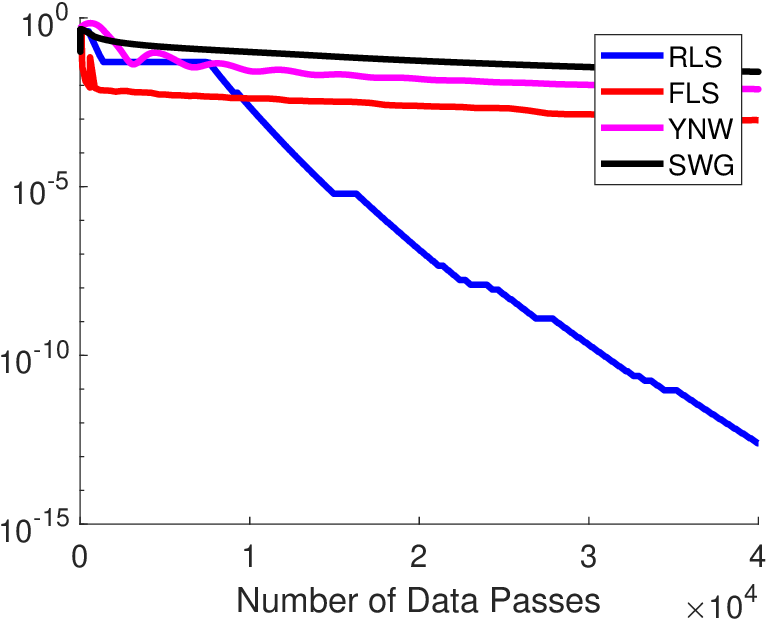}
            & \includegraphics[width=0.28\textwidth]{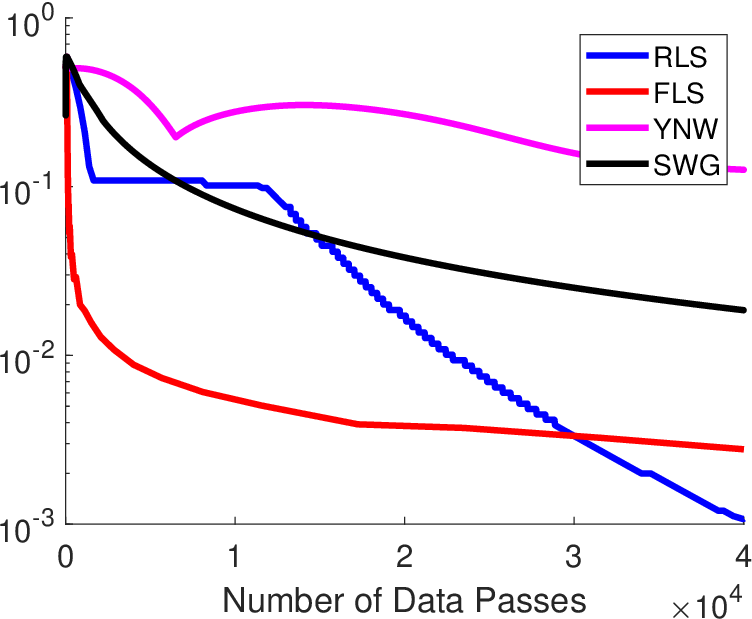}
            & \includegraphics[width=0.28\textwidth]{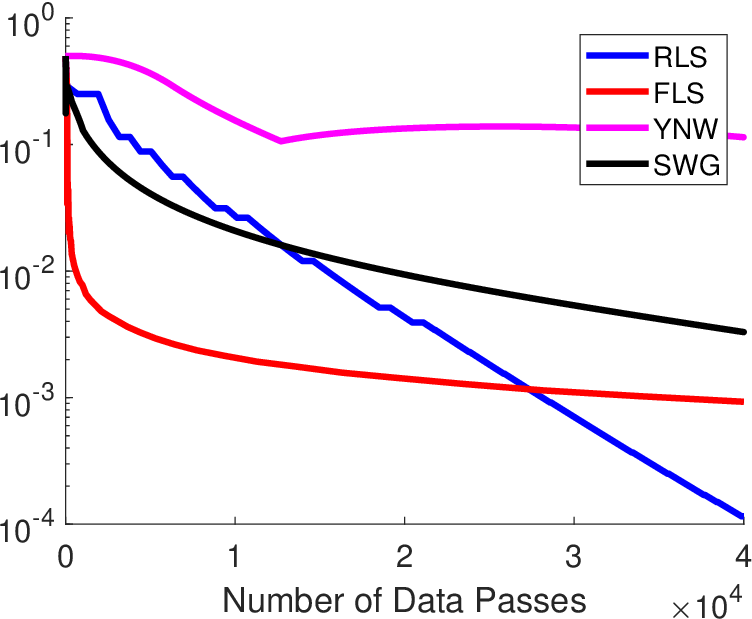} \\
            \raisebox{8ex}{\small\rotatebox[origin=c]{90}{$f(\bx)-f^*$}}
            & \includegraphics[width=0.28\textwidth]{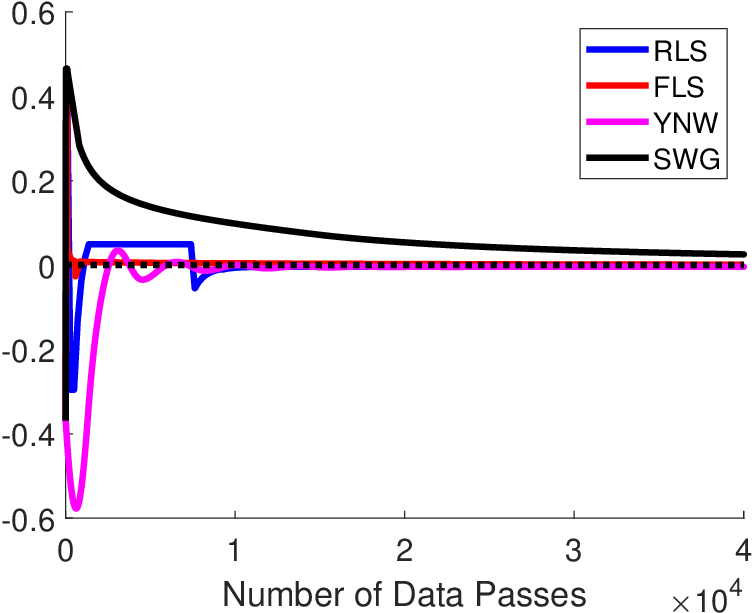}
            & \includegraphics[width=0.28\textwidth]{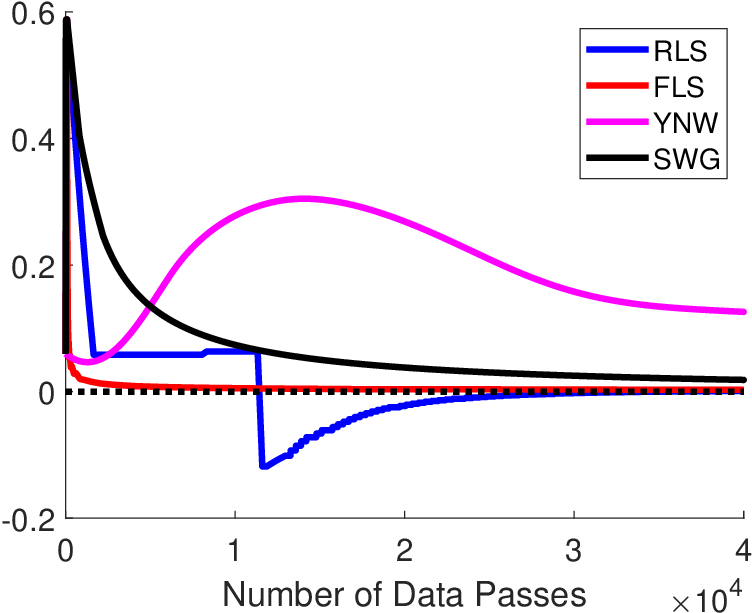}
            & \includegraphics[width=0.28\textwidth]{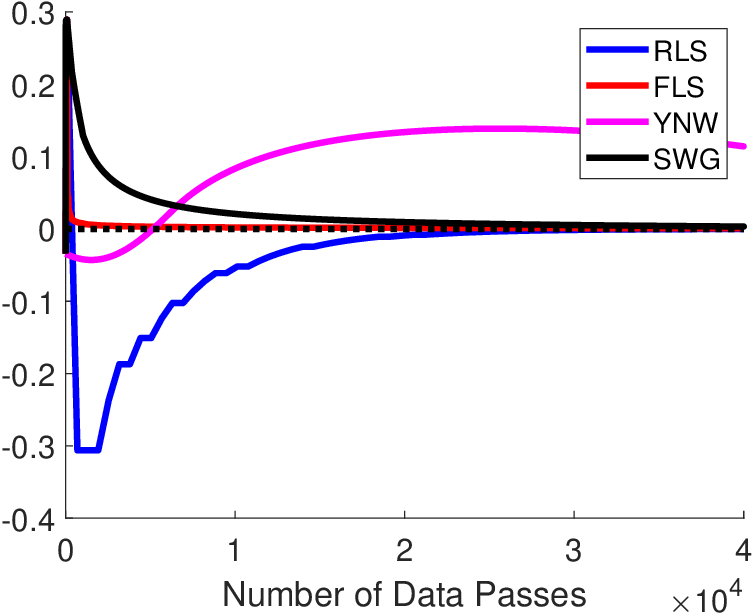} \\
            \raisebox{8ex}{\small\rotatebox[origin=c]{90}{$ \max\limits_{i=1,2}\{f_i(\bx)\}$}}
            & \includegraphics[width=0.28\textwidth]{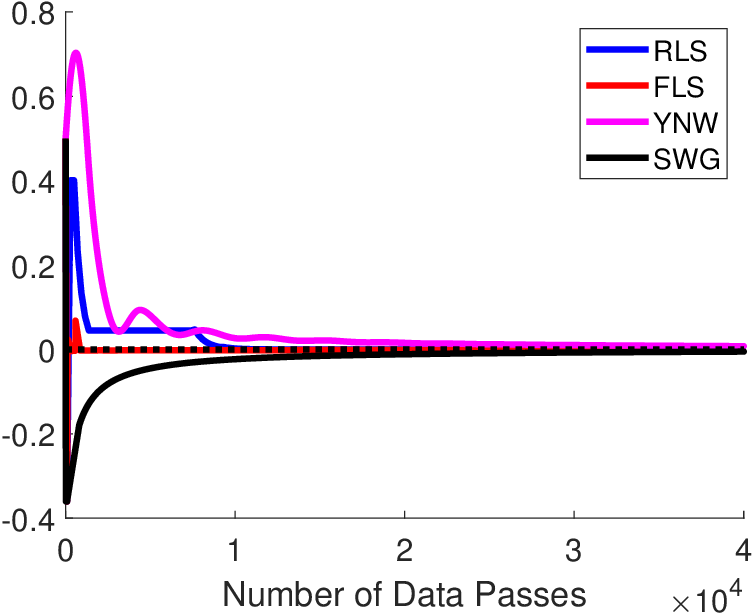}
            & \includegraphics[width=0.28\textwidth]{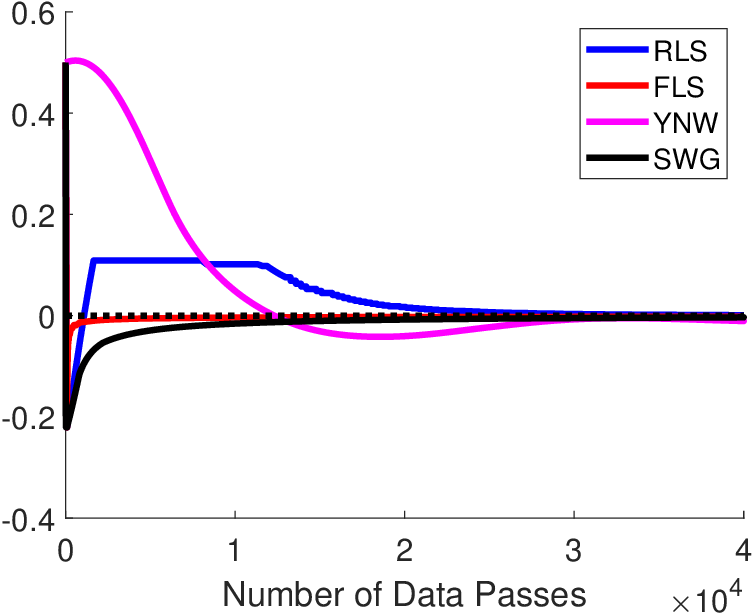}
            & \includegraphics[width=0.28\textwidth]{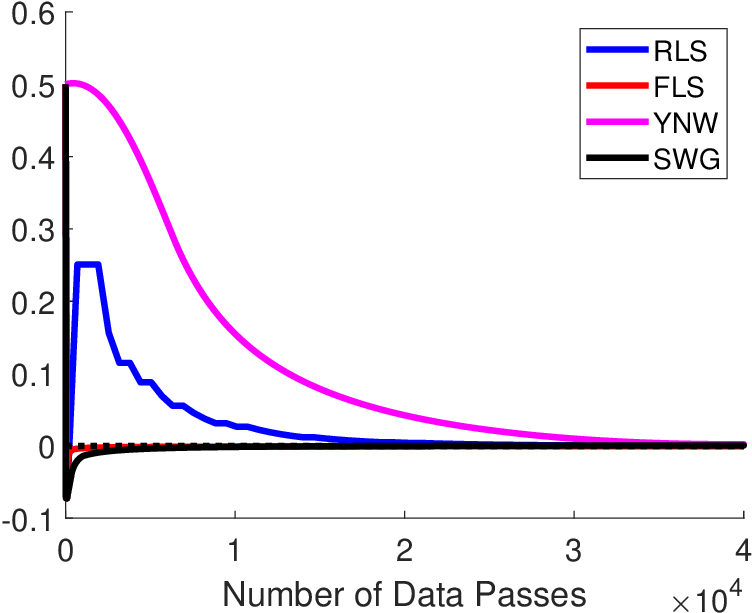} \\
        \end{tabular}%
    }
    \vspace{2ex}
    \caption{Performance of RLS, FLS, YNW, and SWG for solving non-smooth Neyman–Pearson classification problems with Type I error control on larger-scale datasets.}
    \label{fig:NP_nonsmooth}
\end{figure}

\noindent {\bf Results:} Figure~\ref{fig:NP_nonsmooth} illustrates the performance of RLS, FLS, YNW, and SWG over varying data passes. The axes' meanings remain consistent with Figures \ref{fig:fairness} and \ref{fig:fairness_smooth}. 
The plots illustrate that RLS achieves the fastest reduction of the metric $P(\bx;f^*)$ compared to YNW and SWG across all datasets. Although FLS initially reduces $P(\bx;f^*)$ more quickly when $P(\bx;f^*) > 10^{-2}$, RLS outperforms FLS in the high-accuracy regime (i.e., when $P(\bx;f^*)\leq 10^{-2}$), ultimately achieving significantly smaller values of $P(\bx;f^*)$. Moreover, while YNW exhibits substantial optimality gaps on the rcv1.binary and real-sim datasets, both FLS and SWG maintain feasibility throughout the iterations; among these, FLS demonstrates a faster reduction of the optimality gap than SWG. These numerical results, obtained on larger-scale classification problems, confirm that RLS is highly effective in simultaneously reducing the optimality gap and constraint violation, thereby yielding an accurate and feasible solution.

\section{Conclusion}\label{sec:conc}
We develop an adaptive level-set method that is both parameter-free and projection-free for constrained convex optimization, that is, it accelerates under the error bound condition and does not require knowledge of unknown parameters or challenging projections for its execution. This method finds an $\epsilon$-optimal and $\epsilon$-feasible solution by considering a sequence of level-set subproblems that are solved in parallel using standard subgradient oracles with simple updates. These oracles restart based on objective function progress and communicate information between them upon each restart. We show that the iteration complexity of our restarting level set method is only worse by a log-factor in both the smooth and non-smooth settings compared to existing accelerated methods for constrained convex optimization under the error bound condition, all of which rely on either unknown parameters or sophisticated oracles (e.g., involving difficult projection or exact line search). Numerical experiments show that the proposed method exhibits promising performance relative to benchmarks on constrained classification problems. 

While slightly beyond the primary scope of this paper, we explore an extension of our analysis in Appendix~\ref{Appendix:LinearEquality}. In this extension, we showcase how our methodology can adapt to scenarios where linear equality constraints are integrated into \eqref{eq:gco}. It's worth noting that in this extension, our adaptive algorithm is no longer parameter-free, as it requires knowledge of parameter $M$, an upper bound on the norm of sub-gradients. Despite this requirement, given that $M$ is readily computable in many applications, we believe this extension could be of interest to certain readers.

\appendix \label{sec:appendix}
\section{An Example of Error Bound Condition with Exponent $d > 2$}\label{Appendix:Example}
To demonstrate that the error bound condition can hold with exponent $ d > 2 $, we present an illustrative example adapted from \citet{li2015new}. Consider a constrained convex polynomial optimization problem in which the feasible region $\Xc$ is given by
\[
\Xc := \left\{ \bx \in \mathbb{R}^n \,\middle|\, h_i(\bx) \leq 0, \; i = 1,\ldots, \ell \right\},
\]
where each \( h_i \) is a convex polynomial of degree at most \( q \), and suppose that \( \Xc \) is compact. Let each objective and constraint function \( f_i \) for \( i = 0, \ldots, m \) also be a convex polynomial of degree at most \( q \). The optimization problem~\eqref{eq:gco} becomes
\[
\min_{\bx \in \Xc} \left\{ f_0(\bx) \,\middle|\, f_i(\bx) \leq 0, \; i = 1, \ldots, m \right\}.
\]
The corresponding optimal solution set is
\[
\Xc^* := \left\{ \bx \in \Xc \,\middle|\, \max\{f(\bx) - f^*, g(\bx), h(\bx)\} \leq 0 \right\},
\]
where \( f(\bx) := f_0(\bx) \), \( g(\bx) := \max_{i=1,\ldots,m} f_i(\bx) \), and \( h(\bx) := \max_{i=1,\ldots,\ell} h_i(\bx) \).

Under these assumptions, Corollary~3.8 in \citet{li2015new} implies that there exists a constant \( c > 0 \) such that for all \( \bx \in \Xc \),
\begin{align*}
\operatorname{dist}(\bx, \Xc^*) &\leq c \left( [f_0(\bx) - f^*]_+ + \sum_{i=1}^m [f_i(\bx)]_+ + \sum_{i=1}^\ell [h_i(\bx)]_+ \right)^\tau \\
&\leq c (m+1)^\tau \left( \max\{f(\bx) - f^*, g(\bx)\} \right)^\tau,
\end{align*}
where \( [\cdot]_+ := \max\{ \cdot, 0 \} \) and
\[
\tau := \max\left\{ \frac{1}{(q + 1)(3q)^{n + m + \ell}}, \frac{1}{q (6q - 3)^{n + m + \ell}} \right\}.
\]

This demonstrates that the EBC equation \eqref{eq:errorbound} holds with exponent $d = 1/\tau > 2$, demonstrating that values of $d > 2$ can naturally arise even in well-structured convex programs.

\section{Proofs}\label{Appendix:Proofs}

This section contains the proof of lemmas~\ref{lemma:approximate_condition_number} and \ref{lem:CInondecreasing} and presents Lemma~\ref{lem:BoundDkk'} and Lemma~\ref{lem:BoundDk'k'} that formally prove the inequalities \eqref{eq:Dkk'bound} and \eqref{eq:Dk'k'bound} used in the proof of Theorem~\ref{thm:correctnessAndComplexity}.

\emph{\bf Proof of Lemma~\ref{lemma:approximate_condition_number}:}
We prove this lemma by showing	$\tilde r>f^*$ and $0<\tilde\theta(r) \leq \theta$ for any $r<f^*$.
Since $\tilde{\bx}\in\mathcal{X}$ and $g(\tilde{\bx}) < 0$, by the definitions of $\tilde r$ and $H(\cdot)$, we have
\begin{eqnarray}
\label{eq:Hhatr}
	H(\tilde r)\leq \max\{f(\tilde{\bx})-\tilde r,g(\tilde{\bx})\}=g(\tilde{\bx})<0.
\end{eqnarray}
By property 3 in Lemma~\ref{lemma:levelset_properties}, we must have $\tilde r>f^*$ which indicates  $\tilde\theta(r) > 0$ for any $r<f^*$. It follows from the convexity of $H(r)$ and definition of $\theta$ in \eqref{eq:conditionnumber} that 
\begin{eqnarray}
	\label{eq:Hlinearize}
	H(r)\geq H(f^*) -\theta(r-f^*)= -\theta(r-f^*),\quad\forall r.
\end{eqnarray}
This further implies that, for $r<f^*$, 
$
	\tilde\theta(r) = \frac{-g(\tilde{\bx})}{\tilde{r}-r} \leq \frac{-H(\tilde r)}{\tilde{r}-r}\leq \frac{\theta(\tilde{r}-f^*)}{\tilde{r}-r} < \frac{\theta(\tilde{r}-r)}{\tilde{r}-r} = \theta,
$
where the first inequality is from \eqref{eq:Hhatr}, the second from \eqref{eq:Hlinearize} with $r=\tilde r$, and the last from $r < f^* < \tilde{r}$. 		\looseness = -1
\hfill\BlackBox

In the proof of the lemmas below we require the result of Proposition~\ref{prop:rplusP}.

\begin{proposition}\label{prop:rplusP}
The function $\dfrac{\alpha}{2}\Pc(\bx; r) + r$ is an increasing function in $r$ for any $\bx\in\mathcal{X}$.
\end{proposition}
\proof{} The proof of this proposition follows from the convexity of $\frac{\alpha}{2}\Pc(\bx; r)$ in $r$. In particular, let $\xi_r$ denote the sub-gradient of $\frac{\alpha}{2}\Pc(\cdot; r)$ with respect to $r$. It is straightforward to see $\xi_r\in [-1, 0]$ for any $r$. Hence, from convexity of $\frac{\alpha}{2}\Pc(\bx; r)$ in $r$, it follows that for $r \geq r'$, we have $\frac{\alpha}{2}\Pc(\bx; r)\geq \frac{\alpha}{2}\Pc(\bx; r') +\xi_r(r - r') \geq \frac{\alpha}{2}\Pc(\bx; r') -1 \cdot ( r - r')$ which indicates 
$
\frac{\alpha}{2}\Pc(\bx; r) + r \geq \frac{\alpha}{2}\Pc(\bx; r') + r'.
$
\hfill\BlackBox

For the remaining results included in this section, we need to introduce slightly different notations representing $\bx_k^{(t_k)}$ and $r_k$ to account for the number of restarts in each $\mathtt{fom}$ instance. More precisely, we use the notation $\bx_k^{(l_k, t_k)}$ and $r_k^{(l_k)}$ to respectively represent the solution in the $t_k$th iteration of $\mathtt{fom}_k$ after $l_k$ restarts and the level parameter in $\mathtt{fom}_k$ after $l_k$ restarts.

Let $k'$ be the smallest index at which  the inequalities	$\Pc(\bx_{k'}^{(l_{k'},t_{k'})};r_{k'}^{(l_{k'})})\leq B \Pc(\bx_{k'}^{(l_{k'},0)};r_{k'}^{(l_{k'})})$ and $\Pc(\bx_{k'}^{(l_{k'}, 0)}; r_{k'}^{(l_{k'})})\geq 0$ hold (i.e. $k'$ is the index found in Line 7 of Algorithm~\ref{algo:RestartingLS}).
RLS then initiates a restart at this index and also restarts all instances whose indices are greater than $k'$. Suppose this is the $l_k$th restart at $\mathtt{fom}_k$ for $k\geq k'$. To understand how level parameters change after each restart, we reply on the following lemma where we show $r_k^{(l_k-1)}\geq r_k^{(l_k)}$ for each $k\geq k'$. 
	\begin{lemma}
		\label{lem:gapchange}
		Suppose $\alpha$ and $B$ are such that $0 < \alpha < B< 1$. Let $k'$ be the index found in Line 7 of Algorithm~\ref{algo:RestartingLS}. In addition, assume the level parameters $r_k^{(l_k)}$ for $k\geq k'$ are updated as in Line 14, i.e., $r_{k}^{(l_{k})} = r_{k-1}^{(l_{k-1})} + \frac{\alpha}{2}\Pc(\bx_{k-1}^{(l_{k-1},0)}; r_{k-1}^{(l_{k-1})})$ for $k\geq k'+1$. We then have 
		\begin{equation}\label{eq:r.minus.r.bound.k'}
		r_{k'}^{(l_{k'}-1)}-r_{k'}^{(l_{k'})}=0,\end{equation} 
		and 
\begin{equation}\label{eq:r.minus.r.bound}
r_k^{(l_k-1)}-r_k^{(l_k)}\geq 0\quad \text{for } k\geq k'+1.
\end{equation}
		
	\end{lemma}
	\proof{}
		The equation $r_{k'}^{(l_{k'})} = r_{k'}^{(l_{k'}-1)}$ directly follows from the steps of Algorithm~\ref{algo:RestartingLS}. In particular, the level parameters are only updated for $k\geq k'+1$ after each restart at index $k'$. Recall that the index $k'$ is selected such that $\Pc(\bx_{k'}^{(l_{k'}-1, 0)}; r_{k'}^{(l_{k'}-1)})\geq 0$ and 
$
\Pc(\bx_{k'}^{(l_{k'}-1,t_{k'})}; r_{k'}^{(l_{k'}-1)}) \leq B \Pc(\bx_{k'}^{(l_{k'}-1, 0)}; r_{k'}^{(l_{k'}-1)}).
$
Since $\bx_{k'}^{(l_{k'}, 0)} = \argmin_{\bx\in\left\{\bx_{k'}^{(l_{k'}-1,t_{k'})}, \bx_1^{(l_1-1, 0)}, \ldots, \bx_K^{(l_K-1, 0)}\right\}} \Pc(\bx; r_{k'}^{(l_{k'}-1)})$ and $r_{k'}^{(l_{k'})} = r_{k'}^{(l_{k'}-1)}$, we get 
\begin{equation}\label{eq:rescaling.condition.k'}
\Pc(\bx_{k'}^{(l_{k'},0)}; r_{k'}^{(l_{k'})}) \leq B \Pc(\bx_{k'}^{(l_{k'}-1, 0)}; r_{k'}^{(l_{k'}-1)}).
\end{equation}
Now let's consider the index $k'+1$:\looseness = -1
		\begin{eqnarray}
			r_{k'+1}^{(l_{k'+1}-1)}-r_{k'+1}^{(l_{k'+1})}  &=&  r_{k'}^{(l_{k'}-1)}+\dfrac{\alpha}{2} \Pc(\bx_{k'}^{(l_{k'}-1,0)}; r_{k'}^{(l_{k'}-1)})-r_{k'}^{(l_{k'})} - \dfrac{\alpha}{2} \Pc(\bx_{k'}^{(l_{k'},0)}; r_{k'}^{(l_{k'})})\nonumber\\ 
		&=&  \dfrac{\alpha}{2} \Pc(\bx_{k'}^{(l_{k'}-1,0)}; r_{k'}^{(l_{k'}-1)})-\dfrac{\alpha}{2} \Pc(\bx_{k'}^{(l_{k'},0)}; r_{k'}^{(l_{k'})})\nonumber\\[2pt]
		&\geq&  \dfrac{\alpha}{2} \Pc(\bx_{k'}^{(l_{k'}-1,0)}; r_{k'}^{(l_{k'}-1)})-\dfrac{\alpha}{2} B \Pc(\bx_{k'}^{(l_{k'}-1,0)}; r_{k'}^{(l_{k'}-1)})\nonumber \\[2pt]
			&=& \dfrac{\alpha}{2} (1-B)\Pc(\bx_{k'}^{(l_{k'}-1,0)}; r_{k'}^{(l_{k'}-1)})\geq 0 \label{ineq:r_first}
		\end{eqnarray}
		where we used the definitions of $r_{k'+1}$ to obtain the first equality; \eqref{eq:r.minus.r.bound.k'} for the second; and \eqref{eq:rescaling.condition.k'} for the first inequality and the inequality $\Pc(\bx_{k'}^{(l_{k'}-1, 0)};r_{k'}^{(l_{k'}-1)})\geq 0$ for the last;  Similarly for $k\geq k'+2$, we show 
		\begin{align}\label{ineq:r_recursive}
		\nonumber
		r_{k}^{(l_{k}-1)}-r_{k}^{(l_{k})}  &=  r_{k-1}^{(l_{k-1}-1)}+\dfrac{\alpha}{2} \Pc(\bx_{k-1}^{(l_{k-1}-1,0)}; r_{k-1}^{(l_{k-1}-1)})-r_{k-1}^{(l_{k-1})} - \dfrac{\alpha}{2} \Pc(\bx_{k-1}^{(l_{k-1},0)}; r_{k-1}^{(l_{k-1})})\nonumber\\
		&\geq  r_{k-1}^{(l_{k-1}-1)} + \dfrac{\alpha}{2} \Pc(\bx_{k-1}^{(l_{k-1}-1,0)}; r_{k-1}^{(l_{k-1}-1)})- \left(r_{k-1}^{(l_{k-1})}+\dfrac{\alpha}{2} \Pc(\bx_{k-1}^{(l_{k-1}-1,0)}; r_{k-1}^{(l_{k-1})}) \right)\geq  0,\nonumber\looseness = -1
		\end{align}
where the first inequality follows from Line 13 of Algorithm~\ref{algo:RestartingLS} which indicates that $$\Pc(\bx_{k-1}^{(l_{k-1},0)}; r_{k-1}^{(l_{k-1})}) \leq \Pc(\bx_{k-1}^{(l_{k-1}-1,0)}; r_{k-1}^{(l_{k-1})}).$$ The last inequality form Proposition~\ref{prop:rplusP} and the fact that $ r_{k-1}^{(l_{k-1}-1)}\geq r_{k-1}^{(l_{k-1})}$ for $k\geq k'+2$ which can be easily verified by \eqref{ineq:r_first}. 
\hfill\BlackBox

\emph{\bf Proof of Lemma~\ref{lem:CInondecreasing}:} 
To prove this lemma, we claim that if the inequality~\eqref{eq:xkquality} ($\alpha\Pc(\bx_k^{(0)}; r_k) \leq f^* - r_k$) is satisfied at some $k$, it will continue to be satisfied after each restart. By Theorem~\ref{thm:keythm}, it will then follows that none of the indices that are smaller than the current critical index can become a critical index after the restarts. Therefore, the critical index cannot decrease throughout the algorithm. Our claim can be proved by defining the quantity $V_k\left(l_k\right)$, $k = 0,1, 2, \ldots$ as follows:
\begin{equation}\label{eq:V}
 V_k\left(l_k\right):= f^* -r_k^{(l_k)} -\alpha \Pc(\bx_k^{(l_k,0)}; r_k^{(l_k)}).
\end{equation}
We show $V_k(l_k)$ is a non-decreasing function of $l_k$ when the inequality~\eqref{eq:xkquality} holds at index $k$. More precisely, we show the inequality $V_k(l_k)\geq V_k(l_k - 1)$ holds for such $k$. Therefore, assuming \eqref{eq:xkquality} is satisfied in $l_k -1$ restarts, it will continue to hold after next restart since $V_k(l_k)\geq V_k(l_k - 1) \geq 0$.

Fix an index $k\in[0,\tilde K]$ for which we have $\alpha\Pc(\bx_k^{(l_k-1,0)}; r_k^{(l_k-1)}) < f^* - r_k^{(l_k-1)}$. The restarts that occur at larger indices than $k$ do not change $V_k(l_k-1)$ since $\bx_k^{(l_k-1, 0)}$ and $r_k^{(l_k-1)}$ remain the same. Therefore, we only focus on the restarts that either happen (i) at index $k$ or (ii) any smaller index than $k$. 

Case (i): When the restart occurs at index $k$, it means we have found a solution $\bx_{k}^{(l_k-1,t_k)}$ such that $\Pc(\bx_{k}^{(l_k-1,0)}; r_{k}^{(l_k-1)})\geq 0$ and
		$\Pc(\bx_{k}^{(l_k-1,t_k)}; r_{k}^{(l_k-1)}) \leq B \Pc(\bx_{k}^{(l_k-1,0)}; r_{k}^{(l_k-1)}).$
		 Since $\bx_{k}^{(l_k, 0)} = \bx_{k}^{(l_k-1,t_k)}$ and $r_{k}^{(l_k)} = r_{k}^{(l_k-1)}$ which follows from the restarting steps, we get	
		 $\Pc(\bx_{k}^{(l_k,0)}; r_{k}^{(l_k)}) \leq B \Pc(\bx_{k}^{(l_k-1,0)}; r_{k}^{(l_k-1)}).$
		 Using this inequality and the fact that $r_{k}^{(l_k)} = r_{k}^{(l_k-1)}$, it is straightforward to see $V_k(l_k)\geq V_k(l_k-1) + \alpha(1-B) \Pc(\bx_{k}^{(l_k-1,0)}; r_k^{(l_k-1)})$. This inequality further implies that
	\begin{align}
	\nonumber
	V_k(l_k)&\geq V_k(l_k-1) + \alpha(1-B) \Pc(\bx_{k}^{(l_k-1,0)}; r_k^{(l_k-1)}) 
    \geq V_k(l_k-1)+ \alpha(1-B)\theta\left(f^*-r_k^{(l_k-1)}\right)\\\nonumber
	&>  V_k(l_k-1)+ \alpha(1-B)\theta(f^*-r_k^{(l_k-1)} - \alpha\Pc(\bx_k^{(l_k-1,0)}; r_k^{(l_k-1)})) \\ 
	&= \left(1+\alpha(1-B)\theta\right)V_k(l_k-1) > V_k(l_k-1),\label{eq:changeVatk}
	\end{align}
	 where the second inequality follows from $0\leq \alpha\Pc(\bx_k^{(l_k-1,0)}; r_k^{(l_k-1)}) \leq f^* - r_k^{(l_k-1)}$ and Lemma~\ref{theorem:levelset_method_theta} that shows $\theta \leq \frac{H(r_k^{(l_k-1)})}{f^* - r_k^{(l_k-1)}} \leq \frac{\Pc(\bx_k^{(l_k-1,0)}; r_k^{(l_k-1)})}{f^*-r_k^{(l_k-1)}}$, and the third from the inequality $\Pc(\bx_k^{(l_k-1,0)}; r_k^{(l_k-1)})\geq 0$ which holds since the restart has occurred at index $k$.

Case (ii): In this case we have 
	\begin{align}
	\label{eq:changeVlessthank}
	V_k(l_k)=f^* - r_k^{(l_k)} - \alpha\Pc(\bx_k^{(l_k,0)}; r_k^{(l_k)})
	\geq f^* - r_k^{(l_k-1)} - \alpha \Pc(\bx_k^{(l_k-1,0)}; r_k^{(l_k-1)})=V_k(l_k-1).
	\end{align}
	To obtain the above inequality we used Line 13 of Algorithm~\ref{algo:RestartingLS} which indicates that $\Pc(\bx_k^{(l_k,0)}; r_k^{(l_k)}) \leq \Pc(\bx_k^{(l_k-1,0)}; r_k^{(l_k)})$. In addition, following the proof of  Proposition~\ref{prop:rplusP}, one can show that the function $f^* - r - \alpha \Pc(\bx; r)$ is decreasing in $r$ for any $\bx\in\X$. Therefore, the inequality~\eqref{eq:changeVlessthank} holds since $r_k^{(l_k-1)} \geq r_k^{(l_k)}$ which follows from Lemma~\ref{lem:gapchange}. Considering \eqref{eq:changeVatk} and \eqref{eq:changeVlessthank}, our proof is complete. \hfill\BlackBox

In the following lemmas we prove the inequalities \eqref{eq:Dkk'bound} and \eqref{eq:Dk'k'bound} used in the proof of Theorem~\ref{thm:correctnessAndComplexity} in Section~\ref{section:main_algorithm}.

\begin{lemma}\label{lem:BoundDkk'}
Suppose $\alpha$, $B$, and $\epsilon$ are such that $0<\alpha<B<1$ and $\epsilon > 0$.
Consider an index $k$ and $\tilde K$ as defined in~\eqref{eq:tildeK} for $r = r_{\text{ini}}$. Let $D_k$ denote the total number of desirable restarts starting at $\mathtt{fom}_k$ with $k < k^*$ until an $\epsilon$-optimal and $\epsilon$-feasible solution in found, assuming $k^*$ is the critical index at the time of restart. Then we have
\begin{align*}
D_k = \sum_{k^*= k+1}^{\tilde K} D_{kk^*} \leq \ln\left(\dfrac{4-\alpha\theta}{3\alpha\theta}\right)/\ln\left(1+\alpha(1-B)\theta/4\right),
\end{align*}
 where $D_{kk^*}$ is the total number of desirable restarts starting at $\mathtt{fom}_k$ with $k < k^*$.
\end{lemma}
\proof{}
Using the definitions of $D_k$ and $D_{kk^*}$, it is easy to see that $D_k = \sum_{k^*= k+1}^{\tilde K} D_{kk^*}$ because by Lemma~\ref{lem:CInondecreasing}, the critical index $k^*$ is non-decreasing and can vary between $k+1$ and $\tilde K$. Let's define
\begin{equation*}
V'_k(l_k) := \dfrac{\Pc(\bx_k^{(l_k,0)}; r_k^{(l_k)})}{f^* - r_k^{(l_k)} - \frac{\alpha}{4}\Pc(\bx_k^{(l_k,0)}; r_k^{(l_k)})},
\end{equation*}
where the denominator is a modification of the quantity $V_k(l_k)$ in \eqref{eq:V}. 

We show $V'_k(l_k)$ is a non-increasing function of $l_k$. More precisely, we show the inequality $V'_k(l_k)\leq V'_k(l_k - 1)$ holds for each $k < k^*$.  The restarts that occur at larger indices than $k$ do not change $V'_k(l_k-1)$ since $\bx_k^{(l_k-1, 0)}$ and $r_k^{(l_k-1)}$ remain the same. Therefore, we only focus on the restarts that either happen at index $k$ or any smaller index than $k$.

Since $k < k^*$, Theorem~\ref{thm:keythm} indicates that for all $l_k \geq 0$ we have $r_k^{(l_k)} < f^*$ and the inequality \eqref{eq:xkquality} holds at $k$. Hence, we can easily verify that $\Pc(\bx_k^{(l_k,0)}; r_k^{(l_k)})\geq H(r_k^{(l_k)}) > 0$, $V'_k(l_k) > 0$, $r_k^{(l_k)} \geq r_0^{(l_0)}$, and 

\begin{equation}\label{eq:Pdecreases}\Pc(\bx_k^{(l_k,0)}; r_k^{(l_k)}) \leq \Pc(\bx_{k-1}^{(l_{k-1},0)}; r_k^{(l_k)}) \leq \Pc(\bx_{k-1}^{(l_{k-1},0)}; r_{k-1}^{(l_{k-1})}).
\end{equation} 
The above inequality follows from Line 13 of Algorithm \ref{algo:RestartingLS}, Lemma~\ref{lem:gapchange}, and the fact that the function $\Pc(\cdot; r)$ is decreasing in $r$.
Following the same proof as in the proof of Lemma~\ref{lem:CInondecreasing}, one can also show 
\begin{equation}\label{eq:V'4restartatk}
f^* - r_k^{(l_k)} - \dfrac{\alpha}{4}\Pc(\bx_k^{(l_k,0)}; r_k^{(l_k)}) \geq \left(1+\alpha(1-B)\theta/4\right)\left(f^* - r_k^{(l_{k-1})} - \dfrac{\alpha}{4}\Pc\left(\bx_k^{(l_{k-1},0)}; r_k^{(l_{k-1})}\right)\right),
\end{equation}
if the restarts occurs at index $k$ and 

\begin{equation}\label{eq:V'4restartlessthank}
f^* - r_k^{(l_k)} - \dfrac{\alpha}{4}\Pc(\bx_k^{(l_k,0)}; r_k^{(l_k)}) \geq f^* - r_k^{(l_{k-1})} - \dfrac{\alpha}{4}\Pc(\bx_k^{(l_{k-1},0)}; r_k^{(l_{k-1})}),
\end{equation}
if the restarts occurs at an index smaller than $k$. Combining \eqref{eq:V'4restartatk} and \eqref{eq:V'4restartlessthank} with the inequality \eqref{eq:Pdecreases}, we can respectively show
$$
\begin{cases}
    V'_{k}(l_k)\leq \dfrac{1}{\left(1+\alpha(1-B)\theta/4\right)} V'_{k}(l_{k} -1)       & \quad \text{if the restarts occurs at index $k$, }\\[10pt]
    V'_{k}(l_k)\leq V'_{k}(l_{k} -1)   & \quad \text{if the restarts occurs at an index smaller than $k$}.
  \end{cases}
$$

Suppose $\hat l_k$ denotes the smallest restarting iteration at which the index $k$ becomes smaller than a critical index. By recursively applying the definition of $V'$, we obtain 
\begin{eqnarray}
V'_k(l_k) &\leq&\dfrac{\Pc(\bx_k^{(\hat l_k,0)}; r_k^{(\hat l_k)})}{\left(1+\alpha(1-B)\theta/4\right)^{D_k}\left(f^* - r_k^{(\hat l_k)} -\dfrac{\alpha}{4}\Pc(\bx_k^{(\hat l_k,0)}; r_k^{(\hat l_k)})\right)}\nonumber\\
&\leq& \dfrac{1}{\left(1+\alpha(1-B)\theta/4\right)^{D_k}}\cdot\dfrac{(1/\alpha)(f^*-r_k^{(\hat l_k)})}{(3/4)(f^*-r_k^{(\hat l_k)})} = \dfrac{1}{\left(1+\alpha(1-B)\theta/4\right)^{D_k}}\cdot\dfrac{4}{3\alpha},\label{eq:upperboundonV'}
\end{eqnarray}
where the second inequality holds because for any $k< k^*$, we have $\Pc(\bx_k^{(\hat l_k,0)}; r_k^{(\hat l_k)})\leq (f^*-r_k^{(\hat l_k)})/\alpha$.

In addition since $\Pc(\bx_k^{(l_k,0)}; r_k^{(l_k)}) \geq H(r_k^{(l_k)})\geq \theta (f^* - r_k^{(l_k)})$ which holds by Lemma~\ref{theorem:levelset_method_theta}, we can easily verify that 
\begin{equation}\label{eq:lowerboundonV'}
V'_k(l_k) \geq \dfrac{\theta}{1-\frac{\alpha\theta}{4}}.
\end{equation}
Then from the inequalities \eqref{eq:upperboundonV'}  and \eqref{eq:lowerboundonV'} it follows that
$
\dfrac{\theta}{1-\frac{\alpha\theta}{4}}\leq \dfrac{1}{\left(1+\alpha(1-B)\theta/4\right)^{D_k}}\cdot\dfrac{4}{3\alpha},
$
We can then obtain our desirable bound on $D_k$ by taking a logarithmic transformation on both sides of the above inequality and organizing terms.
\hfill\BlackBox
	
\begin{lemma}\label{lem:BoundDk'k'} Suppose $\alpha$, $B$, and $\epsilon$ are such that $0<\alpha<B<1$ and $\epsilon > 0$. Consider a critical index $k^*$ and $\tilde K$ as defined in~\eqref{eq:tildeK} for $r = r_{\text{ini}}$. Let $D_{k^*k^*}$ denote the total number of desirable restarts starting at $\mathtt{fom}_{k^*}$ before the critical index increases or an $\epsilon$-feasible and $\epsilon$-optimal solution is found. Then
		$D_{k^*k^*}\leq \ln\left(2/\alpha\theta\right)/\ln(1/B).$
	\end{lemma}
	
	\proof{}
	Suppose $k^*$ is the critical index while there is no $\epsilon$-feasible and $\epsilon$-optimal solution is found. Consider $k\leq k^*$. Using the updates of $\bx_{k}$ and $r_{k}$ in Algorithm \ref{algo:RestartingLS} and the fact that the function $\Pc(\cdot; r_k^{(l_k)})$ is non-increasing in $r_k^{(l_k)}$ we can obtain
	\begin{eqnarray}\label{eq:boundonP}
	\Pc(\bx_{k}^{(l_{k},0)}; r_{k}^{(l_{k})})\leq \Pc(\bx_{k-1}^{(l_{k-1},0)}; r_{k}^{(l_{k})}) \leq \Pc(\bx_{k-1}^{(l_{k-1},0)}; r_{k-1}^{(l_{k-1})}).
	\end{eqnarray}
	Notice that we used the inequality $r_{k-1}^{(l_{k-1})}\leq r_{k}^{(l_{k})}$ from Theorem~\ref{thm:keythm} to obtain the second inequality above. 
	In addition, 
	\begin{eqnarray}\label{eq:boundonf*-r}
	f^* - r_{k}^{(l_k)} = f^* - r_{k-1}^{(l_{k-1})} -\dfrac{\alpha}{2}\Pc(\bx_{k-1}^{(l_{k-1},0)}; r_{k-1}^{(l_{k-1})})\geq \dfrac{1}{2}\left(f^* - r_{k-1}^{(l_{k-1})}\right),	
 \end{eqnarray}
 where the inequality holds since $k \leq k^*$ and hence Theorem \ref{thm:keythm} indicates that {\small $\Pc(\bx_{k-1}^{(l_{k-1},0)}; r_{k-1}^{(l_{k-1})}) \leq (f^* - r_{k-1}^{(l_{k-1})})/\alpha$.}
Define
		$
		U_{k^*}(l_{k^*}):=\dfrac{\Pc(\bx_{k^*}^{(l_{k^*},0)}; r_{k^*}^{(l_{k^*})})}{f^* - r_{k^*}^{(l_{k^*})}}.
		$ 
We show this quantity is non-increasing at each desirable restart. Especially, it reduces with a factor of $B$ if desirable restarts occur at $\mathtt{fom}_{k^*}$. More precisely, we will prove 
$$
\begin{cases}
    U_{k^*}(l_{k^*})\leq U_{k^*}(l_{k^*} -1)       & \quad \text{if $\mathtt{fom}_{k^*}$ instance is updated by a restart at $\mathtt{fom}_{k}$ with $k<k^*$, }\\[10pt]
    U_{k^*}(l_{k^*})\leq B \cdot U_{k^*}(l_{k^*} -1)  & \quad \text{if $\mathtt{fom}_{k^*}$ instance is updated by a restart at $\mathtt{fom}_{k^*}$}.
  \end{cases}
$$

Finally, by providing a lower bound on $U_{k^*}(l_{k^*})$ for any $l_{k^*}$, we show that this quantity can only be reduced by a finite number of desirable restarts.

First, let's consider a case where $\bx_{k^*}^{(l_{k^*},0)}$ and $r_{k^*}^{(l_{k^*})}$ are generated by a restart at $\texttt{fom}_k$ with $k<k^*$. We then have 
		\begin{align}
		\label{eq:changeU1}
		U_{k^*}(l_{k^*}) &=\dfrac{\Pc(\bx_{k^*}^{(l_{k^*},0)}; r_{k^*}^{(l_{k^*})})}{f^* - r_{k^*}^{(l_{k^*})}}
		\leq\dfrac{\Pc(\bx_{k^*}^{(l_{k^*}-1,0)}; r_{k^*}^{(l_{k^*}-1)})+r_{k^*}^{(l_{k^*}-1)}-r_{k^*}^{(l_{k^*})} }{f^*-r_{k^*}^{(l_{k^*}-1)}+r_{k^*}^{(l_{k^*}-1)}-r_{k^*}^{(l_{k^*})}}\nonumber\\[10pt]
		&\leq\dfrac{\Pc(\bx_{k^*}^{(l_{k^*}-1,0)}; r_{k^*}^{(l_{k^*}-1)}) }{f^*-r_{k^*}^{(l_{k^*}-1)}}=U_{k^*}(l_{k^*}-1),
		\end{align}
		where the first inequality follows from the definition of $\Pc(\bx; r)$ and Line 13 of Algorithm~\ref{algo:RestartingLS}, and the second from $r_{k^*}^{(l_{k^*}-1)}-r_{k^*}^{(l_{k^*})}\geq 0$ (see Lemma~\ref{lem:gapchange}) and the fact that $k^*$ is a critical index and hence
		$
		\Pc(\bx_{k^*}^{(l_{k^*}-1,0)}, r_{k^*}^{(l_{k^*}-1)})\geq \frac{f^*-r_{k^*}^{(l_{k^*}-1)}}{\alpha}\geq f^*-r_{k^*}^{(l_{k^*}-1)}.
		$	
Now suppose  $\bx_{k^*}^{(l_{k^*},0)}$ and $r_{k^*}^{(l_{k^*})}$ are generated by a restart at $\texttt{fom}_{k^*}$, Since every desirable restart from $\mathtt{fom}_{k^*}$ happens when $\Pc(\bx_{k^*}^{(l_{k^*}-1,t_{k^*})}; r_{k^*}^{(l_{k^*}-1)}) \leq B \Pc(\bx_{k^*}^{(l_{k^*}-1,0)};r_{k^*}^{(l_{k^*})})$, we get 
	$
		U_{k^*}(l_{k^*})\leq\frac{\Pc(\bx_{k^*}^{(l_{k^*},0)}; r_{k^*}^{(l_{k^*})})}{f^* - r_{k^*}^{(l_{k^*})}} \leq \frac{\Pc(\bx_{k^*}^{(l_{k^*}-1,t_{k^*})}; r_{k^*}^{(l_{k^*}-1)})}{f^* - r_{k^*}^{(l_{k^*}-1)}}
		\leq\frac{B \Pc(\bx_{k^*}^{(l_{k^*}-1,0)}; r_{k^*}^{(l_{k^*} -1)}) }{f^*-r_{k^*}^{(l_{k^*}-1)}}
  =BU_{k^*}(l_{k^*}-1), 
	$
		where we used $r_{k^*}^{(l_{k^*})}=r_{k^*}^{(l_{k^*}-1)}$. Using this inequality and \eqref{eq:changeU1} recursively, we get
		\begin{align}
		\nonumber
		U_{k^*}(l_{k^*})=\dfrac{\Pc(\bx_{k^*}^{(l_{k^*},0)}; r_{k^*}^{(l_{k^*})})}{f^* - r_{k^*}^{(l_{k^*})}} \leq\dfrac{B^{D_{k^*k^*}}\Pc(\bx_{k^*}^{(\hat l_{k^*},0)}; r_{k^*}^{(\hat l_{k^*})})}{f^* - r_{k^*}^{(\hat l_{k^*})}}&\leq \dfrac{B^{D_{k^*k^*}}\Pc\left(\bx_{k^*}^{(\hat l_{k^*-1},0)}; r_{k^*-1}^{(\hat l_{k^*-1})}\right)}{\dfrac{1}{2}\left(f^* - r_{k^*-1}^{(\hat l_{k^*-1})}\right)}\\
  &\leq \dfrac{2B^{D_{k^*k^*}}}{\alpha}, \label{eq:changeU4}
		\end{align}
		where the first inequality holds since $r_{k^*}^{(\hat l_{k^*})} \geq r_{k^*}^{(l_{k^*})}$ which follows from Lemma~\ref{lem:gapchange}. The second holds by \eqref{eq:boundonP} and \eqref{eq:boundonf*-r} applied to $k= k^*$, and the third by Theorem \ref{thm:keythm} which indicates that $\Pc(\bx_{k^*-1}^{(\hat l_{k^*-1},0)}; r_{k^*-1}^{(\hat l_{k^*-1})}) \leq (f^* - r_{k^*-1}^{(\hat l_{k^*-1})})/\alpha$.

		Moreover, since $\Pc(\bx_{k^*}^{(l_{k^*},0)}; r_{k^*}^{(l_{k^*})}) \geq H(r_{k^*}^{(l_{k^*})})\geq \theta (f^* - r_{k^*}^{(l_{k^*})})$ which holds by Lemma~\ref{theorem:levelset_method_theta},  we get 
		$
		U_{k^*}(l_{k^*})=\frac{\Pc(\bx_{k^*}^{(l_{k^*},0)}; r_{k^*}^{(l_{k^*})})}{f^* - r_{k^*}^{(l_{k^*})}}
		\geq\frac{\theta(f^* - r_{k^*}^{(l_{k^*})})}{f^* - r_{k^*}^{(l_{k^*})}}=\theta.
		$
		This inequality and \eqref{eq:changeU4} together imply
		$
		\theta\leq \frac{2B^{D_{k^*k^*}}}{\alpha},
		$
		which results in the following upper bound on $D_{k^*k^*}$:
		$
		D_{k^*k^*}\leq \ln\left(2/\alpha\theta\right)/\ln(1/B).
		$\looseness = -1
\hfill\BlackBox

	\section{Linear Equality Constraints}\label{Appendix:LinearEquality}
	
As discussed in Section~\ref{section:main_algorithm}, the restarting level set method, described in Algorithm~\ref{algo:RestartingLS}, requires $K+1$ instances of the first-order method. Here, $K$ (defined in \eqref{eq:trueK}) depends on certain unknown parameters $f^*$ and $\theta$. Under Assumption~\ref{assume:strictfeasible}, we established a lower bound on the condition number $\theta$ and an upper bound on $f^*$.  Utilizing these established bounds, we were able to formulate an approximate value for $K$, termed $\tilde K$, seen in \eqref{eq:tildeK}. This $\tilde K$ is something we can calculate and apply in implementing Algorithm~\ref{algo:RestartingLS}. Notably, our subsequent analysis remains applicable even when Assumption~\ref{assume:strictfeasible} is not met exactly where linear equality constraints are in play. This situation renders Assumption~\ref{assume:strictfeasible} inapplicable due to the equality constraints. Despite this, we attempt to find a computable lower bound on $\theta$ in this context. It's important to highlight that, in this particular case, the lower bound, denoted as $\tilde\theta$, hinges on having knowledge of the parameter $M$. This does introduce a limitation. Suppose
\begin{eqnarray}
	f^* := 
	\min_{\bx \in \Xc} \left\{f(\bx):= f_0(\bx)\quad\text{s.t.}\quad g(\bx) := \max_{i=1,\dots,m} f_i(\bx)\leq 0, \bA\bx=\bb \right\},
	\label{eq:gco_linearequality}
\end{eqnarray}
where $f_i$ for $i=0,1,\dots,m$ are convex real-valued functions, $\Xc\subset\mathbb{R}^n$ is a closed convex set onto which the projection mapping is easy to compute,  $\bA\in\mathbb{R}^{l\times n}$ is a full row rank matrix of size $l\times n$ and $\bb\in\mathbb{R}^l$ is a $l$-dimensional vector. In this case, $P(\bx;r)$ is defined as
\begin{align}
	\label{preliminaries:Pfunction_linearequality}
	P(\bx;r):= \max \left\{ f(\bx) -r,g(\bx),\|\bA\bx-\bb\|_\infty\right\}.
\end{align}

In this section, we consider a modification of Assumption~\ref{assume:strictfeasible} as describe below as well as Assumption~\ref{assumptions:nonsmooth}.
\begin{assumption}
	There exists a computable solution  $\tilde{\bx} \in \mathrm{int}(\mathcal{X})$ such that $\bA\tilde{\bx}=\bb$ and $g(\tilde{\bx})< 0$. 
	\label{assume:strictfeasible_linearequality}
\end{assumption}

\begin{lemma}\label{lemma:approximate_condition_number_linearequality} Let $\tilde \bx$ be the feasible solution in Assumption \ref{assume:strictfeasible_linearequality} and $\theta$ be the condition measure defined in~\eqref{eq:conditionnumber}.  It holds that
	$$
	\theta\geq \left( 1+ \frac{ f(\tilde{\bx})-r_{\text{ini}}}{-g(\tilde{\bx}) } + \sqrt{l}\left\|(\bA\bA^\top)^{-1}\bA\right\|\left(\frac{M(f(\tilde{\bx})-r_{\text{ini}}-g(\tilde{\bx}))}{-g(\tilde{\bx}) }+ \frac{f(\tilde{\bx})-r_{\text{ini}}}{\mathrm{dist}(\tilde{\bx},\partial \mathcal{X})}\right)\right)^{-1}.
	$$
\end{lemma}
\proof{}
Let $\bx^*$ be the optimal solution of \eqref{eq:gco_linearequality}. Under Assumption~\ref{assume:strictfeasible_linearequality}, there exist Lagrangian multipliers $\lambda^*\geq 0$ and $\gamma^*\in\mathbb{R}^{l}$ such that $\lambda^* g(\bx^*) = 0$ and 
\begin{align}
	\label{eq:optcond1}
\bzt_f +\lambda^* \bzt_g+\bA^\top \gamma^*+ {\bu}^*=\mathbf{0},
\end{align}
where $\bzt_f\in\partial f(\bx^*)$, $\bzt_g\in\partial g(\bx^*)$, ${\bu}^*\in \mathcal{N}_\Xc(\bx^*)$ and $\mathcal{N}_\Xc(\bx^*)$ is the normal cone of $\Xc$ at $\bx^*$.
Let ${\bar\gamma}^+\in\mathbb{R}_+^{l}$ and ${\bar\gamma}^- \in\mathbb{R}_+^{l}$ are defined such that ${\bar\gamma}^+_i=\gamma^*_i$ and ${\bar\gamma}^-_i=0$ if $\gamma^*_i\geq 0$  and ${\bar\gamma}^+_i=0$ and ${\bar\gamma}^-_i=-\gamma^*_i$ if $\gamma^*_i< 0$ for $i=1,\dots,l$. Let $\mu := 1+\lambda^*+\|{\bar\gamma}^+\|_1+\|{\bar\gamma}^-\|_1$. Then we can normalize \eqref{eq:optcond1} as 
\begin{align}
	\label{eq:optcond2}
	\frac{1}{\mu}\left(\bzt_f +\lambda^* \bzt_g+\bA^\top {\bar\gamma}^+-\bA^\top {\bar\gamma}^-\right)+ {\bv}^*=\mathbf{0},
\end{align}
where ${\bv}^*=\dfrac{1}{\mu}\bu\in \mathcal{N}_\Xc(\bx^*)$. Let $\Delta$ denotes the non-negative probability simplex. From the definition of $H(r)$ it then follows that for any $r$, 
\begin{eqnarray*}
H(r)&=& \min_{\bx\in\Xc}\max \left\{ f(\bx) -r,g(\bx), \max(\bA\bx-\bb),\max(-\bA\bx+\bb)\right\}\\
&=& \min_{\bx\in\Xc}\max_{(\alpha, \lambda, \gamma^+, \gamma^-)\in \Delta} \left\{\alpha (f^*-r) + \lambda g(x) + \gamma^{+} (\bA\bx-\bb) + \gamma^{-} (-\bA\bx+\bb)\right\}\\
&\geq&\dfrac{1}{\mu}\min_{\bx\in\Xc} \left\{f(\bx^*)+\bzt_f^\top(\bx-\bx^*)-r+\lambda^*(g(\bx^*)+\bzt_g^\top(\bx-\bx^*))+(\bx-\bx^*)^\top\bA^\top\gamma^+\right.\\
&& \left.-(\bx-\bx^*)^\top\bA^\top\gamma^-\right\}
\geq \dfrac{1}{\mu} \left(f(\bx^*)-r+\lambda^*g(\bx^*)\right)
= H(f^*) -\dfrac{1}{\mu}\left(r-f^*\right),
\end{eqnarray*}
where the first inequality is established through the convexity of $f$ and $g$, along with the conditions $\bA\bx^* = \bb$ and $\frac{1}{\mu} \cdot (1, \lambda^*, {\bar\gamma}^+, {\bar\gamma}^-) \in\Delta$. The second inequality is a result of \eqref{eq:optcond2}. The final equality can be deduced from the observations that $f(\bx^*)=f^*$, $\lambda^*g(\bx^*)=0$, and $H(f^*)=0$. This indicates that $-\mu^{-1}$ acts as the subgradient of $H(r)$ at $r=f^*$. By referring to the definition of $\theta$ in \eqref{eq:conditionnumber}, we ascertain that $\theta \geq \mu^{-1}$. Consequently, a lower bound for $\theta$ can be established by bounding $\lambda^*$ and $|{\bar\gamma}^+|_1+|{\bar\gamma}^-|_1=|\gamma^*|_1$ from above. To achieve this,  we borrow some of the analysis from the proof of Lemma 3 in \citet{lin2022complexity}.


\noindent {\bf An upper bound on $\lambda^*$:} From the convexity of $g$ and the equality \eqref{eq:optcond1}, it follows that
{\small
\begin{eqnarray}
\lambda^*g(\tilde{\bx}) 
\geq\lambda^*\left[g(\bx^*) +(\tilde{\bx}-\bx^*)^\top\bzt_g\right]=-(\tilde{\bx}-\bx^*)^\top(\bzt_f +\bA^\top \gamma^*+ \bu^*)
=-(\tilde{\bx}-\bx^*)^\top(\bzt_f+ \bu^*),
\label{eq:boundlmabda1}
\end{eqnarray} }
where the first inequality and the last equality hold due to the non-negativity of $\lambda^*$ and $\bA\tilde{\bx}=\bA\bx^*=\bb$, respectively.

In the case of $u^* = 0$, inequality \eqref{eq:boundlmabda1} suggests that 
\begin{eqnarray}
\label{eq:boundlmabdau=0}
- \lambda^*g(\tilde{\bx}) 
\leq (\tilde{\bx} - \bx^*)^\top\bzt_f \leq f(\tilde{\bx}) - f^* \leq f(\tilde{\bx})-r_{\text{ini}}
\end{eqnarray}
where the second inequality arises from convexity of $f$ and the third from $r_{\text{ini}}\leq f^*$. The validity of this inequality also holds when $u^*\neq 0$, and we shortly proceed to demonstrate it. As such, $\lambda^*$ can be bounded above as follows 
\begin{eqnarray}
	\label{eq:boundlambdak}
\lambda^*\leq \dfrac{f(\tilde{\bx})-r_{\text{ini}}}{-g(\tilde{\bx}) }.
\end{eqnarray}

 To establish \eqref{eq:boundlmabdau=0} in the scenario of $u^*\neq 0$, consider that $\bx^*\in\partial\mathcal{\Xc}$. Let $\mathcal{H}$ represent the supporting hyperplane of $\Xc$ at $\bx^*$, i.e. $\big\{\bx\in\mathbb{R}^n\,|\,(\bx-\bx^*)^\top \bu^*=0\big\}$. Given that $\mathrm{dist}(\tilde{\bx},\mathcal{H}) = |(\bx^*-\tilde{\bx})^\top\bu^*|/\|\bu^*\|$ and $\mathrm{dist}(\tilde{\bx},\mathcal{H}) \ge \mathrm{dist}(\tilde{\bx},\partial \mathcal{X})>0$, we get
$
(\bx^*-\tilde{\bx})^\top\bu^*=|(\bx^*-\tilde{\bx})^\top\bu^*|=\mathrm{dist}(\tilde{\bx},\mathcal{H})\|\bu^*\| \ge \mathrm{dist}(\tilde{\bx},\partial \mathcal{X})\|\bu^* \|,
$	
where the first equality follows from 	$\bu^*\in\mathcal{N}_{\mathcal{X}}(\bx^*)$. Applying this inequality alongside \eqref{eq:boundlmabda1} yields
$
-\lambda^*g(\tilde{\bx})  + \mathrm{dist}(\tilde{\bx},\partial \mathcal{X})\|\bu^*\|  \le (\tilde{\bx}-\bx^*)^\top \bzt_f \leq f(\tilde{\bx})-f^*\leq f(\tilde{\bx})-r_{\text{ini}},
$
and subsequently,
\begin{eqnarray}
	\label{eq:boundlambdak2}
- \lambda^* g(\tilde{\bx})\leq f(\tilde{\bx})-r_{\text{ini}}\text{ and }\|\bu^*\| \leq \frac{f(\tilde{\bx})-r_{\text{ini}}}{\mathrm{dist}(\tilde{\bx},\partial \mathcal{X})}.
\end{eqnarray}
\noindent {\bf An upper bound on $\gamma^*$:} Utilizing equation \eqref{eq:optcond1} along with the full row rank of $\bA$, we deduce that 
$
\gamma^*= -(\bA\bA^\top)^{-1}\bA\left(\bzt_f +\lambda^* \bzt_g+ \bu^*\right),
$
which, together with \eqref{eq:boundlambdak}, leads to
\begin{align}
	\nonumber
	\|\gamma^*\|_1\leq\sqrt{l}\|\gamma^*\|
	\leq&\sqrt{l}\|(\bA\bA^\top)^{-1}\bA\|\left(M(1+\lambda^*)+\|\bu^*\|\right)\\\nonumber
	\leq&\sqrt{l}\|(\bA\bA^\top)^{-1}\bA\|\left(\frac{M(f(\tilde{\bx})-r_{\text{ini}}-g(\tilde{\bx}))}{-g(\tilde{\bx}) }+ \frac{f(\tilde{\bx})-r_{\text{ini}}}{\mathrm{dist}(\tilde{\bx},\partial \mathcal{X})}\right).
\end{align}
In the above inequalities, we incorporated \eqref{eq:boundlambdak2} and made use of Assumption~\ref{assumptions:nonsmooth}.
By applying this upper bound of $\|\gamma^*\|_1$ and the upper bound of $\lambda^*$ from \eqref{eq:boundlambdak} to $\mu^{-1}$, we successfully establish the desired lower bound for $\theta$. 
\hfill\BlackBox

\begin{remark}
	Note that Lemma~\ref{lemma:approximate_condition_number_linearequality} encompasses Lemma~\ref{lemma:approximate_condition_number} as a special case. In particular, in cases where linear equality constraints are absent, the introduction of $\gamma^*$ in \eqref{eq:optcond1} becomes unnecessary. Then the aforementioned lemma implies $\theta \geq (1+\lambda^*)^{-1}$. Given that \eqref{eq:boundlambdak} remains applicable, we still have $\theta \geq (1+\lambda^*)^{-1}=\frac{g(\tilde{\bx}) }{f^*-f(\tilde{\bx})+g(\tilde{\bx})}\geq \tilde \theta(r)$ for any $r<f^*$, where $\tilde\theta(r)$ is defined in Lemma~\ref{lemma:approximate_condition_number}. 
\end{remark}

\vskip 0.2in
\bibliography{references}

\begin{thebibliography}{68}
\providecommand{\natexlab}[1]{#1}
\providecommand{\url}[1]{\texttt{#1}}
\expandafter\ifx\csname urlstyle\endcsname\relax
  \providecommand{\doi}[1]{doi: #1}\else
  \providecommand{\doi}{doi: \begingroup \urlstyle{rm}\Url}\fi

\bibitem[Adcock et~al.(2025)Adcock, Colbrook, and
  Neyra-Nesterenko]{adcock2025restarts}
B.~Adcock, M.~J. Colbrook, and M.~Neyra-Nesterenko.
\newblock Restarts subject to approximate sharpness: a parameter-free and
  optimal scheme for first-order methods.
\newblock \emph{Foundations of Computational Mathematics}, pages 1--56, 2025.

\bibitem[Aravkin et~al.(2019)Aravkin, Burke, Drusvyatskiy, Friedlander, and
  Roy]{aravkin2019level}
A.~Aravkin, J.~Burke, D.~Drusvyatskiy, M.~Friedlander, and S.~Roy.
\newblock Level-set methods for convex optimization.
\newblock \emph{Mathematical Programming}, 174\penalty0 (1-2):\penalty0
  359--390, 2019.

\bibitem[Bayandina et~al.(2018)Bayandina, Dvurechensky, Gasnikov, Stonyakin,
  and Titov]{bayandina2018mirror}
A.~Bayandina, P.~Dvurechensky, A.~Gasnikov, F.~Stonyakin, and A.~Titov.
\newblock \emph{Mirror descent and convex optimization problems with non-smooth
  inequality constraints}, pages 181--213.
\newblock Springer, 2018.

\bibitem[Beck and Teboulle(2009)]{beck2009fast}
A.~Beck and M.~Teboulle.
\newblock A fast iterative shrinkage-thresholding algorithm for linear inverse
  problems.
\newblock \emph{SIAM Journal on Imaging Sciences}, 2\penalty0 (1):\penalty0
  183--202, 2009.

\bibitem[Bolte et~al.(2017)Bolte, Nguyen, Peypouquet, and
  Suter]{bolte2017error}
J.~Bolte, P.~Nguyen, J.~Peypouquet, and W.~Suter.
\newblock From error bounds to the complexity of first-order descent methods
  for convex functions.
\newblock \emph{Mathematical Programming}, 165\penalty0 (2):\penalty0 471--507,
  2017.

\bibitem[Boob et~al.(2023)Boob, Deng, and Lan]{boob2023stochastic}
D.~Boob, Q.~Deng, and G.~Lan.
\newblock Stochastic first-order methods for convex and nonconvex functional
  constrained optimization.
\newblock \emph{Mathematical Programming}, 197\penalty0 (1):\penalty0 215--279,
  2023.

\bibitem[Chambolle and Pock(2011)]{chambolle2011first}
A.~Chambolle and T.~Pock.
\newblock A first-order primal-dual algorithm for convex problems with
  applications to imaging.
\newblock \emph{Journal of Mathematical Imaging and Vision}, 40:\penalty0
  120--145, 2011.

\bibitem[Charisopoulos and Davis(2024)]{charisopoulos2022}
V.~Charisopoulos and D.~Davis.
\newblock A superlinearly convergent subgradient method for sharp semismooth
  problems.
\newblock \emph{Mathematics of Operations Research}, 49\penalty0 (3):\penalty0
  1678--1709, 2024.

\bibitem[Davis and Jiang(2024)]{davis2022}
D.~Davis and L.~Jiang.
\newblock A local nearly linearly convergent first-order method for nonsmooth
  functions with quadratic growth.
\newblock \emph{Foundations of Computational Mathematics}, pages 1--82, 2024.

\bibitem[Davis et~al.(2018)Davis, Drusvyatskiy, MacPhee, and
  Paquette]{davis2018}
D.~Davis, D.~Drusvyatskiy, K.~J. MacPhee, and C.~Paquette.
\newblock Subgradient methods for sharp weakly convex functions.
\newblock \emph{Journal of Optimization Theory and Applications}, 179\penalty0
  (3):\penalty0 962--982, 2018.

\bibitem[Davis et~al.(2024)Davis, Drusvyatskiy, and
  Charisopoulos]{davis2019stochastic}
D.~Davis, D.~Drusvyatskiy, and V.~Charisopoulos.
\newblock Stochastic algorithms with geometric step decay converge linearly on
  sharp functions.
\newblock \emph{Mathematical Programming}, 207\penalty0 (1):\penalty0 145--190,
  2024.

\bibitem[Deng et~al.(2024)Deng, Lan, and Lin]{deng2024}
Q.~Deng, G.~Lan, and Z.~Lin.
\newblock Uniformly optimal and parameter-free first-order methods for convex
  and function-constrained optimization.
\newblock \emph{arXiv preprint arXiv:2412.06319}, 2024.

\bibitem[D{\'\i}az and Grimmer(2023)]{diaz2021}
M.~D{\'\i}az and B.~Grimmer.
\newblock Optimal convergence rates for the proximal bundle method.
\newblock \emph{SIAM Journal on Optimization}, 33\penalty0 (2):\penalty0
  424--454, 2023.

\bibitem[Drusvyatskiy and Paquette(2019)]{CompositionsofConvexandSmooth}
D.~Drusvyatskiy and C.~Paquette.
\newblock Efficiency of minimizing compositions of convex functions and smooth
  maps.
\newblock \emph{Mathematical Programming}, 178:\penalty0 503--558, 2019.

\bibitem[Ene et~al.(2021)Ene, Nguyen, and Vladu]{ene2021}
A.~Ene, H.~L. Nguyen, and A.~Vladu.
\newblock Adaptive gradient methods for constrained convex optimization and
  variational inequalities.
\newblock In \emph{Proceedings of the AAAI Conference on Artificial
  Intelligence}, volume~35, pages 7314--7321, 2021.

\bibitem[Estrin and Friedlander(2020)]{Friedlander2020}
R.~Estrin and M.~P. Friedlander.
\newblock A perturbation view of level-set methods for convex optimization.
\newblock \emph{Optimization Letters}, 14\penalty0 (8):\penalty0 1989--2006,
  2020.

\bibitem[Fercoq(2023)]{fercoq2022}
O.~Fercoq.
\newblock Quadratic error bound of the smoothed gap and the restarted averaged
  primal-dual hybrid gradient.
\newblock \emph{Open Journal of Mathematical Optimization}, 4:\penalty0 1--34,
  2023.

\bibitem[Fercoq and Qu(2019)]{fercoq2019adaptive}
O.~Fercoq and Z.~Qu.
\newblock Adaptive restart of accelerated gradient methods under local
  quadratic growth condition.
\newblock \emph{IMA Journal of Numerical Analysis}, 39\penalty0 (4):\penalty0
  2069--2095, 2019.

\bibitem[Fercoq and Qu(2020)]{fercoq2020restarting}
O.~Fercoq and Z.~Qu.
\newblock Restarting the accelerated coordinate descent method with a rough
  strong convexity estimate.
\newblock \emph{Computational Optimization and Applications}, 75\penalty0
  (1):\penalty0 63--91, 2020.

\bibitem[Fercoq et~al.(2019)Fercoq, Alacaoglu, Necoara, and
  Cevher]{fercoq2019almost}
O.~Fercoq, A.~Alacaoglu, I.~Necoara, and V.~Cevher.
\newblock Almost surely constrained convex optimization.
\newblock In \emph{International Conference on Machine Learning}, pages
  1910--1919. PMLR, 2019.

\bibitem[Freund and Lu(2018)]{freund2018new}
R.~Freund and H.~Lu.
\newblock New computational guarantees for solving convex optimization problems
  with first order methods, via a function growth condition measure.
\newblock \emph{Mathematical Programming}, 170\penalty0 (2):\penalty0 445--477,
  2018.

\bibitem[Goh et~al.(2016)Goh, Cotter, Gupta, and
  Friedlander]{goh2016satisfying}
G.~Goh, A.~Cotter, M.~Gupta, and M.~P. Friedlander.
\newblock Satisfying real-world goals with dataset constraints.
\newblock \emph{Advances in Neural Information Processing Systems}, 29, 2016.

\bibitem[Grimmer(2018)]{grimmer2018radial}
B.~Grimmer.
\newblock Radial subgradient method.
\newblock \emph{SIAM Journal on Optimization}, 28\penalty0 (1):\penalty0
  459--469, 2018.

\bibitem[Grimmer(2019)]{grimmer2019convergence}
B.~Grimmer.
\newblock Convergence rates for deterministic and stochastic subgradient
  methods without lipschitz continuity.
\newblock \emph{SIAM Journal on Optimization}, 29\penalty0 (2):\penalty0
  1350--1365, 2019.

\bibitem[Grimmer(2023)]{grimmer2021general}
B.~Grimmer.
\newblock General h{\"o}lder smooth convergence rates follow from specialized
  rates assuming growth bounds.
\newblock \emph{Journal of Optimization Theory and Applications}, 197\penalty0
  (1):\penalty0 51--70, 2023.

\bibitem[Grimmer(2024)]{grimmer2024optimal}
B.~Grimmer.
\newblock On optimal universal first-order methods for minimizing heterogeneous
  sums.
\newblock \emph{Optimization Letters}, 18\penalty0 (2):\penalty0 427--445,
  2024.

\bibitem[Guyon et~al.(2004)Guyon, Gunn, Ben-Hur, and Dror]{guyon2004result}
I.~Guyon, S.~Gunn, A.~Ben-Hur, and G.~Dror.
\newblock Result analysis of the nips 2003 feature selection challenge.
\newblock \emph{Advances in neural information processing systems}, 17, 2004.

\bibitem[Hamedani and Aybat(2021)]{hamedani2021primal}
E.~Y. Hamedani and N.~S. Aybat.
\newblock A primal-dual algorithm with line search for general convex-concave
  saddle point problems.
\newblock \emph{SIAM Journal on Optimization}, 31\penalty0 (2):\penalty0
  1299--1329, 2021.

\bibitem[Iouditski and Nesterov(2014)]{iouditski2014primal}
A.~Iouditski and Y.~Nesterov.
\newblock Primal-dual subgradient methods for minimizing uniformly convex
  functions.
\newblock \emph{arXiv preprint arXiv:1401.1792}, 2014.

\bibitem[Ito and Fukuda(2021)]{ito2021}
M.~Ito and M.~Fukuda.
\newblock Nearly optimal first-order methods for convex optimization under
  gradient norm measure: An adaptive regularization approach.
\newblock \emph{Journal of Optimization Theory and Applications}, 188\penalty0
  (3):\penalty0 770--804, 2021.

\bibitem[Ito et~al.(2023)Ito, Lu, and He]{ito2023}
M.~Ito, Z.~Lu, and C.~He.
\newblock A parameter-free conditional gradient method for composite
  minimization under h{\"o}lder condition.
\newblock \emph{Journal of Machine Learning Research}, 24\penalty0
  (166):\penalty0 1--34, 2023.

\bibitem[Johnstone and Moulin(2020)]{johnstone2020faster}
P.~R. Johnstone and P.~Moulin.
\newblock Faster subgradient methods for functions with h{\"o}lderian growth.
\newblock \emph{Mathematical Programming}, 180\penalty0 (1):\penalty0 417--450,
  2020.

\bibitem[Lan and Zhou(2016)]{lan2016algorithms}
G.~Lan and Z.~Zhou.
\newblock Algorithms for stochastic optimization with expectation constraints.
\newblock \emph{arXiv preprint arXiv:1604.03887}, 2016.

\bibitem[Lan et~al.(2023)Lan, Ouyang, and Zhang]{lan2023optimal}
G.~Lan, Y.~Ouyang, and Z.~Zhang.
\newblock Optimal and parameter-free gradient minimization methods for convex
  and nonconvex optimization.
\newblock \emph{arXiv preprint arXiv:2310.12139}, 2023.

\bibitem[Lewis et~al.(2004)Lewis, Yang, Rose, and Li]{lewis2004rcv1}
D.~D. Lewis, Y.~Yang, T.~G. Rose, and F.~Li.
\newblock Rcv1: A new benchmark collection for text categorization research.
\newblock \emph{Journal of Machine Learning Research}, 5\penalty0
  (Apr):\penalty0 361--397, 2004.

\bibitem[Li et~al.(2015)Li, Mordukhovich, and Pham]{li2015new}
G.~Li, B.~S. Mordukhovich, and T.~Pham.
\newblock New fractional error bounds for polynomial systems with applications
  to h{\"o}lderian stability in optimization and spectral theory of tensors.
\newblock \emph{Mathematical Programming}, 153\penalty0 (2):\penalty0 333--362,
  2015.

\bibitem[Lin and Xiao(2015)]{lin2015adaptive}
Q.~Lin and L.~Xiao.
\newblock An adaptive accelerated proximal gradient method and its homotopy
  continuation for sparse optimization.
\newblock \emph{Computational Optimization and Applications}, 60\penalty0
  (3):\penalty0 633--674, 2015.

\bibitem[Lin et~al.(2018{\natexlab{a}})Lin, Ma, and Yang]{lin2018level}
Q.~Lin, R.~Ma, and T.~Yang.
\newblock Level-set methods for finite-sum constrained convex optimization.
\newblock In \emph{International Conference on Machine Learning}, pages
  3112--3121, 2018{\natexlab{a}}.

\bibitem[Lin et~al.(2018{\natexlab{b}})Lin, Nadarajah, and
  Soheili]{lin2018feasiblelevel}
Q.~Lin, S.~Nadarajah, and N.~Soheili.
\newblock A level-set method for convex optimization with a feasible solution
  path.
\newblock \emph{SIAM Journal on Optimization}, 28\penalty0 (4):\penalty0
  3290--3311, 2018{\natexlab{b}}.

\bibitem[Lin et~al.(2020)Lin, Nadarajah, Soheili, and Yang]{lin2019data}
Q.~Lin, S.~Nadarajah, N.~Soheili, and T.~Yang.
\newblock A data efficient and feasible level set method for stochastic convex
  optimization with expectation constraints.
\newblock \emph{Journal of Machine Learning Research}, 2020.

\bibitem[Lin et~al.(2022)Lin, Ma, and Xu]{lin2022complexity}
Q.~Lin, R.~Ma, and Y.~Xu.
\newblock Complexity of an inexact proximal-point penalty method for
  constrained smooth non-convex optimization.
\newblock \emph{Computational Optimization and Applications}, 82\penalty0
  (1):\penalty0 175--224, 2022.

\bibitem[Liu and Yang(2017)]{liu2017adaptive}
M.~Liu and T.~Yang.
\newblock Adaptive accelerated gradient converging method under h{\"o}lderian
  error bound condition.
\newblock In \emph{Advances in Neural Information Processing Systems}, pages
  3104--3114, 2017.

\bibitem[Necoara et~al.(2019)Necoara, Nesterov, and Glineur]{necoara2019linear}
I.~Necoara, Y.~Nesterov, and F.~Glineur.
\newblock Linear convergence of first order methods for non-strongly convex
  optimization.
\newblock \emph{Mathematical Programming}, 175\penalty0 (1-2):\penalty0
  69--107, 2019.

\bibitem[Nesterov(2013)]{nesterov_nonsmoothminimization}
Y.~Nesterov.
\newblock Gradient methods for minimizing composite functions.
\newblock \emph{Mathematical Programming}, 140\penalty0 (1):\penalty0 125--161,
  2013.

\bibitem[Nesterov(2018)]{nesterov2018lectures}
Y.~Nesterov.
\newblock \emph{Lectures on Convex Optimization}, volume 137.
\newblock Springer, 2018.

\bibitem[Polyak(1969)]{polyak1969}
B.~T. Polyak.
\newblock Minimization of unsmooth functionals.
\newblock \emph{USSR Computational Mathematics and Mathematical Physics},
  9\penalty0 (3):\penalty0 14--29, 1969.

\bibitem[Renegar(2016)]{renegar2016efficient}
J.~Renegar.
\newblock ``efficient” subgradient methods for general convex optimization.
\newblock \emph{SIAM Journal on Optimization}, 26\penalty0 (4):\penalty0
  2649--2676, 2016.

\bibitem[Renegar and Grimmer(2022)]{renegar2018simple}
J.~Renegar and B.~Grimmer.
\newblock A simple nearly optimal restart scheme for speeding up first-order
  methods.
\newblock \emph{Foundations of Computational Mathematics}, 22\penalty0
  (1):\penalty0 211--256, 2022.

\bibitem[Renegar and Zhou(2021)]{renegar2021FeasProb}
J.~Renegar and S.~Zhou.
\newblock A different perspective on the stochastic convex feasibility problem.
\newblock \emph{arXiv preprint arXiv:2108.12029}, 2021.

\bibitem[Roulet and d'Aspremont(2017)]{roulet2017sharpness}
V.~Roulet and A.~d'Aspremont.
\newblock Sharpness, restart and acceleration.
\newblock In \emph{Advances in Neural Information Processing Systems}, pages
  1119--1129, 2017.

\bibitem[Sujanani and Monteiro(2025)]{sujanani2024}
A.~Sujanani and R.~D. Monteiro.
\newblock Efficient parameter-free restarted accelerated gradient methods for
  convex and strongly convex optimization: A. sujanani, rdc monteiro.
\newblock \emph{Journal of Optimization Theory and Applications}, 206\penalty0
  (2):\penalty0 52, 2025.

\bibitem[Wang and Liu(2022)]{wang2022}
T.~Wang and H.~Liu.
\newblock Convergence results of a new monotone inertial forward--backward
  splitting algorithm under the local h{\"o}lder error bound condition.
\newblock \emph{Applied Mathematics \& Optimization}, 85\penalty0 (2):\penalty0
  7, 2022.

\bibitem[Wei et~al.(2018)Wei, Yu, Ling, and Neely]{wei2018solving}
X.~Wei, H.~Yu, Q.~Ling, and M.~Neely.
\newblock Solving non-smooth constrained programs with lower complexity than
  {$O(1/\epsilon)$}: A primal-dual homotopy smoothing approach.
\newblock In \emph{Advances in Neural Information Processing Systems}, pages
  3995--4005, 2018.

\bibitem[Wolfe(1959)]{wolfe1959simplex}
P.~Wolfe.
\newblock The simplex method for quadratic programming.
\newblock \emph{Econometrica: Journal of the Econometric Society}, pages
  382--398, 1959.

\bibitem[Xu(2020)]{xu2018primal}
Y.~Xu.
\newblock Primal-dual stochastic gradient method for convex programs with many
  functional constraints.
\newblock \emph{SIAM Journal on Optimization}, 30\penalty0 (2):\penalty0
  1664--1692, 2020.

\bibitem[Xu(2021{\natexlab{a}})]{xu2017first}
Y.~Xu.
\newblock First-order methods for constrained convex programming based on
  linearized augmented lagrangian function.
\newblock \emph{INFORMS Journal on Optimization}, 3\penalty0 (1):\penalty0
  89--117, 2021{\natexlab{a}}.

\bibitem[Xu(2021{\natexlab{b}})]{xu2017global}
Y.~Xu.
\newblock Iteration complexity of inexact augmented lagrangian methods for
  constrained convex programming.
\newblock \emph{Mathematical Programming}, 185\penalty0 (1):\penalty0 199--244,
  2021{\natexlab{b}}.

\bibitem[Xu and Yang(2018)]{xu2018frank}
Y.~Xu and T.~Yang.
\newblock Frank-wolfe method is automatically adaptive to error bound
  condition.
\newblock \emph{arXiv preprint arXiv:1810.04765}, 2018.

\bibitem[Xu et~al.(2016)Xu, Yan, Lin, and Yang]{xu2016homotopy}
Y.~Xu, Y.~Yan, Q.~Lin, and T.~Yang.
\newblock Homotopy smoothing for non-smooth problems with lower complexity than
  $\mathcal{O}(1/\epsilon)$.
\newblock In \emph{Advances In Neural Information Processing Systems}, pages
  1208--1216, 2016.

\bibitem[Xu et~al.(2017{\natexlab{a}})Xu, Lin, and Yang]{xu2017adaptive}
Y.~Xu, Q.~Lin, and T.~Yang.
\newblock Adaptive svrg methods under error bound conditions with unknown
  growth parameter.
\newblock In \emph{Advances in Neural Information Processing Systems}, pages
  3277--3287, 2017{\natexlab{a}}.

\bibitem[Xu et~al.(2017{\natexlab{b}})Xu, Lin, and Yang]{xu2017stochastic}
Y.~Xu, Q.~Lin, and T.~Yang.
\newblock Stochastic convex optimization: Faster local growth implies faster
  global convergence.
\newblock In \emph{International Conference on Machine Learning}, pages
  3821--3830, 2017{\natexlab{b}}.

\bibitem[Xu et~al.(2017{\natexlab{c}})Xu, Liu, Lin, and Yang]{xu2017admm}
Y.~Xu, M.~Liu, Q.~Lin, and T.~Yang.
\newblock Admm without a fixed penalty parameter: Faster convergence with new
  adaptive penalization.
\newblock In \emph{Advances in Neural Information Processing Systems}, pages
  1267--1277, 2017{\natexlab{c}}.

\bibitem[Yan et~al.(2019)Yan, Xu, Lin, Zhang, and Yang]{yan2019stochastic}
Y.~Yan, Y.~Xu, Q.~Lin, L.~Zhang, and T.~Yang.
\newblock Stochastic primal-dual algorithms with faster convergence than
  $\mathcal{O}(1/\sqrt{T})$ for problems without bilinear structure.
\newblock \emph{arXiv preprint arXiv:1904.10112}, 2019.

\bibitem[Yang and Lin(2018)]{yang2018rsg}
T.~Yang and Q.~Lin.
\newblock Rsg: Beating subgradient method without smoothness and strong
  convexity.
\newblock \emph{The Journal of Machine Learning Research}, 19\penalty0
  (1):\penalty0 236--268, 2018.

\bibitem[Yang et~al.(2017)Yang, Lin, and Zhang]{yang2017richer}
T.~Yang, Q.~Lin, and L.~Zhang.
\newblock A richer theory of convex constrained optimization with reduced
  projections and improved rates.
\newblock In \emph{Proceedings of the 34th International Conference on Machine
  Learning-Volume 70}, pages 3901--3910. JMLR. org, 2017.

\bibitem[Yu et~al.(2017)Yu, Neely, and Wei]{yu2017online}
H.~Yu, M.~Neely, and X.~Wei.
\newblock Online convex optimization with stochastic constraints.
\newblock \emph{Advances in Neural Information Processing Systems}, 30, 2017.

\bibitem[Zhang(2020)]{zhang2020new}
H.~Zhang.
\newblock New analysis of linear convergence of gradient-type methods via
  unifying error bound conditions.
\newblock \emph{Mathematical Programming}, 180\penalty0 (1):\penalty0 371--416,
  2020.

\bibitem[Zhang et~al.(2021)Zhang, Dai, Guo, and Peng]{zhang2019proximal}
H.~Zhang, Y.-H. Dai, L.~Guo, and W.~Peng.
\newblock Proximal-like incremental aggregated gradient method with linear
  convergence under bregman distance growth conditions.
\newblock \emph{Mathematics of Operations Research}, 46\penalty0 (1):\penalty0
  61--81, 2021.

\end{thebibliography}

\end{document}